\documentclass[a4paper, 12pt]{article}

\usepackage{amssymb, amsmath, amscd, amsthm, indentfirst, multirow, nccmath, longtable}
\usepackage[inline]{enumitem} 
\usepackage[all]{xy}

\newtheorem{thm}{Theorem}[section] 
\newtheorem{lem}[thm]{Lemma}

\newtheorem{cor}[thm]{Corollary}

\newcommand{\QED}{{\unskip\nobreak\hfil\penalty50\quad\null\nobreak\hfil
{Q.E.D.}\parfillskip0pt\finalhyphendemerits0\par\medskip}}
\newcommand{\Proof}{\noindent{\bf Proof.}\quad}

\newcommand{\G}{\mathbb{G}}
\newcommand{\Aut}{{\rm Aut}}
\newcommand{\SL}{{\rm SL}}
\newcommand{\GL}{{\rm GL}}
\newcommand{\mf}{\mathfrak}
\newcommand{\Mat}{{\rm Mat}}
\newcommand{\diag}{{\rm diag}}
\newcommand{\cU}{\mathcal{U}}
\newcommand{\cA}{\mathcal{A}}
\newcommand{\Z}{\mathbb{Z}}
\newcommand{\ds}{\displaystyle}
\newcommand{\id}{{\rm id}}
\newcommand{\ol}{\overline}

\allowdisplaybreaks[4]

\begin{document}

\title{\Large Representations of $\G_a \rtimes \G_m$ into $\SL(3, k)$ in positive characteristic}
\author{Ryuji Tanimoto} 
\date{}
\maketitle

\begin{abstract}
Let $k$ be an algebraically closed field of positive characteristic $p$. 
In this article, we classify representations of $\G_a \rtimes \G_m$ into $\SL(3, k)$, 
and thereby we classify fundamental representations of $\G_a$ into $\SL(3, k)$. 
\end{abstract}


\section{Introduction}

Let $k$ be an algebraically closed field, 
let $\G_a$ denote the additive group of $k$ and 
let $\G_m$ denote the multiplicative group of $k$. 
The special linear group $\SL(2, k)$ has the ascending sequence 
\[
 \hspace{-6pt} \xymatrix@C=18pt@M=12pt{ \G_a \ar@<-0.5ex>@{^(->}[r]^{\;\ \raisebox{6pt}{$\iota^+$} } & B}\hspace{-6pt}\quad   \subset \quad   \SL(2, k) 
\]
of subgroups, where $B$ is the Borel subgroup of $\SL(2, k)$ defined by 
\[
 B := 
\left\{ 
\left. 
\left(
\begin{array}{c c}
 a & b \\
 c & d
\end{array}
\right) 
 \in 
\SL(2, k) 
\; \right| \; 
c = 0 \, 
\right\} 
\] 
and $\hspace{-6pt} \xymatrix@C=12pt@M=6pt{\iota^+ : \G_a \ar@{^(->}[r] & B}\hspace{-6pt} $ is the homomorphism defined by 
\[
 \iota^+(t) 
 := 
\left(
\begin{array}{c c}
 1 & t \\
 0 & 1 
\end{array}
\right) . 
\]
Clearly, $B$ is isomorphic to the semidirect product $\G_a \rtimes \G_m$ of $\G_a$ by $\G_m$.

Given a representation of $\SL(2, k)$ into $\SL(n, k)$, we naturally have a representation of $B$ into $\SL(n, k)$ and 
also have a representation of $\G_a$ into $\SL(n, k)$. Conversely, given a representation $u$ of $\G_a$ into $\SL(n, k)$, 
we ask whether or not $u$ can be extended to a representation $\varphi : B \to \SL(n, k)$; and if such a $\varphi$ exists,  
we further ask whether or not $\varphi$ can be extended to a representation $\psi : \SL(2, k) \to \SL(n, k)$. 
If the characteristic of $k$ is zero, any representation $u : \G_a \to \SL(n, k)$ can be extended until a representation 
$\psi : \SL(2, k) \to \SL(n, k)$. 
However, if the characteristic of $k$ is greater than or equal to three, Fauntleroy \cite{Fauntleroy} found a three-dimensional 
representation $u : \G_a \to \SL(3, k)$ such that $u$ cannot be extended to any representation $\psi : \SL(2, k) \to \SL(3, k)$. 
This interesting representation implies that even if we happen to classify representations of $\SL(2, k)$ into $\SL(n, k)$ ($n \geq 3$), 
we cannot know all representations of $\G_a$ into $\SL(n, k)$ in positive characteristic. 
It is a complicated problem to classify representations of $\G_a$ into $\SL(n, k)$ in positive characteristic. 
At present, for $n = 1, 2, 3, 4, 5$, we can give overlapping classifications of $\G_a$ into $\SL(n, k)$ in positive characteristic (cf. \cite{Tanimoto 2020, Tanimoto 2019, Tanimoto 2008}). 
Besides the problem, 
it remains a problem to classify representations of $\G_a \rtimes \G_m$ into $\SL(n, k)$ in positive characteristic.

In this article, we classify representations of $\G_a \rtimes \G_m$ into $\SL(3, k)$ in positive characteristic $p$ (see Theorem 2.3 
and Corollary 2.4), and thereby 
classify representations $u : \G_a \to \SL(3, k)$ so that $u$ can be extended to representations $\psi : \SL(2, k) \to \SL(3, k)$ 
(see Theorem 4.1 and Corollary 4.2). 
As a result, for each $p \geq 2$, we can find new three-dimensional representations $u : \G_a \to \SL(3, k)$ so that $u$ 
cannot be extended to any representation $\psi : \SL(2, k) \to \SL(3, k)$ (see Corollary 4.3). 
\\

\noindent 
{\bf Notations and definitions}  
\medskip

From now on until the last section of this article, we assume that the characteristic $p$ of $k$ is positive.  
Let $k[T]$ be the polynomial ring in one variable over $k$.  
We say that a polynomial $f(T)$ of $k[T]$ is a {\it $p$-polynomial} if $f(T)$ can be written in the form 
\[
 f(T) = \sum_{i \geq 0} a_i \, T^{p^i} \qquad (\, a_i \in k \quad \text{ for all \, $i \geq 0$}\, ) . 
\]
Let $\mf{P}$ denote the set of all $p$-polynomials.

In this article, we consider algebraic groups over $k$ (in particular, $\G_m$, $\G_a$, $\G_a \rtimes \G_m$, $\SL(2, k)$) 
and their representations as algebraic group. 
We simply write $\G_a \rtimes \G_m$ instead of $\G_a \rtimes_\sigma \G_m$, where the homomorphism 
$\sigma : \G_m \to \Aut_k(\G_a)$, $z \mapsto \sigma_z$ is defined by $\sigma_z(t) := z^2 \, t$ for all $t \in \G_a$ 
and the product of elements $(t_1, z_1)$, $(t_2, z_2)$ of $\G_a \rtimes \G_m$ is defined by 
\[
 (t_1, z_1) \cdot (t_2, z_2) := (\, t_1 + \sigma_{z_1} (t_2), \, z_1 \, z_2 \, ) . 
\]
We can identify the semi-direct product $\G_a \rtimes \G_m$ with the Borel subgroup $B$ of $\SL(2, k)$
consisting of all upper triangular matrices of $\SL(2, k)$. 
In fact, consider the isomorphism from $\G_a \rtimes \G_m$ to $B$ defined by 
\[
 (t, \, z) 
 \mapsto 
\left(
\begin{array}{c c}
 z & t \, z^{- 1} \\
 0 & z^{- 1}
\end{array}
\right) 
\left( 
 =
\left(
\begin{array}{c c}
 1 & t \\ 
 0 & 1 
\end{array}
\right) 
\left(
\begin{array}{c c}
 z & 0 \\ 
 0 & z^{- 1} 
\end{array}
\right) 
\right) . 
\]
For an algebraic group $G$, 
we say that two representations $\rho_1 : G \to \GL(n, k)$ and $\rho_2 : G \to \GL(n, k)$ are {\it equivalent} if 
there exists a regular matrix $P$ of $\GL(n, k)$ such that $P^{-1} \, \rho_1(g) \, P = \rho_2(g)$ for all $g \in G$.  
We also say that two representations $\rho_1 : G \to \SL(n, k)$ and $\rho_2 : G \to \SL(n, k)$ of $G$ into $\SL(n, k)$ 
are {\it equivalent} if there exists a regular matrix $P$ of $\GL(n, k)$ such that $P^{-1} \, \rho_1(g) \, P = \rho_2(g)$ for all $g \in G$.

Let $R$ be a commutative ring. 
For $n, n' \geq 1$, we denote by $\Mat_{n,\, n'}(R)$ the set of all $n \times n'$ matrices whose all entries belong to $R$. 
In particular when $n = n'$, we write $\Mat(n, R)$ in place of $\Mat_{n, \, n'}(R)$. 
For a matrix $A$ of $\Mat(n, R)$, we denote by $^\top \! A$ the transpose of $A$. 
For elements $d_i$ $(1\leq i \leq n)$ of $R$, we denote by $\diag(d_1, \ldots, d_n)$ the diagonal matrix 
of $\Mat(n, R)$ whose $(i, i)$-th entries are $d_i$ for all $1 \leq i \leq n$, i.e., 
\[
\diag(d_1, \ldots, d_n) 
 := 
\left(
\begin{array}{c c c c c}
 d_1 & 0 & \cdots & \cdots &  0  \\
 0 & d_2 & 0  &  & \vdots  \\
 \vdots & 0 & \ddots & \ddots & \vdots \\
 \vdots & &  \ddots  & \ddots & 0 \\
 0 & \cdots & \cdots & 0  & d_n  
\end{array} 
\right) . 
\]

\section{Preliminaries}

\subsection{Representations of $\G_m$}

\begin{lem} 
Let $h : \G_m \to \GL(n, k)$ be a representation of $\G_m$. Then there exists a regular matrix $P$ of $\GL(n, k)$ such 
that for all $z \in \G_m$,  
\begin{align*}
 P^{-1} \, h(z) \, P
 = 
\diag(\, z^{\ell_1}, \, z^{\ell_2}, \, \ldots, \, z^{\ell_n} \, ) \qquad  (\, \ell_1 \geq \ell_2 \geq \cdots \geq \ell_n \,) . 
\end{align*}
\end{lem}

\Proof 
The proof is straightforward. 
\QED 

%

\begin{lem} 
Let $h : \G_m \to \GL(n, k)$ and $h' : \G_m \to \GL(n, k)$ be representations of $\G_m$ with the forms 
\begin{align*}
\left\{
\begin{array}{l @{\qquad} l}
  h(z) = \diag( \, z^{\ell_1}, \, z^{\ell_2}, \, \ldots, \, z^{\ell_n} \,)  
 &  (\, \ell_1 \geq \ell_2 \geq \cdots \geq \ell_n \,) , \\ 
  h'(z) = \diag( \, z^{\ell_1'}, \, z^{\ell_2'}, \, \ldots, \, z^{\ell_n'} \,)  
 &  (\, \ell_1' \geq \ell_2' \geq \cdots \geq \ell_n' \,) . 
\end{array}
\right. 
\end{align*} 
Assume that $h$ is equivalent to $h'$. 
Then we have 
\[
 (\, \ell_1 , \, \ell_2 , \, \ldots, \, \ell_n \,)
 = 
 (\, \ell_1' , \, \ell_2' , \, \ldots, \,  \ell_n' \,) . 
\]
\end{lem}

\Proof 
Let $V := k^{\oplus n}$ be the column space of dimension $n$. 
For any integer $\ell$, we let $W_\ell$ and ${W'}_{\ell}$ be the subspaces of $V$ defined by 
\begin{align*}
 W_\ell &:= 
\{
\,
 v \in V 
 \mid 
  h(z) \, v = z^\ell \, v \quad \text{ for all \, $z \in \G_m$}
\,
\}  , \\
{W'}_\ell &:= 
\{
\,
 v \in V  
 \mid 
  h'(z) \, v = z^\ell \, v \quad \text{ for all \, $z \in \G_m$} 
\,
\} . 
\end{align*} 
Since $h$ is equivalent to $h'$, there exists a regular matrix $P$ of $\GL(n, k)$ such that 
$P^{-1} h(z) P = h'(z)$ for all $z \in \G_m$. 
For any integer $\ell$ and for any $v \in V$, we have 
\begin{align*}
 v \in W_\ell 
 & \quad \Longleftrightarrow  \quad 
 h(z) \, v = z^{\ell} \, v \quad \text{ for all \; $z \in \G_m$} \\
 & \quad \Longleftrightarrow  \quad 
 P h'(z) P^{-1} \, v = z^{\ell} \, v \quad \text{ for all \; $z \in \G_m$} \\
 & \quad \Longleftrightarrow  \quad 
 h'(z) P^{-1} \, v = z^{\ell} \, P^{-1} \, v \quad \text{ for all \; $z \in \G_m$} \\
 & \quad \Longleftrightarrow  \quad 
 P^{-1} \, v \in {W'}_\ell .
\end{align*}
Thus $W_\ell$ is isomorphic to ${W'}_\ell$. 
Thus we can obtain the desired equality. 
\QED

\subsection{Representations of $\G_a$ into $\SL(3, k)$}
 
In this Subsection 1.1, we give a classification of representations of $\G_a$ into $\SL(3, k)$. 
For stating this classification, we prepare the following six sets 
$\cU_{[3]}$, $\cU_{[3]}^-$, 
$\cA_{(1, 2)}$, $\cA_{(2, 1)}$, 
$\cA_{(1, 2)}^-$, $\cA_{(2, 1)}^-$.

For $p \geq 3$, 
we denote by $\cU_{[3]}$ the set of all morphisms $u : \G_a \to \SL(3, k)$ of affine $k$-varieties with the form 
\begin{align*}
u(t) 
 =
\left(
\begin{array}{c c c}
 1 & \alpha_1(t) & \frac{\lambda}{2} \, \alpha_1(t)^2 + \alpha_2(t) \\
 0 & 1 & \lambda \, \alpha_1(t) \\
 0 & 0 & 1 
\end{array}
\right)  \qquad 
\left( 
\begin{array}{ll}
 \alpha_1(T) \in \mf{P} \backslash \{ 0 \} ,  &  \alpha_2(T) \in \mf{P} \\
 \lambda \in k \backslash \{ 0 \} & 
\end{array}
\right) .  
\end{align*}
Clearly, any element of $\cU_{[3]}$ is a representation of $\G_a$.

For $p \geq 3$, we denote by $\cU_{[3]}^-$ the set defined by 
\[
 \cU_{[3]}^- 
 := 
\left\{
 u^- : \G_a \to \SL(3, k) \; \left| \;
\begin{array}{l}
\text{there exists $u$ of $\cU_{[3]}$ such that} \\
\text{$u^-(t) = {^\top} \! u(t)$ for all $t \in \G_a$} 
\end{array} 
\right. 
\right\} . 
\]
Clearly, any element of $\cU_{[3]}^-$ is a representation of $\G_a$.

For $p \geq 2$, we denote by $\cA_{(1, 2)}$ the set of all morphisms $u : \G_a \to \SL(3, k)$ of affine $k$-varieties with the form 
\begin{align*}
u(t) 
 =
\left(
\begin{array}{c | c c}
 1 & \alpha_1(t) & \alpha_2(t) \\ 
\hline 
 0 & 1 & 0 \\
 0 & 0 & 1 
\end{array}
\right) 
\qquad 
\bigl( \, \alpha_1(T), \alpha_2(T) \in \mf{P} \, \bigr) . 
\end{align*}
Clearly, any element of $\cA_{(1, 2)}$ is a representation of $\G_a$.

For $p \geq 2$, we denote by $\cA_{(2, 1)}$ the set of all morphisms $u : \G_a \to \SL(3, k)$ of affine $k$-varieties with the form 
\begin{align*}
u(t) 
 =
\left(
\begin{array}{c c | c}
 1 & 0 & \alpha_2(t) \\
 0 & 1 & \alpha_1(t) \\ 
\hline 
 0 & 0 & 1 
\end{array}
\right)  
\qquad 
\bigl( \, \alpha_1(T), \alpha_2(T) \in \mf{P} \, \bigr) . 
\end{align*}
Clearly, any element of $\cA_{(2, 1)}$ is a representation of $\G_a$.

For $p \geq 2$, we denote by $\cA_{(1, 2)}^-$ and $\cA_{(2, 1)}^-$ the sets defined by 
\begin{align*}
 \cA_{(1, 2)}^-
 & := 
\left\{
 u^- : \G_a \to \SL(3, k) \; \left| \;
\begin{array}{l}
\text{there exists $u$ of $\cA_{(2, 1)}$ such that} \\
\text{$u^-(t) = {^\top} \! u(t)$ for all $t \in \G_a$}
\end{array} 
\right. 
\right\}  \\
\intertext{and} 
 \cA_{(2, 1)}^-
 & := 
\left\{
 u^- : \G_a \to \SL(3, k) \; \left| \;
\begin{array}{l}
\text{there exists $u$ of $\cA_{(1, 2)}$ such that} \\
\text{$u^-(t) = {^\top} \! u(t)$ for all $t \in \G_a$}
\end{array} 
\right. 
\right\} . 
\end{align*}
Clearly, any element of $\cA_{(1, 2)}^-$ and $\cA_{(2, 1)}^-$ is a representation of $\G_a$.

\begin{lem}
Let $u : \G_a \to \SL(n, k)$ be a representation. 
Then the following assertions {\rm (1)} and {\rm (2)} hold true: 
\begin{enumerate}[label = {\rm (\arabic*)}]
\item There exists a representation $u^* : \G_a \to \SL(n, k)$ such that $u^*$ is equivalent to $u$ and 
$u^*(t)$ is an upper triangular matrix for each $t \in \G_a$. 

\item There exists a representation $u_* : \G_a \to \SL(n, k)$ such that $u_*$ is equivalent to $u$ and 
$u_*(t)$ is a lower triangular matrix for each $t \in \G_a$. 
\end{enumerate} 
\end{lem}

\Proof  
See \cite[Lemma 1.8]{Tanimoto 2019}.
\QED

We know the following classification of representations of $\G_a$ into $\SL(3, k)$: 

\begin{lem}
Let $u : \G_a \to \SL(3, k)$ be a morphism of affine $k$-varieties such that $u(t)$ is an upper {\rm (}resp. lower{\rm )} 
triangular matrix for each $t \in \G_a$. 
Then the following assertions {\rm (1)} and {\rm (2)} hold true: 
\begin{enumerate}[label = {\rm (\arabic*)}]
\item If $p = 2$, then $u : \G_a \to \SL(3, k)$ is a representation of $\G_a$ if and only if  $u \in \cA_{(1, 2)} \cup \cA_{(2, 1)}$ 
{\rm (}resp. $u \in \cA_{(1, 2)}^- \cup \cA_{(2, 1)}^-${\rm )}. 

\item  If $p \geq 3$, then $u : \G_a \to \SL(3, k)$ is a representation of $\G_a$ if and only if 
$u \in \cU_{[3]} \cup\cA_{(1, 2)} \cup \cA_{(2, 1)}$ 
{\rm (}resp. $u \in \cU_{[3]}^- \cup\cA_{(1, 2)}^- \cup \cA_{(2, 1)}^-${\rm )}. 
\end{enumerate} 
\end{lem}

\Proof See \cite[Theorem 6.1]{Tanimoto 2008}. 
\QED

\subsection{Representations of $\G_a \rtimes \G_m$ into $\SL(n, k)$}

Let $n \geq 1$ and let $\varphi : \G_a \rtimes \G_m \to \SL(n, k)$ be a representation of $\G_a \rtimes \G_m$ into $\SL(n, k)$. 
We can define a representation $h_\varphi : \G_m \to \SL(n, k)$ of $\G_m$ into $\SL(n, k)$ as  
\[
 h_\varphi (z) 
 := 
\varphi 
\left(
\begin{array}{c c}
 z & 0 \\
 0 & z^{-1}
\end{array}
\right) . 
\]
We can define a representation $u_\varphi : \G_a \to \SL(n, k)$ of $\G_a$ into $\SL(n, k)$ as  
\[
 u_\varphi (t) 
 := 
\varphi 
\left(
\begin{array}{c c}
 1 & t \\
 0 & 1
\end{array}
\right) . 
\]

\begin{lem}
Let $\varphi : \G_a \rtimes \G_m \to \GL(n, k)$ be a representation with the form 
\begin{align*}
 h_\varphi(z) 
 = 
\diag(\, z^{\ell_1}, \, z^{\ell_2}, \, \ldots, \, z^{\ell_n} \, ) \qquad  (\, \ell_1 \geq \ell_2 \geq \cdots \geq \ell_n \,) . 
\end{align*}
Let  $\lambda_i$ $(1 \leq i \leq N)$ be positive integers satisfying $\lambda_1 + \lambda_2 + \cdots + \lambda_N  = n$ and 
\begin{align*} 
& \overbrace{ \; \ell_1 = \cdots = \ell_{\lambda_1} \; }^{\lambda_1} 
 \; > \; \overbrace{ \; \ell_{\lambda_1 + 1} = \cdots = \ell_{\lambda_1 + \lambda_2 } \; }^{\lambda_2} 
 \; > \; 
 \cdots \\
& \qquad 
  \; > \; \overbrace{ \; \ell_{\lambda_1 + \lambda_2 + \cdots + \lambda_{N - 1}  + 1} 
  = \cdots = \ell_{\lambda_1 + \lambda_2 + \cdots + \lambda_{N - 1} + \lambda_N } \; }^{\lambda_N}   . 
\end{align*}
Write $u_\varphi(t)$ as 
\begin{align*}
  u_\varphi(t) 
 = 
\left(
\begin{array}{c | c  | c}
 A_{1, \, 1}(t) & \cdots & A_{1, \, N}(t) \\
\hline 
 \vdots & \ddots & \vdots \\ 
\hline 
 A_{N, \, 1}(t) & \cdots & A_{N, \, N}(t) 
\end{array}
\right) , 
\qquad 
A_{i, \, j}(t) \in \Mat_{\lambda_i, \, \lambda_j}(k) \qquad (\, 1 \leq i , j \leq N \, ) . 
\end{align*}
Then $A_{i,\, j}(t) = 0$ for all $ 1 \leq j < i \leq N$, 
and $A_{i, \, i}(t) = I_{\lambda_i} \in \Mat(\lambda_i, k)$ for all $1 \leq i \leq N$. 
In particular, for any $t \in k$, the matrix $u_\varphi(t)$ is an upper triangular matrix. 
\end{lem}

\Proof 
Since 
\[
\left(
\begin{array}{c c}
 z & 0 \\
 0 & z^{-1} 
\end{array}
\right)
\left(
\begin{array}{c c}
 1 & t \\
 0 & 1 
\end{array}
\right)
= 
\left(
\begin{array}{c c}
 1 & z^2 \, t \\
 0 & 1 
\end{array}
\right)
\left(
\begin{array}{c c}
 z & 0 \\
 0 & z^{-1} 
\end{array}
\right) , 
\]
we have 
\begin{align*}
&
 \diag(z^{\ell_1}, \ldots, z^{\ell_n}) 
\left(
\begin{array}{c | c  | c}
 A_{1, \, 1}(t) & \cdots & A_{1, \, N}(t) \\
\hline 
 \vdots & \ddots & \vdots \\ 
\hline 
 A_{N, \, 1}(t) & \cdots & A_{N, \, N}(t) 
\end{array}
\right) \\
& \qquad 
 = 
 \left(
\begin{array}{c | c  | c}
 A_{1, \, 1}(z^2 \, t) & \cdots & A_{1, \, N}(z^2 \, t) \\
\hline 
 \vdots & \ddots & \vdots \\ 
\hline 
 A_{N, \, 1}(z^2 \, t) & \cdots & A_{N, \, N}(z^2 \, t) 
\end{array}
\right) 
\diag(z^{\ell_1}, \ldots, z^{\ell_n}) . 
\end{align*}
Comparing the $(d, e)$-th submatrices of both sides of the above equality, we have 
\[
 z^{\ell_{\lambda_1 + \cdots + \lambda_{d - 1} + 1}} \, A_{d,\, e}(t) 
 = 
 A_{d, \, e}(z^2 \, t) \, z^{ \ell_{\lambda_1 + \cdots + \lambda_{e - 1} + 1} } . 
\]
If $ e < d$, then $\ell_{\lambda_1 + \cdots + \lambda_{e - 1} + 1 } >  \ell_{\lambda_1 + \cdots + \lambda_{d - 1} + 1}$ 
and thereby $A_{d,\, e}(t) = O$; and if $d = e$, then $A_{d, \, e}(t) \in \Mat(\lambda_d, k)$. 
So, each $A_{d, \, d} : \G_a \to \GL(\lambda_d, k)$ is a constant representations of $\G_a$. 
Thus $A_{d, \, d}(t) = I_{\lambda_d}$. 
\QED

\subsection{Representations of $\SL(2, k)$}

Let $n \geq 1$, let $G$ be an algebraic group over $k$, and let $\psi : \SL(2, k) \to G$ be a homomorphism of algebraic groups over $k$.  
We can define a homomorphism $h_\psi : \G_m \to G$ as  
\[
 h_\psi (z) 
 := 
\psi 
\left(
\begin{array}{c c}
 z & 0 \\
 0 & z^{-1}
\end{array}
\right) . 
\]
We can define a homomorphism $u_\psi : \G_a \to G$ as  
\[
 u_\psi (t) 
 := 
\psi 
\left(
\begin{array}{c c}
 1 & t \\
 0 & 1
\end{array}
\right) . 
\]
We can define a homomorphism $u_\psi^- : \G_a \to G$ as  
\[
 u_\psi^-(s) 
 := 
\psi  
\left(
\begin{array}{c c}
 1 & 0 \\
 s & 1
\end{array}
\right) . 
\]

Let $\psi : \SL(2, k) \to \GL(n, k)$ be a representation of $\SL(2, k)$. 
Let $V := k^{\oplus n}$ be the $n$-dimensional column vector space. 
So, the representation $\psi$ yields a linear action of $\SL(2, k)$ on $V$. 
For any $\ell \in \Z$, we let $V_\ell$ be the subspace of $V$ defined by 
\[
 V_\ell 
 := 
\left\{ \, v \in V  \; \left| \;   
h_\psi(z) \, v = z^\ell \, v \text{ \; for all \, $z \in \G_m$} 
\right. 
\right\} . 
\]
Clearly, $h_\psi$ is trivial if and only if $V = V_0$.

\begin{lem}
Let $\psi : \SL(2, k) \to \GL(n, k)$ be a representation. 
Then the following assertions {\rm (1)} and {\rm (2)} hold true: 
\begin{enumerate}[label = {\rm (\arabic*)}]
\item $\ds  V = \bigoplus_{\ell \in \Z} V_\ell$. 

\item 
Let 
\[
 J 
 := 
 \psi \left(
\begin{array}{r r}
 0 & 1 \\
 -1 & 0
\end{array}
\right)
\]
be the $k$-linear isomorphism from $V$ to itself. 
Then we have $J^2 = - \id_V$ and $J(V_\ell) = V_{- \ell}$ for all $\ell \in \Z$. 
In particular, the $k$-vector spaces $V_\ell$ and $V_{- \ell}$ are isomorphic, i.e., 
\[
 V_\ell \cong V_{- \ell} .
\]
\end{enumerate} 
\end{lem}

\Proof 
The proofs of assertions (1) and (2) are straightforward. 
\QED

\begin{lem}
Let $\psi : \SL(2, k) \to \GL(n, k)$ be a representation. 
Then we have $\psi(\SL(2, k)) \subset \SL(n, k)$. 
\end{lem}

\Proof 
$V$ has the decomposition satisfying the following conditions (1) and (2): 
\begin{enumerate}[label = {\rm (\arabic*)}]
\item $V = V_{m_1} \oplus V_{m_2} \oplus \cdots \oplus V_{m_r}$ \, ($m_1 > m_2 > \cdots > m_r$).

\item $V_i \ne 0$ for all $i \in \{m_1, \, m_2, \ldots, m_r \}$. 
\end{enumerate} 
We know from Lemma 1.6 that $m_i = - m_{r - i + 1}$ for all $1 \leq i \leq r$, and 
$\dim_k V_{m_i} = \dim_k V_{m_{r - i + 1}}$ for all $1 \leq i \leq r$. 
Thus for all $z \in \G_m$, we have 
\begin{align*}
 \det h_\psi(z)  =  \prod_{i = 1}^r z^{m_i \, \dim_k V_{m_i} }  = z^{d } , 
\end{align*}
where 
\[
 d := \sum_{i  = 1}^r  m_i \, \dim_k V_{m_i} . 
\]
Then 
\begin{align*}
2 \, d
 & = \sum_{i  = 1}^r  m_i \, \dim_k V_{m_i}  + \sum_{i  = 1}^r  ( - m_{r - i + 1} ) \, \dim_k V_{m_{r - i + 1}} \\
 & = \sum_{i  = 1}^r  m_i \, \dim_k V_{m_i}  -  \sum_{i  = 1}^r  m_i \, \dim_k V_{m_i}  \\
 & = 0 , 
\end{align*} 
which implies $\det h(z) = 1$. 
Let $D(a)$ be the affine open subset of $\SL(2, k)$ defined by 
\[
 D(a) 
 := 
\left\{ 
\left. 
\left(
\begin{array}{c c}
 a & b \\
 c & d 
\end{array} 
\right) 
\in \SL(2, k) 
\; \right| \; 
a \ne 0 
\right\} . 
\]
For any element $\left(
\begin{array}{c c}
 a & b \\
 c & d 
\end{array} 
\right) $ of $D(a)$, we have 
\[
 \left(
\begin{array}{c c}
 a & b \\
 c & d 
\end{array} 
\right) 
 = 
\left(
\begin{array}{c c}
 1 & 0 \\
 \frac{c}{a} & 1 
\end{array} 
\right) 
\left(
\begin{array}{c c}
 a & 0 \\
 0 & a^{-1} 
\end{array} 
\right) 
\left(
\begin{array}{c c}
 1 & \frac{b}{a} \\
 0 & 1
\end{array} 
\right) , 
\]
which implies 
\[
 \det \, \psi \left(
\begin{array}{c c}
 a & b \\
 c & d 
\end{array} 
\right) 
= 1 . 
\]
Thus $\psi(D(a)) \subset \SL(n, k)$. Since $\psi$ is a continuous map, we have 
\[
 \psi(\ol{D(a)}) \subset \ol{\SL(n, k)}
\]
and thereby have the desired inclusion $\psi(\SL(2, k)) \subset \SL(n, k)$. 
\QED

\begin{lem}
Let $\psi : \SL(2, k) \to \GL(n, k)$ be a representation with the form 
\begin{align*}
 h_\psi(z) 
 = 
\diag(\, z^{\ell_1}, \, z^{\ell_2}, \, \ldots, \, z^{\ell_n} \, ) \qquad  (\, \ell_1 \geq \ell_2 \geq \cdots \geq \ell_n \,) . 
\end{align*}
Let  $\lambda_i$ $(1 \leq i \leq N)$ be positive integers satisfying $\lambda_1 + \lambda_2 + \cdots + \lambda_N  = n$ and 
\begin{align*} 
& \overbrace{ \; \ell_1 = \cdots = \ell_{\lambda_1} \; }^{\lambda_1} 
 \; > \; \overbrace{ \; \ell_{\lambda_1 + 1} = \cdots = \ell_{\lambda_1 + \lambda_2 } \; }^{\lambda_2} 
 \; > \; 
 \cdots \\
& \qquad 
  \; > \; \overbrace{ \; \ell_{\lambda_1 + \lambda_2 + \cdots + \lambda_{N - 1}  + 1} 
  = \cdots = \ell_{\lambda_1 + \lambda_2 + \cdots + \lambda_{N - 1} + \lambda_N } \; }^{\lambda_N}   . 
\end{align*}
Then the following assertions {\rm (1)} and {\rm (2)} hold true: 
\begin{enumerate}[label = {\rm (\arabic*)}]
\item
Write $u_\psi(t)$ as 
\begin{align*}
  u_\psi(t) 
 = 
\left(
\begin{array}{c | c  | c}
 A_{1, \, 1}(t) & \cdots & A_{1, \, N}(t) \\
\hline 
 \vdots & \ddots & \vdots \\ 
\hline 
 A_{N, \, 1}(t) & \cdots & A_{N, \, N}(t) 
\end{array}
\right) , 
\qquad 
A_{i, \, j}(t) \in \Mat_{\lambda_i, \, \lambda_j}(k) \qquad (\, 1 \leq i , j \leq N \, ) . 
\end{align*}
Then $A_{i,\, j}(t) = 0$ for all $ 1 \leq j < i \leq N$, 
and $A_{i, \, i}(t) = I_{\lambda_i} \in \Mat(\lambda_i, k)$ for all $1 \leq i \leq N$. 
In particular, for any $t \in k$, the matrix $u_\psi(t)$ is an upper triangular matrix.

\item 
Write $u_\psi^-(s)$ as 
\begin{align*}
  u_\psi^-(s) 
 = 
\left(
\begin{array}{c | c  | c}
 B_{1, \, 1}(s) & \cdots & B_{1, \, N}(s) \\
\hline 
 \vdots & \ddots & \vdots \\ 
\hline 
 B_{N, \, 1}(s) & \cdots & B_{N, \, N}(s) 
\end{array}
\right) , 
\qquad 
B_{i, \, j}(s) \in \Mat_{\lambda_i, \, \lambda_j}(k) \qquad (\, 1 \leq i , j \leq N \, ) . 
\end{align*}
Then $B_{i,\, j}(s) = 0$ for all $ 1 \leq j < i \leq N$, 
and $B_{i, \, i}(t) = I_{\lambda_i} \in \Mat(\lambda_i, k)$ for all $1 \leq i \leq N$. 
In particular, for any $s \in k$, the matrix $u_\psi^-(s)$ is a lower triangular matrix. 
\end{enumerate} 
\end{lem}

\Proof 
See the proof of Lemma 1.5. 
\QED

\begin{lem}
Let $\psi : \SL(2, k) \to \GL(n, k)$ be a representation such that $u_\psi$ is trivial. 
Then $u_\psi^-$ is trivial and  $h_\psi$ is also trivial. 
\end{lem}

\Proof 
Note that 
\[
\left(
\begin{array}{r r}
 0 & - 1 \\
 1 & 0
\end{array}
\right) 
\left(
\begin{array}{r r}
 1 & - t \\
 0 & 1
\end{array}
\right) 
\left(
\begin{array}{r r}
 0 & 1 \\
 - 1 & 0
\end{array}
\right) 
= 
\left(
\begin{array}{r r}
 1 & 0 \\
 t & 1
\end{array}
\right) 
\qquad 
\text{ for all \quad $t \in \G_a$. } 
\]
Applying $\psi$ to the above equality, we have 
\[
\psi 
\left(
\begin{array}{r r}
 0 & - 1 \\
 1 & 0
\end{array}
\right) 
\, 
\psi 
\left(
\begin{array}{r r}
 1 & - t \\
 0 & 1
\end{array}
\right) 
\, 
\psi 
\left(
\begin{array}{r r}
 0 & 1 \\
 - 1 & 0
\end{array}
\right) 
= 
\psi
\left(
\begin{array}{r r}
 1 & 0 \\
 t & 1
\end{array}
\right) 
\qquad 
\text{ for all \quad $t \in \G_a$, } 
\] 
and thereby have $u^-(t) = I_n$ for all $t \in \G_a$. .

Note that 
\[
\left(
\begin{array}{c c}
 1 & 1 \\
 0 & 1
\end{array}
\right) 
\left(
\begin{array}{c c}
 1 & 0 \\
 \gamma & 1
\end{array}
\right) 
 = 
\left(
\begin{array}{c c}
 1 & 0 \\
 \frac{\gamma}{1 + \gamma} & 1
\end{array}
\right)  
\left(
\begin{array}{c c}
 1 + \gamma & 0 \\
 0  & \frac{1}{1 + \gamma} 
\end{array}
\right) 
\left(
\begin{array}{c c}
 1 & \frac{1}{1 + \gamma} \\
 0 & 1
\end{array}
\right) 
\qquad 
\text{ for all \quad $\gamma \in k \backslash \{ -1 \}$.  } 
\]
Applying $\psi$ to the above equality, we have 
$I_n \cdot I_n = I_n \cdot h_\psi(1 + \gamma) \cdot I_n$ for all $\gamma \in  k \backslash \{ -1 \}$. 
Thus $h_\psi$ is trivial. 
\QED

\begin{lem}
Let $G$ be an algebraic group over $k$. 
Let $\psi_i : \SL(2, k) \to G$ $(i = 1, 2)$ be homomorphisms of algebraic groups 
satisfying the following conditions {\rm (1)} and {\rm (2)}: 
\begin{enumerate}[label = {\rm (\arabic*)}]
\item For any 
$\left(
\begin{array}{c c}
 a & b \\
 0 & d
\end{array}
\right) 
\in \SL(2, k)$, 
the equality 
$\psi_1 \left(
\begin{array}{c c}
 a & b \\
 0 & d
\end{array}
\right)
= 
\psi_2 \left(
\begin{array}{c c}
 a & b \\
 0 & d
\end{array}
\right)
$ 
holds true.

\item For any 
$\left(
\begin{array}{c c}
 1 & 0 \\
 s & 1
\end{array}
\right) 
\in \SL(2, k)$, 
the equality 
$\psi_1 \left(
\begin{array}{c c}
 1 & 0 \\
 s & 1
\end{array}
\right) 
= 
\psi_2 \left(
\begin{array}{c c}
 1 & 0 \\
 s & 1
\end{array}
\right) 
$ 
holds true. 
\end{enumerate} 
Then we have $\psi_1 = \psi_2$. 
\end{lem}

\Proof 
Choose an arbitrary matrix 
$
\left(
\begin{array}{c c}
 a & b \\
 c & d
\end{array}
\right)
$ 
of $\SL(2, k)$. 
In the case where $a = 0$, we have $b \, c = - 1$ and 
\begin{align*}
\psi_1 
\left(
\begin{array}{r r}
 0 & b \\
 c & d
\end{array}
\right)
 & = 
\psi_1 
\left(
\begin{array}{r r}
 1 & 1 \\
 0 & 1
\end{array}
\right) 
\psi_1 
\left(
\begin{array}{r r}
 1 & 0 \\
 -1 & 1
\end{array}
\right) 
\psi_1 
\left(
\begin{array}{r r}
 1 & 1 \\
 0 & 1
\end{array}
\right) 
\psi_1 \left(
\begin{array}{r r}
 - c & - d \\
 0 & b
\end{array}
\right)  \\
& = 
\psi_2 
\left(
\begin{array}{r r}
 1 & 1 \\
 0 & 1
\end{array}
\right) 
\psi_2 
\left(
\begin{array}{r r}
 1 & 0 \\
 -1 & 1
\end{array}
\right) 
\psi_2 
\left(
\begin{array}{r r}
 1 & 1 \\
 0 & 1
\end{array}
\right) 
\psi_2 
\left(
\begin{array}{r r}
 - c & - d \\
 0 & b
\end{array}
\right)  \\
 & = 
\psi_2 
\left(
\begin{array}{r r}
 0 & b \\
 c & d
\end{array}
\right) . 
\end{align*} 
In the case where $a \ne 0$, we have the equality 
\[
\left(
\begin{array}{r r}
 a & b \\
 c & d
\end{array}
\right)
 = 
\left(
\begin{array}{r r}
 1 & 0 \\
 \frac{c}{a} & 1
\end{array}
\right)
\left(
\begin{array}{c c}
 a & 0 \\
 0 & a^{-1} 
\end{array}
\right)
\left(
\begin{array}{r r}
 1 & \frac{b}{a} \\
 0 & 1
\end{array}
\right) . 
\]
Thus,  
\begin{align*}
\psi_1 
\left(
\begin{array}{r r}
 a & b \\
 c & d
\end{array}
\right)
 &= 
\psi_1 
\left(
\begin{array}{r r}
 1 & 0 \\
 \frac{c}{a} & 1
\end{array}
\right) 
\psi_1 
\left(
\begin{array}{c c}
 a & 0 \\
 0 & a^{-1} 
\end{array}
\right) 
\psi_1
\left(
\begin{array}{r r}
 1 & \frac{b}{a} \\
 0 & 1
\end{array}
\right) \\
 &= 
\psi_2
\left(
\begin{array}{r r}
 1 & 0 \\
 \frac{c}{a} & 1
\end{array}
\right) 
\psi_2 
\left(
\begin{array}{c c}
 a & 0 \\
 0 & a^{-1} 
\end{array}
\right) 
\psi_2 
\left(
\begin{array}{r r}
 1 & \frac{b}{a} \\
 0 & 1
\end{array}
\right) \\
 & =
\psi_2 
\left(
\begin{array}{r r}
 a & b \\
 c & d
\end{array}
\right) . 
\end{align*}
\QED

\section{Representations of $\G_a \rtimes \G_m$ into $\SL(3, k)$}

\subsection{Candidates for classifying representations of $\G_a \rtimes \G_m$ into $\SL(3, k)$}

\begin{lem}
Let $\varphi : \G_a \rtimes \G_m \to \SL(3, k)$ be a representation such that 
$h_\varphi$ has the form
\[
h_\varphi(z) 
= 
\diag(\, z^{\ell_1}, \, z^{\ell_2}, \, z^{\ell_3} \, ) 
\qquad 
 (\, \ell_1 \geq \ell_2 \geq \ell_3 \,) . 
\]
Then the following assertions {\rm (1)}, {\rm (2)}, {\rm(3)} hold true: 
\begin{enumerate}[label = {\rm (\arabic*)}]

\item 
Assume $\ell_1 > \ell_2 > \ell_3$. 

\begin{enumerate}[label = {\rm (1.\arabic*)}]  
\item If $u_\varphi \in \cU_{[3]}$ where $p \geq 3$, then we can express $u_\varphi$ as  
\begin{align*}
u_\varphi(t) 
&= 
\left(
\begin{array}{c c c}
 1 & c_1 \, t^{p^{e_1}} & \frac{1}{2} \, \lambda \, c_1^2 \, t^{2 \, p^{e_1}} \\
 0 & 1 & \lambda \, c_1 \, t^{p^{e_1}}  \\
 0 & 0 & 1
\end{array}
\right) 
\qquad 
(\, c_1 \in k \backslash \{0\} , \quad \lambda \in k \backslash \{ 0 \} , \quad e_1 \geq 0 \,) , 
\end{align*}
the $3$-tuples of $(\ell_1, \ell_2, \ell_3)$ of $\ell_1$, $\ell_2$, $\ell_3$ appearing in $h_\varphi(z)$ as
\[
(\ell_1, \; \ell_2, \; \ell_3) 
 = 
( 2 \, p^{e_1}, \;  0, \;  -2 \, p^{e_1} )
\] 
and the representation $\varphi$ as 
\begin{align*} 
\varphi
\left(
\begin{array}{c c}
 a & b \\
 0 & d 
\end{array}
\right) 
 & =
\left(
\begin{array}{c c c}
 a^{2 \, p^{e_1}}  & c_1 \, a^{p^{e_1}} \, b^{p^{e_1}}  & \frac{1}{2} \, \lambda \, c_1^2 \, b^{2 \, p^{e_1}} \\
 0  &  1  &   \lambda \, c_1  \,  b^{p^{e_1}} \,  d^{p^{e_1}} \\
 0 & 0 & d^{2 \, p^{e_1}} 
\end{array}  
\right) . 
\end{align*}

\item If $u_\varphi \in \cA_{(1, 2)}$ where $p \geq 2$, then we can express $u_\varphi$ as 
\begin{align*}
u_\varphi(t) 
 & = 
\left(
\begin{array}{c c c}
 1 & c_1 \, t^{p^{e_1}} & c_2 \, t^{p^{e_2}} \\
 0 & 1 & 0 \\
 0 & 0 & 1
\end{array}
\right) 
\qquad 
(\, c_1, c_2 \in k , \quad e_1, e_2 \geq 0 \,) , 
\end{align*} 
and we can express $\varphi$ by separating the following four cases {\rm (1.2.a), (1.2.b), (1.2.c), (1.2.d)}:  
\begin{enumerate}[label = {\rm (1.2.\alph*)}, leftmargin = * ]
\item If $c_1 = 0$ and $c_2 = 0$, then $\ell_1 > 0 > \ell_3$ and 
\begin{fleqn}[36pt] 
\begin{align*}
\varphi 
\left(
\begin{array}{c c}
 a & b \\
 0 & d 
\end{array}
\right) 
& =
\left\{
\begin{array}{l @{\qquad} l} 
\left(
\begin{array}{c c c}
 a^{\ell_1} & 0 & 0 \\
 0 & a^{\ell_2} & 0 \\
 0 & 0 & d^{ - \ell_3} 
\end{array} 
\right)
 &  (\, \ell_2 \geq 0 \, ) , \\ [24pt] 
\left(
\begin{array}{c c c}
 a^{\ell_1} & 0 & 0 \\
 0 & d^{- \ell_2} & 0 \\
 0 & 0 & d^{ - \ell_3} 
\end{array} 
\right)
 &  (\, \ell_2 < 0 \, ) . 
\end{array}
\right. 
\end{align*} 
\end{fleqn}

\item If $c_1 \ne 0$ and $c_2 = 0$, then $\ell_1 - \ell_2 = 2 \, p^{e_1}$, $\ell_1 > 0 > \ell_3$ and 
\begin{fleqn}[36pt] 
\begin{align*}
\varphi
\left(
\begin{array}{c c}
 a & b \\
 0 & d 
\end{array}
\right) 
& =
\left\{
\begin{array}{l @{\qquad} l} 
\left(
\begin{array}{c c c}
 a^{\ell_1} & c_1 \, a^{\ell_2 + p^{e_1}} \, b^{p^{e_1}}  & 0 \\
 0 & a^{\ell_2} & 0 \\
 0 & 0 & d^{ - \ell_3} 
\end{array} 
\right) 
&  ( \, \ell_2 \geq 0 \, ) ,  \\  [24pt] 
\left(
\begin{array}{c c c}
 a^{\ell_1} & c_1 \, a^{p^{e_1}} \, b^{p^{e_1}} \, d^{- \ell_2}  & 0 \\
 0 & d^{-\ell_2} & 0 \\
 0 & 0 & d^{ - \ell_3} 
\end{array}
\right) 
& (\,  \ell_2 < 0 \, ) .  
\end{array} 
\right. 
\end{align*}
\end{fleqn}

\item If $c_1 = 0$ and $c_2 \ne 0$, then $\ell_1 - \ell_3 = 2 \, p^{e_2}$, $\ell_1 > 0 > \ell_3$ and  
\begin{fleqn}[36pt] 
\begin{align*}
\varphi
\left(
\begin{array}{c c}
 a & b \\
 0 & d 
\end{array}
\right) 
& =
\left\{
\begin{array}{l @{\qquad} l} 
\left(
\begin{array}{c c c}
 a^{\ell_1} & 0 &  c_2 \, a^{p^{e_2}} \, b^{p^{e_2}} \, d^{- \ell_3} \\
 0 & a^{\ell_2} & 0 \\
 0 & 0 &  d^{ - \ell_3} 
\end{array} 
\right) 
&  ( \, \ell_2 \geq 0 \, ) ,  \\  [24pt] 
\left(
\begin{array}{c c c}
 a^{\ell_1} & 0 & c_2 \, a^{p^{e_2}} \, b^{p^{e_2}} \, d^{- \ell_3}   \\
 0 & d^{-\ell_2} & 0 \\
 0 & 0 & d^{ - \ell_3} 
\end{array}
\right) 
& (\,  \ell_2 < 0 \, ) .  
\end{array} 
\right. 
\end{align*}
\end{fleqn}

\item If $c_1 \ne 0$ and $c_2 \ne 0$, then we have 
\[
(\ell_1, \; \ell_2, \; \ell_3) 
 =
 \left( 
  \frac{2 \, p^{e_1} + 2 \, p^{e_2}}{3} , \quad
  \frac{-4 \, p^{e_1} + 2 \, p^{e_2}}{3} , \quad 
 \frac{2 \, p^{e_1} - 4 \, p^{e_2}}{3} 
\right) . 
\]
So, $\ell_1 > 0$, $\ell_2 \geq 0$, $\ell_3 < 0$ and $e_2 > e_1 \geq 0$. 
And we can express $\varphi$ as 
\begin{fleqn}[36pt] 
\begin{align*}
\varphi
\left(
\begin{array}{c c}
 a & b \\
 0 & d 
\end{array}
\right) 
  & = 
\left(
\begin{array}{c c c}
 a^{\ell_1}  &  c_1 \,a^{\ell_1} \, b^{p^{e_1}} \, d^{p^{e_1}}  &  c_2 \,a^{\ell_1} \, b^{p^{e_2}} \, d^{p^{e_2}}   \\
 0  &  a^{\ell_2}  &  0  \\
 0 & 0 & d^{- \ell_3} 
\end{array}  
\right) .  
\end{align*}
\end{fleqn} 
\end{enumerate}

\item If $u_\varphi \in \cA_{(2, 1)}$ where $p \geq 2$, then we can express $u_\varphi$ as 
\begin{align*}
u_\varphi(t) 
 & = 
\left(
\begin{array}{c c c}
 1 & 0 & c_2 \, t^{p^{e_2}} \\
 0 & 1 & c_1 \, t^{p^{e_1}} \\
 0 & 0 & 1
\end{array}
\right) 
\qquad 
(\, c_1, c_2 \in k , \quad e_1, e_2 \geq 0 \,) , 
\end{align*} 
and we can express $\varphi$ by separating the following four cases {\rm (1.3.a), (1.3.b), (1.3.c), (1.3.d)}:  
\begin{enumerate}[label = {\rm (1.3.\alph*)}, leftmargin = * ]
\item If $c_1 = 0$ and $c_2 = 0$, then $\ell_1 > 0 > \ell_3$ and 
\begin{fleqn}[36pt] 
\begin{align*}
\varphi 
\left(
\begin{array}{c c}
 a & b \\
 0 & d 
\end{array}
\right) 
& =
\left\{
\begin{array}{l @{\qquad} l} 
\left(
\begin{array}{c c c}
 a^{\ell_1} & 0 & 0 \\
 0 & a^{\ell_2} & 0 \\
 0 & 0 & d^{ - \ell_3} 
\end{array} 
\right)
 &  (\, \ell_2 \geq 0 \, ) , \\  [24pt] 
\left(
\begin{array}{c c c}
 a^{\ell_1} & 0 & 0 \\
 0 & d^{- \ell_2} & 0 \\
 0 & 0 & d^{ - \ell_3} 
\end{array} 
\right)
 &  (\, \ell_2 < 0 \, ) . 
\end{array}
\right. 
\end{align*} 
\end{fleqn}

\item If $c_1 \ne 0$ and $c_2 = 0$, then $\ell_2 - \ell_3 = 2 \, p^{e_1}$, $\ell_1 > 0 > \ell_3$ and 
\begin{fleqn}[36pt] 
\begin{align*}
\varphi
\left(
\begin{array}{c c}
 a & b \\
 0 & d 
\end{array}
\right) 
& =
\left\{
\begin{array}{l @{\qquad} l} 
\left(
\begin{array}{c c c}
 a^{\ell_1} & 0  & 0 \\
 0 & a^{\ell_2} & c_1 \, a^{ p^{e_1}} \, b^{p^{e_1}} \, d^{- \ell_3}  \\
 0 & 0 & d^{ - \ell_3} 
\end{array} 
\right) 
&  ( \, \ell_2 \geq 0 \, ) ,  \\  [24pt] 
\left(
\begin{array}{c c c}
 a^{\ell_1} & 0  & 0 \\
 0 & d^{-\ell_2} & c_1 \, a^{ p^{e_1}} \, b^{p^{e_1}} \, d^{- \ell_3} \\
 0 & 0 & d^{ - \ell_3} 
\end{array}
\right) 
& (\,  \ell_2 < 0 \, ) .  
\end{array} 
\right. 
\end{align*}
\end{fleqn}

\item If $c_1 = 0$ and $c_2 \ne 0$, then $\ell_1 - \ell_3 = 2 \, p^{e_2}$, $\ell_1 > 0 > \ell_3$ and 
\begin{fleqn}[36pt] 
\begin{align*}
\varphi
\left(
\begin{array}{c c}
 a & b \\
 0 & d 
\end{array}
\right) 
& =
\left\{
\begin{array}{l @{\qquad} l} 
\left(
\begin{array}{c c c}
 a^{\ell_1} & 0  & c_2 \, a^{ p^{e_2}} \, b^{p^{e_2}}  \, d^{-\ell_3} \\
 0 & a^{\ell_2} & 0  \\
 0 & 0 & d^{ - \ell_3} 
\end{array} 
\right) 
&  ( \, \ell_2 \geq 0 \, ) ,  \\  [24pt] 
\left(
\begin{array}{c c c}
 a^{\ell_1} & 0  & c_2 \, a^{p^{e_2}} \, b^{p^{e_2}} \, d^{-\ell_3}   \\
 0 & d^{-\ell_2} & 0 \\
 0 & 0 & d^{ - \ell_3} 
\end{array}
\right) 
& (\,  \ell_2 < 0 \, ) .  
\end{array} 
\right. 
\end{align*}
\end{fleqn}

\item If $c_1 \ne 0$ and $c_2 \ne 0$, then we have 
\[
(\ell_1, \; \ell_2, \; \ell_3) 
 = 
\left( 
  \frac{ - 2 \, p^{e_1} + 4 \, p^{e_2}}{3} , \quad 
  \frac{4 \, p^{e_1} - 2 \, p^{e_2}}{3} , \quad 
  \frac{ - 2 \, p^{e_1} - 2 \, p^{e_2}}{3}  
\right) . 
\]
So, $\ell_1 > 0$, $\ell_2 \leq 0$, $\ell_3 < 0$ and $e_2 > e_1 \geq 0$. 
And we can express $\varphi$ as 
\begin{fleqn}[36pt] 
\begin{align*}
\varphi
\left(
\begin{array}{c c}
 a & b \\
 0 & d 
\end{array}
\right) 
 & =
\left(
\begin{array}{c c c}
 a^{\ell_1}  &  0  &  c_2 \,a^{p^{e_2}} \, b^{p^{e_2}} \, d^{- \ell_3}   \\
 0  &  d^{ - \ell_2}  &  c_1 \, a^{p^{e_1}} \, b^{p^{e_1}} \, d^{- \ell_3}  \\
 0 & 0 & d^{ - \ell_3} 
\end{array}  
\right) . 
\end{align*} 
\end{fleqn} 
\end{enumerate}
\end{enumerate}

\item 
Assume $\ell_1 = \ell_2 > \ell_3$. 
Then we have $u_\varphi \in \cA_{(2, 1)}$, we can express $u_\varphi$ as 
\begin{align*}
 u_\varphi(t) 
& = 
\left(
\begin{array}{c c c}
 1 & 0 & c_2 \, t^{p^{e_2}} \\
 0 & 1 & c_1 \, t^{p^{e_1}} \\
 0 & 0 & 1
\end{array}
\right) 
\qquad 
(\, c_1, c_2 \in k  , \quad e_1, e_2 \geq 0 \,) , 
\end{align*} 
and we can express $\varphi$ by separating the following four cases {\rm (2.a), (2.b), (2.c), (2.d)}:  
\begin{enumerate}[label = {\rm (2.\alph*)}, leftmargin = * ]
\item If $c_1 = 0$ and $c_2 = 0$, then $\ell_1 = \ell_2 > 0 > \ell_3$ and  
\begin{fleqn}[36pt] 
\begin{align*}
\varphi
\left(
\begin{array}{c c}
 a & b \\
 0 & d 
\end{array}
\right) 
 =
\left(
\begin{array}{c c c}
 a^{\ell_1}  &  0  &  0  \\
 0  &  a^{\ell_1}  &  0  \\
 0 & 0 & d^{- \ell_3} 
\end{array}  
\right) . 
\end{align*}
\end{fleqn}

\item If $c_1 \ne  0$ and $c_2 = 0$, then $\ell_2 - \ell_3 = 2 \, p^{e_1}$, $\ell_1 = \ell_2 > 0 > \ell_3$ and  
\begin{fleqn}[36pt] 
\begin{align*} 
\varphi
\left(
\begin{array}{c c}
 a & b \\
 0 & d 
\end{array}
\right) 
  & = 
\left(
\begin{array}{c c c}
 a^{\ell_1}  &  0  &  0  \\
 0  &  a^{\ell_1}  & c_1 \, a^{\ell_1} \, b^{p^{e_1}} \, d^{ p^{e_1} }  \\
 0 & 0 & d^{- \ell_3} 
\end{array}  
\right) . 
\end{align*}
\end{fleqn}

\item If $c_1 =  0$ and $c_2 \ne 0$, then $\ell_1 - \ell_3 = 2 \, p^{e_2}$, $\ell_1 = \ell_2 > 0 > \ell_3$ and  
\begin{fleqn}[36pt] 
\begin{align*} 
\varphi
\left(
\begin{array}{c c}
 a & b \\
 0 & d 
\end{array}
\right) 
  & = 
\left(
\begin{array}{c c c}
 a^{\ell_1}  &  0  &  c_2 \, a^{\ell_1} \, b^{p^{e_2}} \, d^{ p^{e_2} }   \\
 0  &  a^{\ell_1}  & 0 \\
 0 & 0 & d^{- \ell_3} 
\end{array}  
\right) . 
\end{align*}
\end{fleqn}

\item If $c_1 \ne  0$ and $c_2 \ne 0$, then 
\[
 p = 3 ,  \qquad e_1 = e_2  \geq 1, \qquad 
 (\ell_1, \; \ell_2, \; \ell_3) 
 = ( 2 \, p^{e_1 - 1} , \; 2 \, p^{e_1 - 1} , \; - 4 \, p^{e_1 - 1} ) , 
\] 
and 
\begin{fleqn}[36pt] 
\begin{align*} 
\varphi
\left(
\begin{array}{c c}
 a & b \\
 0 & d 
\end{array}
\right) 
 & = 
\left(
\begin{array}{c c c}
 a^{2 \,  p^{e_1- 1}}   &  0  &  c_2 \, a^{2 \, p^{e_1 - 1}} \, b^{p^{e_1}} \, d^{p^{e_1}}  \\
 0  & a^{2 \, p^{e_1- 1}}  & c_1 \, a^{2 \, p^{e_1 - 1}} \, b^{p^{e_1}} \, d^{p^{e_1}} \\
 0 & 0 & d^{4 \, p^{e_1 - 1}} 
\end{array}
\right) . 
\end{align*}
\end{fleqn} 
\end{enumerate}

\item 
Assume $\ell_1 > \ell_2 = \ell_3$. 
Then we have $u_\varphi \in \cA_{(1, 2)}$, we can express $u_\varphi$ as 
\begin{align*} 
u_\varphi(t) 
 & = 
\left(
\begin{array}{c c c}
 1 & c_1 \, t^{p^{e_1}}  & c_2 \, t^{p^{e_2}} \\
 0 & 1 & 0 \\
 0 & 0 & 1
\end{array}
\right)  
\qquad 
(\, c_1, c_2 \in k  , \quad e_1, e_2 \geq 0 \,) , 
\end{align*} 
and we can express $\varphi$ by separating the following four cases {\rm (3.a), (3.b), (3.c), (3.d)}:  
\begin{enumerate}[label = {\rm (3.\alph*)}, leftmargin = * ]
\item If $c_1 = 0$ and $c_2 = 0$, then $\ell_1 > 0 > \ell_2 = \ell_3$ and 
\begin{fleqn}[36pt]
\begin{align*}
\varphi
\left(
\begin{array}{c c}
 a & b \\
 0 & d 
\end{array}
\right) 
 =
\left(
\begin{array}{c c c}
 a^{\ell_1}  &  0  &  0  \\
 0  &  d^{- \ell_2}  &  0  \\
 0 & 0 & d^{- \ell_2} 
\end{array}  
\right) . 
\end{align*}
\end{fleqn}

\item If $c_1 \ne 0$ and $c_2 = 0$, then $\ell_1 - \ell_2 = 2 \, p^{e_1}$, $\ell_1 > 0 > \ell_2 = \ell_3$ and 
\begin{fleqn}[36pt]
\begin{align*} 
\varphi
\left(
\begin{array}{c c}
 a & b \\
 0 & d 
\end{array}
\right) 
 & = 
\left(
\begin{array}{c c c}
 a^{\ell_1}  &  c_1 \, a^{\ell_1} \,  b^{p^{e_1}} \, d^{ p^{e_1} }      &  0  \\
 0  &  d^{ - \ell_2}  & 0\\
 0 & 0 & d^{- \ell_2} 
\end{array}  
\right)   . 
\end{align*}
\end{fleqn}

\item If $c_1 = 0$ and $c_2 \ne 0$, then $\ell_1 - \ell_3 = 2 \, p^{e_2}$, $\ell_1 > 0 > \ell_2 = \ell_3$ and 
\begin{fleqn}[36pt]
\begin{align*} 
\varphi
\left(
\begin{array}{c c}
 a & b \\
 0 & d 
\end{array}
\right) 
  & = 
\left(
\begin{array}{c c c}
 a^{\ell_1}  &  0  &  c_2 \, a^{\ell_1} \, b^{p^{e_2}} \, d^{ p^{e_2} }\\
 0  &  d^{ - \ell_2}  & 0   \\
 0 & 0 & d^{- \ell_2} 
\end{array}  
\right)   . 
\end{align*}
\end{fleqn}

\item If $c_1 \ne 0$ and $c_2 \ne 0$, then 
\[
 p = 3 , \qquad
 e_1 \geq 1 , \qquad 
 (\ell_1, \; \ell_2, \; \ell_3) 
 = 
 (4 \, p^{e_1 - 1} , \;    - 2 \, p^{e_1 - 1} , \;  - 2 \, p^{e_1 - 1} ) , 
\] 
and 
\begin{fleqn}[36pt]
\begin{align*}
\varphi
\left(
\begin{array}{c c}
 a & b \\
 0 & d 
\end{array}
\right)  
 & = 
\left(
\begin{array}{c c c}
  a^{4 \, p^{e_1 - 1}}   &   c_1 \, a^{4 \, p^{e_1 - 1}} \, b^{p^{e_1}} \, d^{p^{e_1}}    &  c_2 \, a^{4 \, p^{e_1 - 1}} \, b^{p^{e_1}} \, d^{p^{e_1}}  \\
 0  & d^{2 \, p^{e_1- 1}}  & 0 \\
 0 & 0 & d^{2 \,  p^{e_1- 1}} 
\end{array}
\right) . 
\end{align*}
\end{fleqn} 
\end{enumerate}

\item Assume $\ell_1 = \ell_2 = \ell_3$. Then $u_\varphi$ is the trivial representation of $\G_a$, 
and $\varphi$ is the trivial representation of $\G_a \rtimes \G_m$. 
\end{enumerate}

\end{lem}

\subsubsection{Proof of assertion (1) of Lemma 2.1}

\paragraph{(1.1)} 

Applying $\varphi$ to the equality 
\[
\left(
\begin{array}{c c}
 z & 0 \\
 0 & z^{-1}  
\end{array}
\right)
\left(
\begin{array}{c c}
 1 & t \\
 0 & 1  
\end{array}
\right)
=
\left(
\begin{array}{c c}
 1 & z^2 \, t \\
 0 & 1  
\end{array}
\right)
\left(
\begin{array}{c c}
 z & 0 \\
 0 & z^{-1}  
\end{array}
\right) ,
\]
we have 
\begin{align}
& 
\left( 
\begin{array}{ccc}
 z^{\ell_1} & 0 & 0 \\
 0 & z^{\ell_2} & 0 \\
 0 & 0 & z^{\ell_3} 
\end{array}
\right) 
\left( 
\begin{array}{ccc}
 1 & \alpha_1(t) & \frac{\lambda}{2} \, \alpha_1(t)^2 + \alpha_2(t)  \\
 0 & 1 & \lambda \, \alpha_1(t) \\
 0 & 0 & 1 
\end{array}
\right)  \notag  \\
& \qquad 
 = 
\left( 
\begin{array}{ccc}
 1 & \alpha_1(z^2 \, t) & \frac{\lambda}{2} \, \alpha_1(z^2 \, t)^2 + \alpha_2(z^2 \, t)  \\
 0 & 1 & \lambda \, \alpha_1(z^2 \, t) \\
 0 & 0 & 1 
\end{array}
\right) 
\left( 
\begin{array}{ccc}
 z^{\ell_1} & 0 & 0 \\
 0 & z^{\ell_2} & 0 \\
 0 & 0 & z^{\ell_3} 
\end{array}
\right) . \tag{$\ast$}
\end{align} 
Comparing the $(1, 2)$-th entries of both sides of the equality $(\ast)$, we have 
$z^{\ell_1} \, \alpha_1(t) = \alpha_1(z^2 \, t) \, z^{\ell_2}$. 
Thus $z^{\ell_1 - \ell_2} \, \alpha_1(t) = \alpha_1(z^2 \, t)$. 
So, $\alpha_1(t)$ is a $p$-monomial.  
We can express $\alpha_1(t)$ as 
$\alpha_1(t) = c_1 \, t^{p^{e_1}}$ ($c_1 \in k \backslash \{ 0 \}$, $e_1 \geq 0$).  
So, $\ell_1 - \ell_2 = 2 \, p^{e_1}$.  
Comparing the $(2, 3)$-th entries of both sides of the equality $(\ast)$, we have 
$z^{\ell_2} \, \lambda \, \alpha_1(t) = \lambda \, \alpha_1(z^2 \, t) \, z^{\ell_3}$, which implies 
$\ell_2 - \ell_3 = 2 \, p^{e_1}$.   
Comparing the $(1, 3)$-th entries of both sides of the equality $(\ast)$, we have $ z^{4 \, p^{e_1}} \alpha_2(t)  = \alpha_2(z^2 \, t)$. 
So, $\alpha_2(t)$ is a $p$-monomial.  
We can express $\alpha_2(t)$ as $\alpha_2(t) = c_2 \, t^{p^{e_2}}$ ($c_2 \in k$, $e_2 \geq 0$).  
So, $c_2 \,  z^{4 \, p^{e_1}}  \, t^{p^{e_2}}  =  c_2 \, z^{2 \, p^{e_2}}  \, t^{p^{e_2}}$, which implies $c_2 = 0$ (since $p \geq 3$). 
Thus $u$ has the desired form. 

Since $\ell_1 + \ell_2 + \ell_3 = 0$, we have $(\ell_2 + 2 \, p^{e_1}) + \ell_2 + (\ell_2 - 2 \, p^{e_1}) = 0$ and thereby have 
$\ell_1 = 2 \, p^{e_1}$, $\ell_2 = 0$, $\ell_3 = - 2 \, p^{e_1}$.

Now, $u_\varphi$ and $(\ell_1, \ell_2, \ell_3)$ have the desired forms. 
Thus we can express $\varphi$ as the desired form. 
In fact, 
\begin{align*}
\varphi
\left(
\begin{array}{c c}
 a & b \\
 0 & d
\end{array}
\right) 
 & = 
\varphi
\left( \; 
\left(
\begin{array}{c c}
 a & 0 \\
 0 & d
\end{array}
\right) 
\left(
\begin{array}{c c}
 1 & \frac{b}{a} \\
 0 & 1
\end{array}
\right) 
\; \right)  \\
 & =  
\left(
\begin{array}{c c c}
 a^{2 \, p^{e_1}}  &  0  &  0  \\
 0  &  1  &  0  \\
 0 & 0 & d^{ 2 \, p^{e_1}} 
\end{array}  
\right) 
\left(
\begin{array}{c c c}
 1 & c_1 \, ( \frac{b}{a} )^{p^{e_1}} & \frac{1}{2} \, \lambda \, c_1^2 \,  ( \frac{b}{a} )^{2 \, p^{e_1}} \\
 0 & 1 & \lambda \, c_1 \, ( \frac{b}{a} )^{p^{e_1}}  \\
 0 & 0 & 1
\end{array}
\right) \\
 & =
\left(
\begin{array}{c c c}
 a^{2 \, p^{e_1}}  & c_1 \, a^{p^{e_1}} \, b^{p^{e_1}}  & \frac{1}{2} \, \lambda \, c_1^2 \, b^{2 \, p^{e_1}} \\
 0  &  1  &   \lambda \, c_1  \,  b^{p^{e_1}} \,  d^{p^{e_1}} \\
 0 & 0 & d^{2 \, p^{e_1}} 
\end{array}  
\right) . 
\end{align*}

\paragraph{(1.2)}

We have 
\begin{align}
& 
\left( 
\begin{array}{ccc}
 z^{\ell_1} & 0 & 0 \\
 0 & z^{\ell_2} & 0 \\
 0 & 0 & z^{\ell_3} 
\end{array}
\right) 
\left( 
\begin{array}{ccc}
 1 & \alpha_1(t) & \alpha_2(t)  \\
 0 & 1 & 0 \\
 0 & 0 & 1 
\end{array}
\right) 
 = 
\left( 
\begin{array}{ccc}
 1 & \alpha_1(z^2 \, t) &  \alpha_2(z^2 \, t)  \\
 0 & 1 & 0 \\
 0 & 0 & 1 
\end{array}
\right) 
\left( 
\begin{array}{ccc}
 z^{\ell_1} & 0 & 0 \\
 0 & z^{\ell_2} & 0 \\
 0 & 0 & z^{\ell_3} 
\end{array}
\right) . 
\tag{$\ast$}
\end{align} 
Comparing the $(1, 2)$-th entries of both sides of the equality $(\ast)$, we have $z^{\ell_1} \, \alpha_1(t) = \alpha_1(z^2 \, t) \, z^{\ell_2}$. 
So, $\alpha_1(t) = c_1 \, t^{p^{e_1}}$ for some $c_1 \in k$ and $e_1 \geq 0$. 
If $c_1 \ne 0$, then $\ell_1 - \ell_2 = 2 \, p^{e_1}$. 
Comparing the $(1, 3)$-th entries of both sides of the above equality $(\ast)$, 
we have $z^{\ell_1} \, \alpha_2(t) = \alpha_2(z^2 \, t) \, z^{\ell_3}$. 
So, $\alpha_2(t) = c_2 \, t^{p^{e_2}}$ for some $c_2 \in k$ and $e_2 \geq 0$. 
If $c_2 \ne 0$, then $\ell_1 - \ell_3 = 2 \, p^{e_2}$. 
Thus $u_\varphi$ has the desired form. 

Since $\ell_1 + \ell_2 + \ell_3 = 0$ and $\ell_1 > \ell_2 > \ell_3$, we have $\ell_1 > 0$ and $\ell_3 < 0$. 

We express $\varphi$ by separating the following four cases (1.2.a), (1.2.b), (1.2.c), (1.2.d): 
\begin{enumerate}[label = {\rm (1.2.\alph*)}, leftmargin = *]
\item $c_1 = 0$ and $c_2 = 0$. 

\item $c_1 \ne 0$ and $c_2 = 0$.  

\item $c_1 = 0$ and $c_2 \ne 0$.  

\item $c_1 \ne 0$ and $c_2 \ne 0$. 
\end{enumerate}

In the case (1.2.a), $\varphi$ already has the desired forms.

In the case (1.2.b), 
\begin{align*}
\varphi
\left(
\begin{array}{c c}
 a & b \\
 0 & d 
\end{array}
\right) 
 & =   
\left(
\begin{array}{c c c}
 a^{\ell_1} & 0 & 0 \\
 0 & a^{\ell_2} & 0 \\
 0 & 0 & a^{\ell_3} 
\end{array} 
\right)
\left(
\begin{array}{c c c}
 1 & c_1 \, (\frac{b}{a})^{p^{e_1}} & 0 \\
 0 &1 & 0 \\
 0 & 0 & 1
\end{array} 
\right)
= 
\left(
\begin{array}{c c c}
 a^{\ell_1} & c_1 \, a^{\ell_1 - p^{e_1}} \, b^{p^{e_1}}  & 0 \\
 0 & a^{\ell_2} & 0 \\
 0 & 0 & a^{\ell_3} 
\end{array} 
\right) \\
 & = 
\left(
\begin{array}{c c c}
 a^{\ell_1} & c_1 \, a^{\ell_2 + p^{e_1}} \, b^{p^{e_1}}  & 0 \\
 0 & a^{\ell_2} & 0 \\
 0 & 0 & a^{\ell_3} 
\end{array} 
\right)  \\
& =
\left\{
\begin{array}{l @{\qquad} l} 
\left(
\begin{array}{c c c}
 a^{\ell_1} & c_1 \, a^{\ell_2 + p^{e_1}} \, b^{p^{e_1}}  & 0 \\
 0 & a^{\ell_2} & 0 \\
 0 & 0 & d^{ - \ell_3} 
\end{array} 
\right) 
&  ( \, \ell_2 \geq 0 \, ) ,  \\
\left(
\begin{array}{c c c}
 a^{\ell_1} & c_1 \, a^{p^{e_1}} \, b^{p^{e_1}} \, d^{- \ell_2}  & 0 \\
 0 & d^{-\ell_2} & 0 \\
 0 & 0 & d^{ - \ell_3} 
\end{array}
\right) 
& (\,  \ell_2 < 0 \, ) .  
\end{array} 
\right. 
\end{align*}

In the case (1.2.c), 
\begin{align*}
\varphi
\left(
\begin{array}{c c}
 a & b \\
 0 & d 
\end{array}
\right) 
 & =   
\left(
\begin{array}{c c c}
 a^{\ell_1} & 0 & 0 \\
 0 & a^{\ell_2} & 0 \\
 0 & 0 & a^{\ell_3} 
\end{array} 
\right)
\left(
\begin{array}{c c c}
 1 & 0 & c_2 \, (\frac{b}{a})^{p^{e_2}}  \\
 0 &1 & 0 \\
 0 & 0 & 1
\end{array} 
\right)
= 
\left(
\begin{array}{c c c}
 a^{\ell_1} & 0 &  c_2 \, a^{\ell_1 - p^{e_2}} \, b^{p^{e_2}}   \\
 0 & a^{\ell_2} & 0 \\
 0 & 0 & a^{\ell_3} 
\end{array} 
\right) \\
 & = 
\left(
\begin{array}{c c c}
 a^{\ell_1} & 0 & c_2 \, a^{\ell_3 + p^{e_2}} \, b^{p^{e_2}}  \\
 0 & a^{\ell_2} & 0 \\
 0 & 0 & a^{\ell_3} 
\end{array} 
\right)  \\
& =
\left\{
\begin{array}{l @{\qquad} l} 
\left(
\begin{array}{c c c}
 a^{\ell_1} & 0 &  c_2 \, a^{p^{e_2}} \, b^{p^{e_2}} \, d^{- \ell_3} \\
 0 & a^{\ell_2} & 0 \\
 0 & 0 &  d^{ - \ell_3} 
\end{array} 
\right) 
&  ( \, \ell_2 \geq 0 \, ) ,  \\
\left(
\begin{array}{c c c}
 a^{\ell_1} & 0 & c_2 \, a^{p^{e_2}} \, b^{p^{e_2}} \, d^{- \ell_3}   \\
 0 & d^{-\ell_2} & 0 \\
 0 & 0 & d^{ - \ell_3} 
\end{array}
\right) 
& (\,  \ell_2 < 0 \, ) .  
\end{array} 
\right. 
\end{align*}

In the case (1.2.d), since $\ell_1 + \ell_2 + \ell_3 = 0$, we have 
$\ell_1 = (2 \, p^{e_1} + 2 \, p^{e_2})/3$,  
$\ell_2 = (-4 \, p^{e_1} + 2 \, p^{e_2})/3$, 
$\ell_3 = ( 2 \, p^{e_1} - 4 \, p^{e_2} ) / 3$. 
Since $\ell_1 - \ell_3 > \ell_1 - \ell_2$, we have $e_2 > e_1$. 
Since $p^{e_2} \geq 2 \, p^{e_1}$, we have $\ell_2 \geq 0$. 
So, we can express $\varphi$ with the desired form.

\paragraph{(1.3)}

We have 
\begin{align*}
& 
\left( 
\begin{array}{ccc}
 z^{\ell_1} & 0 & 0 \\
 0 & z^{\ell_2} & 0 \\
 0 & 0 & z^{\ell_3} 
\end{array}
\right) 
\left( 
\begin{array}{ccc}
 1 & 0 & \alpha_2(t)  \\
 0 & 1 & \alpha_1(t) \\
 0 & 0 & 1 
\end{array}
\right) 
 = 
\left( 
\begin{array}{ccc}
 1 & 0 &  \alpha_2(z^2 \, t)  \\
 0 & 1 & \alpha_1(z^2 \, t) \\
 0 & 0 & 1 
\end{array}
\right) 
\left( 
\begin{array}{ccc}
 z^{\ell_1} & 0 & 0 \\
 0 & z^{\ell_2} & 0 \\
 0 & 0 & z^{\ell_3} 
\end{array}
\right) . 
\tag{$\ast$}
\end{align*} 
Comparing the $(1, 3)$-th entries of both sides of the above equality $(\ast)$, 
we have $z^{\ell_1} \, \alpha_2(t) = \alpha_2(z^2 \, t) \, z^{\ell_3}$. 
So, $\alpha_2(t) = c_2 \, t^{p^{e_2}}$ for some $c_2 \in k$ and $e_2 \geq 0$. 
If $c_2 \ne 0$, then $\ell_1 - \ell_3 = 2 \, p^{e_2}$. 
Comparing the $(2,3)$-th entries of both sides of the above equality $(\ast)$, 
we have $z^{\ell_2} \, \alpha_1(t) = \alpha_1(z^2 \, t) \, z^{\ell_3}$. 
So, $\alpha_1(t) = c_1 \, t^{p^{e_1}}$ for some $c_1 \in k$ and $e_1 \geq 0$. 
If $c_1 \ne 0$, then $\ell_2 - \ell_3 = 2 \, p^{e_1}$. 
Thus $u_\varphi$ has the desired form.

Since $\ell_1 + \ell_2 + \ell_3 = 0$ and $\ell_1 > \ell_2 > \ell_3$, we have $\ell_1 > 0$ and $\ell_3 < 0$. 

We express $\varphi$ by separating the following four cases (1.3.a), (1.3.b), (1.3.c), (1.3.d): 
\begin{enumerate}[label = {\rm (1.3.\alph*)}, leftmargin = *]
\item $c_1 = 0$ and $c_2 = 0$.  

\item $c_1 \ne 0$ and $c_2 = 0$.  

\item $c_1 = 0$ and $c_2 \ne 0$. 

\item $c_1 \ne 0$ and $c_2 \ne 0$. 
\end{enumerate}

In the case (1.3.a), $\varphi$ already has the desired forms.

In the case (1.3.b), 
\begin{align*}
\varphi
\left(
\begin{array}{c c}
 a & b \\
 0 & d 
\end{array}
\right) 
 & =   
\left(
\begin{array}{c c c}
 a^{\ell_1} & 0 & 0 \\
 0 & a^{\ell_2} & 0 \\
 0 & 0 & a^{\ell_3} 
\end{array} 
\right)
\left(
\begin{array}{c c c}
 1 & 0 & 0 \\
 0 & 1 & c_1 \, (\frac{b}{a})^{p^{e_1}} \\
 0 & 0 & 1
\end{array} 
\right)
= 
\left(
\begin{array}{c c c}
 a^{\ell_1} & 0  & 0 \\
 0 & a^{\ell_2} & c_1 \, a^{\ell_2 - p^{e_1}} \, b^{p^{e_1}} \\
 0 & 0 & a^{\ell_3} 
\end{array} 
\right) \\
 & = 
\left(
\begin{array}{c c c}
 a^{\ell_1} & 0  & 0 \\
 0 & a^{\ell_2} & c_1 \, a^{\ell_3 + p^{e_1}} \, b^{p^{e_1}} \\
 0 & 0 & a^{\ell_3} 
\end{array} 
\right)  \\
& =
\left\{
\begin{array}{l @{\qquad} l} 
\left(
\begin{array}{c c c}
 a^{\ell_1} & 0  & 0 \\
 0 & a^{\ell_2} & c_1 \, a^{ p^{e_1}} \, b^{p^{e_1}} \, d^{- \ell_3}  \\
 0 & 0 & d^{ - \ell_3} 
\end{array} 
\right) 
&  ( \, \ell_2 \geq 0 \, ) ,  \\
\left(
\begin{array}{c c c}
 a^{\ell_1} & 0  & 0 \\
 0 & d^{-\ell_2} & c_1 \, a^{ p^{e_1}} \, b^{p^{e_1}} \, d^{- \ell_3} \\
 0 & 0 & d^{ - \ell_3} 
\end{array}
\right) 
& (\,  \ell_2 < 0 \, ) .  
\end{array} 
\right. 
\end{align*}

In the case (1.3.c), 
\begin{align*}
\varphi
\left(
\begin{array}{c c}
 a & b \\
 0 & d 
\end{array}
\right) 
 & =   
\left(
\begin{array}{c c c}
 a^{\ell_1} & 0 & 0 \\
 0 & a^{\ell_2} & 0 \\
 0 & 0 & a^{\ell_3} 
\end{array} 
\right)
\left(
\begin{array}{c c c}
 1 & 0 & c_2 \, (\frac{b}{a})^{p^{e_2}} \\
 0 & 1 & 0 \\
 0 & 0 & 1
\end{array} 
\right)
= 
\left(
\begin{array}{c c c}
 a^{\ell_1} & 0  & c_2 \, a^{\ell_1 - p^{e_2}} \, b^{p^{e_2}} \\
 0 & a^{\ell_2} & 0 \\
 0 & 0 & a^{\ell_3} 
\end{array} 
\right) \\
 & = 
\left(
\begin{array}{c c c}
 a^{\ell_1} & 0  & c_2 \, a^{\ell_3 + p^{e_2}} \, b^{p^{e_2}}  \\
 0 & a^{\ell_2} & 0 \\
 0 & 0 & a^{\ell_3} 
\end{array} 
\right)  \\
& =
\left\{
\begin{array}{l @{\qquad} l} 
\left(
\begin{array}{c c c}
 a^{\ell_1} & 0  & c_2 \, a^{ p^{e_2}} \, b^{p^{e_2}}  \, d^{-\ell_3} \\
 0 & a^{\ell_2} & 0  \\
 0 & 0 & d^{ - \ell_3} 
\end{array} 
\right) 
&  ( \, \ell_2 \geq 0 \, ) ,  \\
\left(
\begin{array}{c c c}
 a^{\ell_1} & 0  & c_2 \, a^{p^{e_2}} \, b^{p^{e_2}} \, d^{-\ell_3}   \\
 0 & d^{-\ell_2} & 0 \\
 0 & 0 & d^{ - \ell_3} 
\end{array}
\right) 
& (\,  \ell_2 < 0 \, ) .  
\end{array} 
\right. 
\end{align*}

In the case (1.3.d), 
since $\ell_1 + \ell_2 + \ell_3 = 0$, we have 
$\ell_1 = ( - 2 \, p^{e_1} + 4 \, p^{e_2} )/ 3$, 
$\ell_2 = ( 4 \, p^{e_1} - 2 \, p^{e_2} ) / 3$, 
$\ell_3 = ( - 2 \, p^{e_1} - 2 \, p^{e_2} ) / 3$. 
Since $\ell_1 - \ell_3 > \ell_2 - \ell_3$, we have $e_2 > e_1$.  
Since $p^{e_2}  \geq 2 \, p^{e_1}$, we have $\ell_2 \leq 0$. 
So, we can express $\varphi$ with the desired form.

\subsubsection{Proof of assertion (2) of Lemma 2.1}

We know from Lemma 1.5 that $u \in \cA_{(2, 1)}$. 
So, we have 
\begin{align*}
& 
\left( 
\begin{array}{ccc}
 z^{\ell_1} & 0 & 0 \\
 0 & z^{\ell_2} & 0 \\
 0 & 0 & z^{\ell_3} 
\end{array}
\right) 
\left( 
\begin{array}{ccc}
 1 & 0 & \alpha_2(t)  \\
 0 & 1 & \alpha_1(t) \\
 0 & 0 & 1 
\end{array}
\right) 
 = 
\left( 
\begin{array}{ccc}
 1 & 0 &  \alpha_2(z^2 \, t)  \\
 0 & 1 & \alpha_1(z^2 \, t) \\
 0 & 0 & 1 
\end{array}
\right) 
\left( 
\begin{array}{ccc}
 z^{\ell_1} & 0 & 0 \\
 0 & z^{\ell_2} & 0 \\
 0 & 0 & z^{\ell_3} 
\end{array}
\right) . 
\end{align*} 
As in the case (1.3), for each $i = 1, 2$, we can express $\alpha_i$ as 
$\alpha_i(t) = c_i \, t^{p^{e_i}}$ for some $c_i \in k$ and $e_i \geq 0$. 
Also, we have $\ell_1 - \ell_3 = 2 \, p^{e_2}$ provided that $c_2 \ne 0$,  
and $\ell_2 - \ell_3 = 2 \, p^{e_1}$ provided that $c_1 \ne 0$. 
Thus $u_\varphi$ has the desired form.

Since $\ell_1 + \ell_2 + \ell_3 = 0$ and $\ell_1 = \ell_2 > \ell_3$, we have $\ell_1 = \ell_2 > 0$ and $\ell_3 < 0$. 

We express $\varphi$ by separating the following four cases (2.a), (2.b), (2.c), (2.d): 
\begin{enumerate}[label = {\rm (2.\alph*)}, leftmargin = *]
\item $c_1 = 0$ and $c_2 = 0$. 

\item $c_1 \ne 0$ and $c_2 = 0$. 

\item $c_1 = 0$ and $c_2 \ne 0$. 

\item $c_1 \ne 0$ and $c_2 \ne 0$. 
\end{enumerate}

In the case (2.a), $\varphi$ already has the desired form.

In the cases (2.b) and (2.c), we can express $\varphi$ with the desired forms.

In the case (2.d), we have $ 2 \, p^{e_2} = \ell_1 - \ell_3 = \ell_2 - \ell_3 = 2 \, p^{e_1}$, which implies $e_1 = e_2$. 
Since $\ell_1 + \ell_2 + \ell_3 = 0$, we have $3 \, \ell_3 =  -4 \, p^{e_1}$, 
which implies $p = 3$ and $e_1 \geq 1$. 
Thus 
$\ell_3 = - 4 \, p^{e_1 - 1} = - 4 \cdot 3^{e_1 - 1}$ and 
$\ell_1 = \ell_3 + 2 \, p^{e_1} = - 4 \, p^{e_1 - 1} + 2 \, p^{e_1} = 2 \cdot 3^{e_1 - 1}$.  
So, we can express $\varphi$ with the desired form.

\subsubsection{Proof of assertion (3) of Lemma 2.1}

We know from Lemma 1.5 that $u_\varphi \in \cA_{(1, 2)}$. 
So, we have 
\begin{align*}
& 
\left( 
\begin{array}{ccc}
 z^{\ell_1} & 0 & 0 \\
 0 & z^{\ell_2} & 0 \\
 0 & 0 & z^{\ell_3} 
\end{array}
\right) 
\left( 
\begin{array}{ccc}
 1 & \alpha_1(t) & \alpha_2(t)  \\
 0 & 1 & 0 \\
 0 & 0 & 1 
\end{array}
\right) 
 = 
\left( 
\begin{array}{ccc}
 1 & \alpha_1(z^2 \, t) &  \alpha_2(z^2 \, t)  \\
 0 & 1 & 0 \\
 0 & 0 & 1 
\end{array}
\right) 
\left( 
\begin{array}{ccc}
 z^{\ell_1} & 0 & 0 \\
 0 & z^{\ell_2} & 0 \\
 0 & 0 & z^{\ell_3} 
\end{array}
\right) . 
\end{align*} 
As in the case (1.2), for each $i = 1, 2$, we can express $\alpha_i$ as 
$\alpha_i(t) = c_i \, t^{p^{e_i}}$ for some $c_i \in k$ and $e_i \geq 0$. 
Also, we have $\ell_1 - \ell_2 = 2 \, p^{e_1}$ provided that $c_1 \ne 0$, and $\ell_1 - \ell_3 = 2 \, p^{e_2}$ provided that $c_2 \ne 0$. 
Thus $u_\varphi$ has the desired form.

Since $\ell_1 + \ell_2 + \ell_3 = 0$ and $\ell_1 > \ell_2 = \ell_3$, we have $\ell_1 > 0$ and $\ell_2 = \ell_3 < 0$. 

We express $\varphi$ by separating the following four cases (3.a), (3.b), (3.c), (3.d): 
\begin{enumerate}[label = {\rm (3.\alph*)}, leftmargin = *]
\item $c_1 = 0$ and $c_2 = 0$. 

\item $c_1 \ne 0$ and $c_2 = 0$. 

\item $c_1 = 0$ and $c_2 \ne 0$. 

\item $c_1 \ne 0$ and $c_2 \ne 0$. 
\end{enumerate}

In the case (3.a), $\varphi$ already has the desired form.

In the cases (3.b) and (3.c), we can express $\varphi$ with the desired forms.

In the case where (3.d), we have $2 \, p^{e_1} = \ell_1 - \ell_2 = \ell_1 - \ell_3 = 2 \, p^{e_2}$, which implies $e_1 = e_2$. 
Since $\ell_1 + \ell_2 + \ell_3 = 0$, we have $3 \, \ell_1 =  4 \, p^{e_1} $, 
which implies $p  =3$ and $e_1 \geq 1$. 
Thus $\ell_1 = 4 \, p^{e_1 - 1}$, $\ell_2 = - 2 \, p^{e_1 - 1}$ and $\ell_3 =  -2 \, p^{e_1 - 1}$. 
So, we can express $\varphi$ with the desired form.

\subsubsection{Proof of assertion (4) of Lemma 2.1} 

We know from Lemma 1.5 that $u_\varphi$ is trivial. 
Since $\ell_1 + \ell_2 + \ell_3 = 0$ and $\ell_1 = \ell_2 = \ell_3$, 
we have $\ell_1 = \ell_2 = \ell_3 = 0$. Thus $\varphi$ is trivial.

\subsubsection{On the characteristic $p$ in (1.2.d) and (1.3.d)}

\begin{lem}
Let $p$ be a prime number. Let $e_1$ and $e_2$ be integers satisfying $e_2 > e_1 \geq 0$. 
Then $( 2 \, p^{e_1} + 2 \, p^{e_2} ) /3$ is an integer if and only if 
one of the following conditions {\rm (1)} and {\rm (2)} holds true: 
\begin{enumerate}[label = {\rm (\arabic*)}]
\item $p = 3$. 

\item $p \equiv -1 \pmod 3$ and $e_1 - e_2 \equiv 1 \pmod 2$. 
\end{enumerate} 
\end{lem}

\Proof The proof is straightforward. 
\QED

\subsection{A classification of representations of $\G_a \rtimes \G_m$ into $\SL(3, k)$}

We can define representations $\varphi^* : \G_a \rtimes \G_m \to \SL(3, k)$ 
by separating the following cases (1)$^*$, (2)$^*$, (3)$^*$, (4)$^*$ 
(These cases correspond to the cases (1), (2), (3), (4) in Lemma 2.1):

\begin{enumerate}[label = {\rm (\arabic*)$^*$}, leftmargin = * ]

\item Let $\ell_1$, $\ell_2$, $\ell_3$ be integers satisfying $\ell_1 + \ell_2 + \ell_3 = 0$, 
$\ell_1 > \ell_2 > \ell_3$ and $\ell_1 > 0 > \ell_3$. 

\begin{enumerate}[label = {\rm (\Roman*)$^*$}, leftmargin = *]  

\item If $p \geq 3$ and there exists an integer $e_1 \geq 0$ such that 
$(\ell_1, \ell_2, \ell_3) = (2 \, p^{e_1}, 0, - 2 \, p^{e_1})$, 
then we can define $\varphi^* : \G_a \rtimes \G_m \to \SL(3, k)$ where $p \geq 2$ as 
\begin{align*} 
\varphi^*
\left(
\begin{array}{c c}
 a & b \\
 0 & d 
\end{array}
\right) 
 & :=
\left(
\begin{array}{c c c}
 a^{\ell_1}  &  a^{p^{e_1}} \, b^{p^{e_1}}  & \frac{1}{2} \,  b^{2 \, p^{e_1}} \\
 0  & a^{\ell_2} &   b^{p^{e_1}} \,  d^{p^{e_1}} \\
 0 & 0 & a^{\ell_3}  
\end{array}  
\right) 
 =
\left(
\begin{array}{c c c}
 a^{2 \, p^{e_1}}  &  a^{p^{e_1}} \, b^{p^{e_1}}  & \frac{1}{2} \,  b^{2 \, p^{e_1}} \\
 0  &  1  &   b^{p^{e_1}} \,  d^{p^{e_1}} \\
 0 & 0 & d^{2 \, p^{e_1}} 
\end{array}  
\right) . 
\end{align*}

\item With no additional assumption on $\ell_1$, $\ell_2$, $\ell_3$, we can define $\varphi^* : \G_a \rtimes \G_m \to \SL(3, k)$ as 
\begin{fleqn}[36pt] 
\begin{align*}
\varphi^* 
\left(
\begin{array}{c c}
 a & b \\
 0 & d 
\end{array}
\right) 
& :=
\left\{
\begin{array}{l @{\qquad} l} 
\left(
\begin{array}{c c c}
 a^{\ell_1} & 0 & 0 \\
 0 & a^{\ell_2} & 0 \\
 0 & 0 & d^{ - \ell_3} 
\end{array} 
\right)
 &  (\, \ell_2 \geq 0 \, ) , \\  [24pt] 
\left(
\begin{array}{c c c}
 a^{\ell_1} & 0 & 0 \\
 0 & d^{- \ell_2} & 0 \\
 0 & 0 & d^{ - \ell_3} 
\end{array} 
\right)
 &  (\, \ell_2 < 0 \, ) . 
\end{array}
\right. 
\end{align*} 
\end{fleqn}

\item If there exists an integer $e_1 \geq 0$ such that $\ell_1 - \ell_2 = 2 \, p^{e_1}$, 
then we can define $\varphi^* : \G_a \rtimes \G_m \to \SL(3, k)$ where $p \geq 2$ as 
\begin{fleqn}[36pt] 
\begin{align*}
\varphi^*
\left(
\begin{array}{c c}
 a & b \\
 0 & d 
\end{array}
\right) 
& :=
\left\{
\begin{array}{l @{\qquad} l} 
\left(
\begin{array}{c c c}
 a^{\ell_1} & a^{\ell_2 + p^{e_1}} \, b^{p^{e_1}}  & 0 \\
 0 & a^{\ell_2} & 0 \\
 0 & 0 & d^{ - \ell_3} 
\end{array} 
\right) 
&  ( \, \ell_2 \geq 0 \, ) ,  \\  [24pt] 
\left(
\begin{array}{c c c}
 a^{\ell_1} &  a^{p^{e_1}} \, b^{p^{e_1}} \, d^{- \ell_2}  & 0 \\
 0 & d^{-\ell_2} & 0 \\
 0 & 0 & d^{ - \ell_3} 
\end{array}
\right) 
& (\,  \ell_2 < 0 \, ) .  
\end{array} 
\right. 
\end{align*}
\end{fleqn}

\item If there exists an integer $e_2 \geq 0$ such that $\ell_1 - \ell_3 = 2 \, p^{e_2}$,  
then we can define $\varphi^* : \G_a \rtimes \G_m \to \SL(3, k)$ where $p \geq 2$ as 
\begin{fleqn}[36pt] 
\begin{align*}
\varphi^*
\left(
\begin{array}{c c}
 a & b \\
 0 & d 
\end{array}
\right) 
& :=
\left\{
\begin{array}{l @{\qquad} l} 
\left(
\begin{array}{c c c}
 a^{\ell_1} & 0 &  a^{p^{e_2}} \, b^{p^{e_2}} \, d^{- \ell_3} \\
 0 & a^{\ell_2} & 0 \\
 0 & 0 &  d^{ - \ell_3} 
\end{array} 
\right) 
&  ( \, \ell_2 \geq 0 \, ) ,  \\  [24pt] 
\left(
\begin{array}{c c c}
 a^{\ell_1} & 0 &  a^{p^{e_2}} \, b^{p^{e_2}} \, d^{- \ell_3}   \\
 0 & d^{-\ell_2} & 0 \\
 0 & 0 & d^{ - \ell_3} 
\end{array}
\right) 
& (\,  \ell_2 < 0 \, ) .  
\end{array} 
\right. 
\end{align*}
\end{fleqn}

\item If there exist integers $e_1, e_2$ such that $e_2 > e_1 \geq 0$ and 
\[
(\ell_1, \; \ell_2, \; \ell_3) 
 =
 \left( 
  \frac{2 \, p^{e_1} + 2 \, p^{e_2}}{3} , \quad
  \frac{-4 \, p^{e_1} + 2 \, p^{e_2}}{3} , \quad 
 \frac{2 \, p^{e_1} - 4 \, p^{e_2}}{3} 
\right) ,  
\]
then we can define $\varphi^* : \G_a \rtimes \G_m \to \SL(3, k)$ as 
\begin{fleqn}[36pt] 
\begin{align*}
\varphi^*
\left(
\begin{array}{c c}
 a & b \\
 0 & d 
\end{array}
\right) 
  & := 
\left(
\begin{array}{c c c}
 a^{\ell_1}  &  a^{\ell_1} \, b^{p^{e_1}} \, d^{p^{e_1}}  &  a^{\ell_1} \, b^{p^{e_2}} \, d^{p^{e_2}}   \\
 0  &  a^{\ell_2}  &  0  \\
 0 & 0 & d^{- \ell_3} 
\end{array}  
\right) .  
\end{align*}
\end{fleqn} 
Clearly, $\ell_1 > 0$, $\ell_2 \geq 0$, $\ell_3 < 0$. And $p = 3$ or $1 + p^{e_2 - e_1}$ is a multiple of $3$.

\item 
Let $e_1 \geq 0$ and assume $\ell_2 - \ell_3 = 2 \, p^{e_1}$.  
Then we can define $\varphi^* : \G_a \rtimes \G_m \to \SL(3, k)$ as 
\begin{fleqn}[36pt] 
\begin{align*}
\varphi^*
\left(
\begin{array}{c c}
 a & b \\
 0 & d 
\end{array}
\right) 
& :=
\left\{
\begin{array}{l @{\qquad} l} 
\left(
\begin{array}{c c c}
 a^{\ell_1} & 0  & 0 \\
 0 & a^{\ell_2} & a^{ p^{e_1}} \, b^{p^{e_1}} \, d^{- \ell_3}  \\
 0 & 0 & d^{ - \ell_3} 
\end{array} 
\right) 
&  ( \, \ell_2 \geq 0 \, ) ,  \\  [24pt] 
\left(
\begin{array}{c c c}
 a^{\ell_1} & 0  & 0 \\
 0 & d^{-\ell_2} &  a^{ p^{e_1}} \, b^{p^{e_1}} \, d^{- \ell_3} \\
 0 & 0 & d^{ - \ell_3} 
\end{array}
\right) 
& (\,  \ell_2 < 0 \, ) .  
\end{array} 
\right. 
\end{align*}
\end{fleqn}

\item Let $e_1, e_2$ be integers satisfying $e_2 > e_1 \geq 0$. 
Assume 
\[
(\ell_1, \; \ell_2, \; \ell_3) 
 = 
\left( 
  \frac{ - 2 \, p^{e_1} + 4 \, p^{e_2}}{3} , \quad 
  \frac{4 \, p^{e_1} - 2 \, p^{e_2}}{3} , \quad 
  \frac{ - 2 \, p^{e_1} - 2 \, p^{e_2}}{3}  
\right) . 
\]
Clearly, $\ell_1 > 0$, $\ell_2 \leq 0$ and $\ell_3 < 0$. 
Then we can define $\varphi^* : \G_a \rtimes \G_m \to \SL(3, k)$ as 
\begin{fleqn}[36pt] 
\begin{align*}
\varphi^* 
\left(
\begin{array}{c c}
 a & b \\
 0 & d 
\end{array}
\right) 
 & :=
\left(
\begin{array}{c c c}
 a^{\ell_1}  &  0  &  a^{p^{e_2}} \, b^{p^{e_2}} \, d^{- \ell_3}   \\
 0  &  d^{- \ell_2}  &  a^{p^{e_1}} \, b^{p^{e_1}} \,  d^{- \ell_3}  \\
 0 & 0 & d^{ - \ell_3} 
\end{array}  
\right) . 
\end{align*} 
\end{fleqn} 
\end{enumerate}

\item Let $\ell_1$, $\ell_2$, $\ell_3$ be integers satisfying $\ell_1 + \ell_2 + \ell_3 = 0$ and $\ell_1 = \ell_2 > 0 > \ell_3$. 

\begin{enumerate}[label = {\rm (\Roman*)$^*$}, leftmargin = *, start = 8 ]
\item 
With no additional assumption on $\ell_1$, $\ell_2$, $\ell_3$, we can define $\varphi^* : \G_a \rtimes \G_m \to \SL(3, k)$ where $p \geq 2$ as 
\begin{fleqn}[36pt] 
\begin{align*}
\varphi^* 
\left(
\begin{array}{c c}
 a & b \\
 0 & d 
\end{array}
\right) 
 :=
\left(
\begin{array}{c c c}
 a^{\ell_1}  &  0  &  0  \\
 0  &  a^{\ell_1}  &  0  \\
 0 & 0 & d^{- \ell_3} 
\end{array}  
\right) . 
\end{align*}
\end{fleqn}

\item 
If there exists an integer $e_1 \geq 0$ such that $\ell_1 - \ell_3 = 2 \, p^{e_1}$, 
then we can define $\varphi^* : \G_a \rtimes \G_m \to \SL(3, k)$ where $p \geq 2$ as 
\begin{fleqn}[36pt] 
\begin{align*} 
\varphi^* 
\left(
\begin{array}{c c}
 a & b \\
 0 & d 
\end{array}
\right) 
  & := 
\left(
\begin{array}{c c c}
 a^{\ell_1}  &  0  &  0  \\
 0  &  a^{\ell_1}  &  a^{\ell_1} \, b^{p^{e_1}} \, d^{ p^{e_1} }  \\
 0 & 0 & d^{- \ell_3} 
\end{array}  
\right) . 
\end{align*}
\end{fleqn}

\end{enumerate}

\item Let $\ell_1$, $\ell_2$, $\ell_3$ be integers satisfying $\ell_1 + \ell_2 + \ell_3 = 0$ and $\ell_1 > 0 > \ell_2 = \ell_3$.

\begin{enumerate}[label = {\rm (\Roman*)$^*$}, leftmargin = * , start = 10]
\item 
With no additional assumption on $\ell_1$, $\ell_2$, $\ell_3$, we can define $\varphi^* : \G_a \rtimes \G_m \to \SL(3, k)$ as 
\begin{fleqn}[36pt]
\begin{align*}
\varphi^* 
\left(
\begin{array}{c c}
 a & b \\
 0 & d 
\end{array}
\right) 
 :=
\left(
\begin{array}{c c c}
 a^{\ell_1}  &  0  &  0  \\
 0  &  d^{- \ell_2}  &  0  \\
 0 & 0 & d^{- \ell_2} 
\end{array}  
\right) . 
\end{align*}
\end{fleqn}

\item 
If there exists an integer $e_1 \geq 0$ such that $\ell_1 - \ell_2 = 2 \, p^{e_1}$, 
then we can define $\varphi^* : \G_a \rtimes \G_m \to \SL(3, k)$ where $p \geq 2$ as 
\begin{fleqn}[36pt]
\begin{align*} 
\varphi^* 
\left(
\begin{array}{c c}
 a & b \\
 0 & d 
\end{array}
\right) 
 & := 
\left(
\begin{array}{c c c}
 a^{\ell_1}  &   a^{\ell_1} \,  b^{p^{e_1}} \, d^{ p^{e_1} }      &  0  \\
 0  &  d^{ - \ell_2}  & 0\\
 0 & 0 & d^{- \ell_2} 
\end{array}  
\right)   . 
\end{align*}
\end{fleqn}

\end{enumerate}

\item Let $\ell_1$, $\ell_2$, $\ell_3$ be integers satisfying $\ell_1 = \ell_2 = \ell_3 = 0$.

\begin{enumerate}[label = {\rm (\Roman*)$^*$}, leftmargin = *, start = 12]
\item We can define $\varphi^* :\G_a \rtimes \G_m \to \SL(3, k)$ where $p \geq 2$ as the trivial representation, i.e.,
\begin{fleqn}[36pt]
\begin{align*}
\varphi^* 
\left(
\begin{array}{c c}
 a & b \\
 0 & d 
\end{array}
\right)  
 & := 
\left(
\begin{array}{c c c}
 a^{\ell_1} & 0 & 0  \\
 0  & a^{\ell_2} & 0 \\
 0 & 0 & a^{\ell_3}
\end{array}
\right) 
= 
\left(
\begin{array}{c c c}
 1 & 0 & 0  \\
 0 & 1 & 0 \\
 0 & 0 & 1
\end{array}
\right) . 
\end{align*}
\end{fleqn} 
\end{enumerate} 

\end{enumerate} 

Clearly, any $\varphi^*$ with one of the forms (I)$^*$ -- (XII)$^*$ satisfies 
\[
\varphi^*
\left(
\begin{array}{c c}
 z & 0 \\
 0 & z^{-1}
\end{array}
\right)
 = 
\left(
\begin{array}{c c c}
 z^{\ell_1} & 0 & 0 \\
 0 & z^{\ell_2} & 0 \\
 0 & 0 & z^{\ell_3} 
\end{array}
\right) . 
\]

Let
\begin{align*}
\Lambda^{(1)^*} & := 
\bigl\{ \, 
 {\rm (I)^*}, \;  {\rm (II)^*}, \;  {\rm (III)^*}, \;  {\rm (IV)^*}, \;  {\rm (V)^*},  \; {\rm (VI)^*}, \;  {\rm (VII)^*} 
\, \bigr\} , \\
\Lambda^{(2)^*} & := 
\bigl\{ \, 
 {\rm (VIII)^*}, \;  {\rm (IX)^*}
\, \bigr\}, \\
\Lambda^{(3)^*} & := 
\bigl\{ \, 
 {\rm (X)^*}, \; {\rm (XI)^*}
\, \bigr\} , \\ 
\Lambda^{(4)^*} & := 
\bigl\{ \, 
 {\rm (XII)^*} 
\, \bigr\} , \\ 
\Lambda^* & := \Lambda^{(1)^*} \cup \Lambda^{(2)^*} \cup \Lambda^{(3)^*} \cup \Lambda^{(4)^*} . 
\end{align*}

For any $\lambda \in \Lambda^*$, we can define a set $R^*(\lambda)$ as 
\begin{align*}
 R^*(\lambda) := \{ \, \varphi : \G_a \rtimes \G_m \to \SL(3, k) \mid \text{$\varphi$ is a representation with the form $\lambda$} \, \} . 
\end{align*}

\begin{thm}
The following assertions {\rm (1)} and {\rm (2)} hold true: 
\begin{enumerate}[label = {\rm (\arabic*)}]
\item Any representation $\varphi : \G_a \rtimes \G_m \to \SL(3, k)$ 
is equivalent to a representation $\varphi^* : \G_a \rtimes \G_m \to \SL(3, k)$ with one of the forms $\lambda$ 
of $\Lambda^*$.

\item Let $\phi \in R^*(\lambda)$ and $\psi \in R^*(\mu)$, where $\lambda, \mu \in \Lambda^*$. 
Write $h_\phi(z) = (z^{\ell_1},  z^{\ell_2},  z^{\ell_3})$ and $h_{\psi}(z) = (z^{m_1},  z^{m_2},  z^{m_3})$.   
\begin{enumerate}[label = {\rm (2.\arabic*)}]
\item 
If $\lambda \ne \mu$, then $\phi$ and $\psi$ are not equivalent. 

\item If $\lambda = \mu$, then the following conditions {\rm (i)}, {\rm (ii)}, {\rm (iii)}, {\rm (iv)} are equivalent: 
\begin{enumerate}[label = {\rm (\roman*)}]
\item $\phi$ and $\psi$ are equivalent. 

\item $h_\phi$ and $h_{\psi}$ are equivalent. 

\item $(\ell_1, \, \ell_2, \, \ell_3) = (m_1, \, m_2, \, m_3)$. 

\item $\phi = \psi$. 
\end{enumerate} 
\end{enumerate}
\end{enumerate} 
\end{thm}

\begin{cor}
There exists a one-to-one correspondence between the set of all representations of $\G_a \rtimes \G_m$ 
into $\SL(3, k)$ up to equivalence and the set $\coprod_{\lambda \in \Lambda^*} R^*(\lambda)$, i.e., 
\[
\{\, \varphi : \G_a \rtimes \G_m \to \SL(3, k) \mid \text{$\varphi$ is a representation}\, \} / \sim 
 \; \cong  \: 
 \coprod_{\lambda \in \Lambda^*} R^*(\lambda) . 
\]
\end{cor}

\Proof 
Let $S := \{\, \varphi : \G_a \rtimes \G_m \to \SL(3, k) \mid \text{$\varphi$ is a representation}\, \} / \sim$ and let 
$T := \coprod_{\lambda \in \Lambda^*} R^*(\lambda)$. 

We can define a map $f : S \to T$, as follows: 
Let $\Phi$ be an element of $S$. 
There exists a representation $\varphi : \G_a \rtimes \G_m \to \SL(3, k)$ so that 
$\Phi$ is the equivalence class of $\varphi$, i.e., $\Phi = [\varphi]$.  
We know from assertion (1) of Theorem 2.3 that there exists a representation $\varphi^* : \G_a \rtimes \G_m \to \SL(3, k)$ 
so that $\varphi$ is equivalent to $\varphi^*$ and $\varphi^* \in R^*(\lambda)$ for some $\lambda \in \Lambda^*$. 
So, $\Phi = [\varphi^*]$. 
Such a $\lambda$ is uniquely determined by $\Phi$ (see assertion (2.1) of Theorem 2.3) and 
$\varphi^*$ is also uniquely determined by $\Phi$ (see assertion (2.2) of Theorem 2.3). 
Thus we can define $f : S \to T$ as $f(\Phi) : = \varphi^*$. 

Conversely, we can define a map $g : T \to S$ as $g(\varphi^*) := [\varphi^*]$.  

We can prove $g \circ f = \id_S$ and $f \circ g = \id_T$. 
This completes the proof. 
\QED

\subsubsection{Proof of assertion (1) of Theorem 2.3 } 

Using Lemma 1.1, we may assume from the first that $h_\varphi$ has the form 
\[
h_\varphi(z)
 = 
\diag(\,  z^{\ell_1}, \, z^{\ell_2}, \, z^{\ell_3} \,) 
\qquad 
(\
 \ell_1 \geq \ell_2 \geq \ell_3 
\,) . 
\]
We know from Lemma 2.1 that $\varphi$ has one of the forms (1.1), (1.2.a), (1.2.b), (1.2.c), (1.2.d), (1.3.a), (1.3.b), (1.3.c), (1.3.d), 
(2.a), (2.b), (2.c), (2.d), (3.a), (3.b), (3.c), (3.d), (4). 
For any $(\varphi, P)$ in each line of the following table, 
we can calculate the form of $\varphi^* : \G_a \rtimes \G_m \to \SL(3, k)$ defined by 
\[
\varphi^* \left( 
\begin{array}{c c}
 a & b \\
 c & d 
\end{array}
\right) : = P^{-1} \, \varphi 
\left( 
\begin{array}{c c}
 a & b \\
 c & d 
\end{array}
\right) \, P . 
\]
\begin{longtable}{| @{\hspace{12pt}}   l @{\hspace{24pt}} |  @{\hspace{12pt}}  l   @{\hspace{24pt}}   |   @{\hspace{12pt}}  l  @{\hspace{24pt}}  | }
\hline 
 $\varphi$
 & $P$
 & $\varphi^*$ \\
\hline
(1.1)
 & 
$
\left(
\begin{array}{c c c}
 1 & 0 & 0 \\
 0 & 1 / c_1 & 0 \\
 0 & 0 & 1 / (\lambda \, c_1^2 )
\end{array}
\right)
$
 & (I)$^*$ \\ 
\hline 
(1.2.a)
 & 
$
\left(
\begin{array}{c c c}
 1 & 0 & 0 \\
 0 & 1 & 0 \\
 0 & 0 & 1 
\end{array}
\right)
$
 & (II)$^*$ \\ 
\hline 
(1.2.b)
 & 
$
\left(
\begin{array}{c c c}
 1 & 0 & 0 \\
 0 & 1 / c_1 & 0 \\
 0 & 0 & 1 
\end{array}
\right)
$
 & (III)$^*$  \\ 
\hline 
(1.2.c)
& 
$
\left(
\begin{array}{c c c}
 1 & 0 & 0 \\
 0 & 1 & 0 \\
 0 & 0 & 1 / c_2 
\end{array}
\right)
$
 & (IV)$^*$ \\ 
\hline 
(1.2.d)
 & 
$
\left(
\begin{array}{c c c}
 1 & 0 & 0 \\
 0 & 1/ c_1 & 0 \\
 0 & 0 & 1 / c_2   
\end{array}
\right)
$
 & (V)$^*$ \\ 
\hline 
(1.3.a)
 &
$
\left(
\begin{array}{c c c}
 1 & 0 & 0 \\
 0 & 1 & 0 \\
 0 & 0 & 1  
\end{array}
\right)
$
 & (II)$^*$ \\ 
\hline 
(1.3.b)
 & 
$
\left(
\begin{array}{c c c}
 1 & 0 & 0 \\
 0 & c_1 & 0 \\
 0 & 0 & 1  
\end{array}
\right)
$
 & (VI)$^*$ \\ 
\hline 
(1.3.c)
 & 
$
\left(
\begin{array}{c c c}
 c_2 & 0 & 0 \\
 0 & 1 & 0 \\
 0 & 0 & 1  
\end{array}
\right)
$
 & (IV)$^*$ \\ 
\hline 
(1.3.d)
 & 
$
\left(
\begin{array}{c c c}
 c_2 & 0 & 0 \\
 0 & c_1 & 0 \\
 0 & 0 & 1  
\end{array}
\right)
$
 & (VII)$^*$ \\ 
\hline 
(2.a)
 & 
$
\left(
\begin{array}{c c c}
 1 & 0 & 0 \\
 0 & 1 & 0 \\
 0 & 0 & 1  
\end{array}
\right)
$
 & (VIII)$^*$ \\ 
\hline 
(2.b)
 & 
$
\left(
\begin{array}{c c c}
 1 & 0 & 0 \\
 0 & 1 & 0 \\
 0 & 0 & 1 / c_1  
\end{array}
\right)
$
 & (IX)$^*$ \\ 
\hline 
(2.c)
 & 
$
\left(
\begin{array}{c c c}
 0 & 1 & 0 \\
 1 & 0 & 0 \\
 0 & 0 & 1 / c_2  
\end{array}
\right)
$
 & (IX)$^*$, $e_1 : = e_2$ \\  
\hline 
(2.d)
 & 
$
\left(
\begin{array}{c c c}
 c_2 & c_2 & 0 \\
 0 & c_1 & 0 \\
 0 & 0 & 1  
\end{array}
\right)
$
 & (IX)$^*$ \\ 
\hline 
(3.a)
 & 
$
\left(
\begin{array}{c c c}
 1 & 0 & 0 \\
 0 & 1 & 0 \\
 0 & 0 & 1
\end{array}
\right)
$
 & (X)$^*$ \\ 
\hline 
(3.b)
 & 
$
\left(
\begin{array}{c c c}
 1 & 0 & 0 \\
 0 & 1 / c_1 & 0 \\
 0 & 0 & 1 
\end{array}
\right)
$
  & (XI)$^*$ \\ 
\hline 
(3.c)
 & 
$
\left(
\begin{array}{c c c}
 1 & 0 & 0 \\
 0 & 0 & 1 \\
 0 & 1 / c_2 & 0  
\end{array}
\right)
$
  & (XI)$^*$, $e_1 := e_2$ \\ 
\hline 
(3.d)
 & 
$
\left(
\begin{array}{c c c}
 1 & 0 & 0 \\
 0 & 1 / c_1 & - 1/ c_1 \\
 0 & 0 & 1 / c_2 
\end{array}
\right)
$
 & (XI)$^*$ \\ 
\hline 
(4) 
 & 
$
\left(
\begin{array}{c c c}
 1 & 0 & 0 \\
 0 & 1 & 0 \\
 0 & 0 & 1
\end{array}
\right)
$
 & (XII)$^*$ \\ 
\hline 
\end{longtable}

\subsubsection{Proof of assertion (2) of Theorem 2.3}

Let $V = k^{\oplus n}$ be the column vector space of dimension $n$, 
and let $V' := k^{\oplus n}$ be the row vector space of dimension $n$. 
For a representation $\varphi : \G_a \rtimes \G_m \to \SL(n, k)$ 
and for an integer $\ell$, we denote by $V^\varphi_\ell$ the subspace 
\begin{fleqn}[36pt] 
\begin{align*} 
 V^\varphi_\ell := 
\left\{
 v \in V 
\; 
\left|
\;
 \varphi
\left(
\begin{array}{c c} 
 a & b \\
 0 & d
\end{array}
\right) 
\, v  = a^\ell \, v \text{ \; for all \, 
$\left(
\begin{array}{c c} 
 a & b \\
 0 & d
\end{array}
\right)  \in \G_a \rtimes \G_m$} 
\right. 
\right\} 
\end{align*}
\end{fleqn} 
and denote by $(V')^\varphi_\ell$ the subspace defined by 
\begin{fleqn}[36pt] 
\begin{align*} 
 (V')^\varphi_\ell := 
\left\{
 v' \in V' 
\; 
\left|
\; 
v' \, 
 \varphi
\left(
\begin{array}{c c} 
 a & b \\
 0 & d
\end{array}
\right) 
 = a^\ell \, v' \text{ \; for all \, 
$\left(
\begin{array}{c c} 
 a & b \\
 0 & d
\end{array}
\right)  \in \G_a \rtimes \G_m$} 
\right. 
\right\} . 
\end{align*}
\end{fleqn}

\begin{lem}
Let $\varphi_i: \G_a \rtimes \G_m \to \SL(n, k)$ $(i = 1, 2)$ be two representations such that 
$\varphi_1$ is equivalent to $\varphi_2$. 
Let $\ell$ be an integer. 
Then $V^{\varphi_1}_\ell$ is isomorphic to $V^{\varphi_2}_\ell$ as vector space, 
and $(V')^{\varphi_1}_\ell$ is isomorphic to $(V')^{\varphi_2}_\ell$ as vector space. 
\end{lem}

\Proof There exists a regular representation $P$ of $\GL(n, k)$ so that 
\[
 P^{-1} 
\varphi_1
\left(
\begin{array}{c c}
 a & b \\
 0 & d
\end{array}
\right)
 P 
 =
 \varphi_2
\left(
\begin{array}{c c}
 a & b \\
 0 & d
\end{array}
\right)
\qquad 
\text{
 for all \quad $\left(
\begin{array}{c c}
 a & b \\
 0 & d
\end{array}
\right) \in \G_a \rtimes \G_m$. }
\]
For any $v \in V$, we have 
\begin{align*}
 v \in V^{\varphi_1}_\ell  
 &  \quad \Longleftrightarrow \quad \varphi_1 (A) \, v  = a^{\ell} \, v \text{ \; for all \, $A \in \G_a \rtimes \G_m$}  \\
 &  \quad \Longleftrightarrow \quad P \varphi_2(A) P^{- 1} \, v  = a^{\ell} \, v \text{ \; for all \, $A \in \G_a \rtimes \G_m$}  \\ 
 &  \quad \Longleftrightarrow \quad  \varphi_2 (A) \, P^{- 1} \, v  = a^{\ell} \, P^{- 1} \,  v \text{ \; for all \, $A \in \G_a \rtimes \G_m$}  \\ 
 &  \quad \Longleftrightarrow \quad   P^{- 1} \, v \in V^{\varphi_2}_\ell . 
\end{align*}
Thus $V^{\varphi_1}_\ell$ is isomorphic to $V^{\varphi_2}_\ell$. 
Similarly, $(V')^{\varphi_1}_\ell$ is isomorphic to $(V')^{\varphi_2}_\ell$. 
\QED

Let $\varphi : \G_a \rtimes \G_m \to \SL(n, k)$ be a representation such that 
\[
 \varphi
\left(
\begin{array}{c c}
 z & 0 \\
 0 & z^{-1}
\end{array}
\right) 
 = 
\left(
\begin{array}{c c c}
 z^{\ell_1} & 0 & 0 \\
 0 & z^{\ell_2} & 0 \\
 0 & 0 & z^{\ell_3} 
\end{array}
\right)
\qquad 
(\, \ell_1 \geq \ell_2 \geq \ell_3 \, ) . 
\]

If $\ell_1 > \ell_2 > \ell_3$, we can define elements $d_{(1)^*}(\varphi)$ and $d'_{(1)^*}(\varphi)$ of $\Z^3_{\geq 0}$ as 
\begin{fleqn}[36pt]
\begin{align*}
 d_{(1)^*}(\varphi) &:= (\, \dim_k \, (V)^\varphi_{\ell_1}, \; \dim_k \, (V)^\varphi_{\ell_2}, \; \dim_k \, (V)^\varphi_{\ell_3} \,) ,  \\
 d'_{(1)^*}(\varphi) &:= (\, \dim_k \, (V')^\varphi_{\ell_1}, \; \dim_k \, (V')^\varphi_{\ell_2}, \; \dim_k \, (V')^\varphi_{\ell_3} \,) . 
\end{align*} 
\end{fleqn}

If $\ell_1 = \ell_2 > 0 > \ell_3$, 
we can define elements $d_{(2)^*}(\varphi)$ and $d'_{(2)^*}(\varphi)$ of $\Z^2_{\geq 0}$ as 
\begin{fleqn}[36pt] 
\begin{align*}
 d_{(2)^*}(\varphi) &:= (\, \dim_k \, (V)^\varphi_{\ell_1}, \;   \dim_k \, (V)^\varphi_{\ell_3} \,) ,  \\
 d'_{(2)^*}(\varphi) &:= (\, \dim_k \, (V')^\varphi_{\ell_1}, \; \dim_k \, (V')^\varphi_{\ell_3} \,) . 
\end{align*} 
\end{fleqn}

If $\ell_1 > 0 >  \ell_2 =  \ell_3$, 
we can define elements $d_{(3)^*}(\varphi)$ and $d'_{(3)^*}(\varphi)$ of $\Z^2_{\geq 0}$ as 
\begin{fleqn}[36pt]
\begin{align*}
 d_{(3)^*}(\varphi) &:= (\, \dim_k \, (V)^\varphi_{\ell_1}, \;   \dim_k \, (V)^\varphi_{\ell_2} \,) ,  \\
 d'_{(3)^*}(\varphi) &:= (\, \dim_k \, (V')^\varphi_{\ell_1}, \; \dim_k \, (V')^\varphi_{\ell_2} \,) . 
\end{align*} 
\end{fleqn}

\begin{lem}
We have the following {\rm (1)$^*$}, {\rm (2)$^*$}, {\rm (3)$^*$}: 
\begin{enumerate}[label = {\rm (\arabic*)$^*$}]
\item Let $\ell_1$, $\ell_2$, $\ell_3$ be integers satisfying 
$\ell_1 + \ell_2 + \ell_3 = 0$, $\ell_1 > \ell_2 > \ell_3$ and $\ell_1 > 0 > \ell_3$. 
Let $\varphi^* : \G_a \rtimes \G_m \to \SL(3, k)$ be a representation with one of the 
forms {\rm (I)$^*$}, {\rm (II)$^*$}, {\rm (III)$^*$}, {\rm (IV)$^*$}, {\rm (V)$^*$}, {\rm (VI)$^*$}, {\rm (VII)$^*$}.  
\begin{enumerate}[label = {\rm (1.\arabic*)$^*$}, leftmargin = *]
\item 
If $\varphi^*$ has the form $\rm (I)^*$, 
then $d_{(1)^*}(\varphi^*) = (1, \, 0, \, 0) $, $d'_{(1)^*}(\varphi^*) = (0, \, 0, \, 1)$.

\item 
If $\varphi^*$ has the form $\rm (II)^*$, 
then $d_{(1)^*}(\varphi^*) = (1, \, 1, \, 1) $, $d'_{(1)^*}(\varphi^*) = (1, \, 1, \, 1)$.

\item 
If $\varphi^*$ has the form $\rm (III)^*$, 
then $d_{(1)^*}(\varphi^*) = (1, \, 0, \, 1) $, $d'_{(1)^*}(\varphi^*) = (0, \, 0, \, 1)$.

\item 
If $\varphi^*$ has the form $\rm (IV)^*$, 
then $d_{(1)^*}(\varphi^*) = (1, \, 1, \, 0) $, $d'_{(1)^*}(\varphi^*) = (1, \, 1, \, 1)$.

\item 
If $\varphi^*$ has the form $\rm (V)^*$, 
then $d_{(1)^*}(\varphi^*) = (1, \, 0, \, 0) $, $d'_{(1)^*}(\varphi^*) = (0, \, 1, \, 1)$.

\item 
If $\varphi^*$ has the form $\rm (VI)^*$, 
then $d_{(1)^*}(\varphi^*) = (1, \, 1, \, 0) $, $d'_{(1)^*}(\varphi^*) = (1, \, 0, \, 1)$.

\item 
If $\varphi^*$ has the form $\rm (VII)^*$, 
then $d_{(1)^*}(\varphi^*) = (1, \, 1, \, 0) $, $d'_{(1)^*}(\varphi^*) = (0, \, 0, \, 1)$. 
\end{enumerate}

\item  Let $\ell_1$, $\ell_2$, $\ell_3$ be integers satisfying 
$\ell_1 + \ell_2 + \ell_3 = 0$ and $\ell_1 = \ell_2 > 0 > \ell_3$. 
Let $\varphi^* : \G_a \rtimes \G_m \to \SL(3, k)$ be a representation with one of the 
forms {\rm (VIII)$^*$} and {\rm (IX)$^*$}.  
\begin{enumerate}[label = {\rm (2.\arabic*)$^*$}, leftmargin = *, start = 1]
\item 
If $\varphi^*$ has the form $\rm (VIII)^*$, 
then $d_{(2)^*}(\varphi^*) = (2, \, 1) $, $d'_{(2)^*}(\varphi^*) = (2, \, 1)$.

\item 
If $\varphi^*$ has the form $\rm (IX)^*$, 
then $d_{(2)^*}(\varphi^*) = (2, \, 0) $, $d'_{(2)^*}(\varphi^*) = (1, \, 1)$. 
\end{enumerate}

\item Let $\ell_1$, $\ell_2$, $\ell_3$ be integers satisfying $\ell_1 + \ell_2 + \ell_3 = 0$ 
and $\ell_1 > 0 > \ell_2 = \ell_3$. 
Let $\varphi^* : \G_a \rtimes \G_m \to \SL(3, k)$ be a representation with one of the 
forms {\rm (X)$^*$} and {\rm (XI)$^*$}.  
\begin{enumerate}[label = {\rm (3.\arabic*)$^*$}, leftmargin = *, start = 1]
\item 
If $\varphi^*$ has the form $\rm (X)^*$, 
then $d_{(3)^*}(\varphi^*) = (1, \, 2) $, $d'_{(3)^*}(\varphi^*) = (1, \, 2)$.

\item 
If $\varphi^*$ has the form $\rm (XI)^*$, 
then $d_{(3)^*}(\varphi^*) = (1, \, 1) $, $d'_{(3)^*}(\varphi^*) = (0, \, 2)$.

\end{enumerate}

\end{enumerate} 
\end{lem}

\Proof The proof is straightforward. 
\QED

Now, we give a proof of assertion (2) of Theorem 2.3. 

We first prove (2.1). 
If $\lambda \in \Lambda^{(i)^*}$ and $\mu \in \Lambda^{(j)^*}$, where $i, j \in  \{  1, \, 2, \, 3 , \, 4  \}$ with $i \ne j$, 
we know from Lemma 1.2 that $\phi$ and $\psi$ are not equivalent. 
If $\lambda, \mu \in \Lambda^{(i)^*}$ and $\lambda \ne \mu$, 
we know from Lemmas 2.5 and 2.6 that $\phi$ and $\psi$ are not equivalent.

We next prove (2.2). 
The implication (i) $\Longrightarrow$ (ii) is clear. 
The implication (ii) $\Longrightarrow$ (iii) follows from Lemma 1.2. 
The implications (iii) $\Longrightarrow$ (iv) and (iv) $\Longrightarrow$ (i) are clear.

\section{Fundamental representations of $\G_a \rtimes \G_m$ into $\SL(3, k)$ }

\subsection{Candidates for classifying fundamental representations of $\G_a \rtimes \G_m$ into $\SL(3, k)$}

A representation $\varphi : \G_a \rtimes \G_m \to \SL(n, k)$ of $\G_a \rtimes \G_m$ 
is said to be {\it fundamental} if there exists a representation $\psi : \SL(2, k) \to \SL(n, k)$ such that 
the following diagram commutes: 
\[
\xymatrix@R=36pt@C=36pt@M=8pt{
 \G_a \rtimes \G_m \ar[r]^\varphi \ar@{^(->}[d]_\iota &  \SL(3, k) \\
 \SL(2, k) \ar[ru]_\psi 
}
\]
where $\iota : \G_a \rtimes \G_m \to \SL(2, k)$ is the injective homomorphism defined by 
\[
\iota(t, z)
 := 
\left(
\begin{array}{c c}
 z & t \, z^{- 1} \\
 0 & z^{- 1}
\end{array}
\right)  . 
\]
So, identifying an element $(a \, b, \, a)$ of $\G_a \rtimes \G_m$ with the element 
$\left(
\begin{array}{c c}
 a & b \\
 0 & d 
\end{array}
\right) $ 
of $\SL(2, k)$, 
we have 
\[
 \psi
\left(
\begin{array}{c c}
 a & b \\
 0 & d 
\end{array}
\right) 
 = 
 \varphi
\left(
\begin{array}{c c}
 a & b \\
 0 & d 
\end{array}
\right) . 
\]

\begin{lem}
Let $\varphi : \G_a \rtimes \G_m \to \SL(n, k)$ be a fundamental representation. 
Let $P$ be a regular matrix of $\GL(n, k)$ and let $\varphi^* : \G_a \rtimes \G_m \to \SL(n, k)$ be the representation 
defined by 
\[
 \varphi^*
\left(
\begin{array}{c c}
 a & b \\
 c & d 
\end{array}
\right) 
 := 
P^{-1} 
\varphi
\left(
\begin{array}{c c}
 a & b \\ 
 c & d 
\end{array}
\right) 
P . 
\]
Then $\varphi^*$ is a fundamental representation of $\G_a \rtimes \G_m$ into $\SL(n, k)$. 
\end{lem}

\Proof Let $f : \SL(3, k) \to \SL(3, k)$ be a homomorphism defined by $f(A) := P^{-1} A P$. 
Clearly, $\varphi^* = f \circ \varphi$. 
Then $\varphi^* = f \circ \psi \circ \iota$, which implies $\varphi^*$ is a fundamental representation 
of $\G_a \rtimes \G_m$. 
\QED

\begin{lem}
Let $\varphi^* : \G_a \rtimes \G_m \to \SL(3, k)$ be a representation of $\G_a \rtimes \G_m$ 
with one of the forms {\rm (I)$^*$ -- (XII)$^*$}. 
Assume $\varphi^*$ is fundamental. 
Then $\varphi^*$ has one of the forms {\rm (I)$^*$ -- (VII)$^*$} and  {\rm (XII)$^*$}. 
Furthermore, we have the following:

\begin{enumerate}[label = {\rm (\arabic*)}, leftmargin = *]

\item If $\varphi^*$ has the form {\rm (II)$^*$}, then $\ell_1 > 0$ and $\ell_3 = - \ell_1$.

\item If $\varphi^*$ has the form {\rm (III)$^*$}, then $(\ell_1, \ell_2, \ell_3) = (2 \, p^{e_1}, 0 , - 2 \, p^{e_1})$.

\item If $\varphi^*$ has the form {\rm (IV)$^*$}, then $(\ell_1, \ell_2, \ell_3) = (p^{e_2}, 0 , - p^{e_2})$.

\item If $\varphi^*$ has the form {\rm (V)$^*$}, then 
$p = 2$ and $(\ell_1, \ell_2, \ell_3) = (2 \, p^{e_1}, 0 , - 2 \, p^{e_1})$.

\item If $\varphi^*$ has the form {\rm (VI)$^*$}, then $(\ell_1, \ell_2, \ell_3) = (2 \, p^{e_1}, 0 , - 2 \, p^{e_1})$.

\item 
If $\varphi^*$ has the form {\rm (VII)$^*$}, then $p = 2$ and $(\ell_1, \ell_2, \ell_3) = (2 \, p^{e_1}, 0 , - 2 \, p^{e_1})$.

\end{enumerate}

\end{lem}

\Proof 
We know from Lemma 1.6 that the integers $\ell_1$, $\ell_2$, $\ell_3$ ($\ell_1 \geq \ell_2 \geq \ell_3$) 
satisfies $\ell_3 = - \ell_1$ and $\ell_2 = 0$. 
So, if $\ell_1 \ne \ell_2$, then $\ell_2 \ne \ell_3$. Then $\varphi^*$ has one of the forms (I)$^*$ -- (VII)$^*$. 
If $\ell_1 = \ell_2$, then $\ell_1 = \ell_2 = \ell_3 = 0$. We know from Lemma 1.5 that $u_{\varphi^*}$ is trivial, 
and thereby $\varphi^*$ has the form (XII)$^*$.

The proofs of (1), (2), (3), (5) are straightforward. 
\begin{enumerate}[label = {\rm (\roman*)}, leftmargin = *]
\item[(4)] If $\varphi^*$ has the form (V)$^*$, then 
$\ell_1 = (2 \, p^{e_1} + 2 \, p^{e_2} ) / 3$, 
$\ell_2 = ( -4 \, p^{e_1} + 2 \, p^{e_2} ) / 3$, 
$\ell_3 = ( 2 \, p^{e_1} - 4 \, p^{e_2} ) / 3$, 
where $e_2 > e_1 \geq 0$. 
Since $\ell_1 = - \ell_3$ and $\ell_2 = 0$, we have $4 \, p^{e_1} = 2 \, p^{e_2}$ and $p^{e_2} = 2 \, p^{e_1}$, which implies 
$p = 2$ and $e_2 = e_1 + 1$. Thus $\ell_1 = 2 \, p^{e_1}$ and $\ell_3 = - 2 \, p^{e_1}$. 

\item[(6)]  
If $\varphi^*$ has the form (VII)$^*$, then 
$\ell_1 = ( - 2 \, p^{e_1} + 4 \, p^{e_2} ) /3$, 
$\ell_2 = ( 4 \, p^{e_1} - 2 \, p^{e_2} ) / 3$, 
$\ell_3 =  ( - 2 \, p^{e_1} - 2 \, p^{e_2} ) / 3$, 
where $e_2 > e_1 \geq 0$. 
Since $\ell_1 = - \ell_3$ and $\ell_2 = 0$. Thus $- 4 \, p^{e_1} = - 2 \, p^{e_2}$ and $p^{e_2} = 2 \, p^{e_1}$, which implies 
$p = 2$ and $e_2 = e_1 + 1$. Thus $\ell_1 = 2 \, p^{e_1}$ and $\ell_3 = - 2 \, p^{e_1}$. 
 
\end{enumerate} 
\QED

\begin{lem}
The following assertions {\rm (1)} and {\rm (2)} hold true: 
\begin{enumerate}[label = {\rm (\arabic*)}]
\item 
$
\left(
\begin{array}{c c}
 1 & 1 \\
 0 & 1
\end{array} 
\right) 
\left(
\begin{array}{c c}
 1 & 0 \\
 \gamma & 1
\end{array} 
\right)
 = 
\left(
\begin{array}{c c}
 1 & 0 \\
 \frac{\gamma}{1 + \gamma} & 1
\end{array} 
\right) 
\left(
\begin{array}{c c}
 1 + \gamma & 0 \\
 0 & \frac{1}{1 + \gamma}
\end{array} 
\right)
\left(
\begin{array}{c c}
 1 & \frac{1}{1 + \gamma} \\
 0 & 1
\end{array} 
\right)
$
for all $\gamma \in k \backslash \{ - 1 \}$. 

\item Let $\varphi : \G_a \rtimes \G_m \to \SL(n, k)$ be a fundamental representation. 
Then we have
\begin{fleqn}[36pt] 
\begin{align*} 
 u(1) \, u^-(\gamma) = u^-\left( \frac{\gamma}{1 + \gamma} \right) \, h(1 + \gamma) \, u\left(\frac{1}{1 + \gamma} \right)
\qquad 
\text{ for all \quad  $\gamma \in k \backslash \{ - 1 \}$.  }
\end{align*}
\end{fleqn}
\end{enumerate} 
\end{lem}

\Proof 
The proofs of assertions (1) and (2) are straightforward. 
\QED

\subsubsection{(I)$^*$}

\begin{lem}
Let $\varphi^* : \G_a \rtimes \G_m \to \SL(3, k)$ be a representation with the form 
\begin{fleqn}[36pt] 
\begin{align*} 
\varphi^*
\left(
\begin{array}{c c}
 a & b \\
 0 & d 
\end{array}
\right) 
 & =
\left(
\begin{array}{c c c}
 a^{2 \, p^{e_1}}  &  a^{p^{e_1}} \, b^{p^{e_1}}  & \frac{1}{2} \,  b^{2 \, p^{e_1}} \\
 0  &  1  &   b^{p^{e_1}} \,  d^{p^{e_1}} \\
 0 & 0 & d^{2 \, p^{e_1}} 
\end{array}  
\right) 
\qquad 
(\, e_1 \geq 0 \, ) . 
\end{align*}
\end{fleqn}
Then the following assertions {\rm (1)}, {\rm (2)}, {\rm (3)} hold true: 
\begin{enumerate}[label = {\rm (\arabic*)}]
\item $\varphi^*$ is fundamental. 

\item Let $\psi^* : \SL(2, k) \to \SL(3, k)$ be a representation of $\SL(2, k)$ so that $\varphi^* = \psi^* \circ \iota$. 
Then we have 
\begin{fleqn}[36pt] 
\begin{align*}
 u_{\psi^*}^-(s) 
 = 
\left(
\begin{array}{c c c}
 1 & 0 & 0 \\
 2 \, s^{p^{e_1}}  & 1 & 0 \\
 2 \, s^{2 \, p^{e_1}} & 2 \, s^{p^{e_1}} & 1 
\end{array}
\right) .
\end{align*}
\end{fleqn} 

\item There exists a unique representation $\psi^* : \SL(2, k) \to \SL(3, k)$ such that 
$\varphi^* = \psi^* \circ \iota$. 
\end{enumerate} 
\end{lem}

\Proof 
(1) 
Consider the representation $\psi^* : \SL(2, k) \to \SL(3, k)$ defined by 
\begin{fleqn}[36pt] 
\begin{align*} 
 \psi^*
\left(
\begin{array}{c c}
 a & b \\
 c & d
\end{array}
\right) 
 := 
\left(
\begin{array}{c c c}
 a^{2 \, p^{e_1}} & a^{p^{e_1}} \, b^{p^{e_1}}  & \frac{1}{2} \, b^{2 \, {p^{e_1}}}  \\
 2 \, a^{p^{e_1}} \, c^{p^{e_1}} & a^{p^{e_1}} \, d^{p^{e_1}} + b^{p^{e_1}} \, c^{p^{e_1}} &  b^{p^{e_1}} \, d^{p^{e_1}} \\
 2 \, c^{2 \, p^{e_1}} & 2 \, c^{p^{e_1}} \, d^{p^{e_1}}  & d^{2 \, p^{e_1}}  
\end{array}
\right) . 
\end{align*}
\end{fleqn}

(2) 
We have 
\begin{fleqn}[36pt] 
\begin{align*}
 u_{\psi^*}(t) 
  = 
\left(
\begin{array}{c c c}
 1 & t^{p^{e_1}} & \frac{1}{2} \, t^{2 \, p^{e_1}}  \\
 0 & 1 & t^{p^{e_1}} \\
 0 & 0 & 1
\end{array}
\right) . 
\end{align*}
\end{fleqn} 
We can express $u^-_{\psi^*} : \G_a \to \SL(3, k)$ as  
\begin{fleqn}[36pt] 
\begin{align*}
 u^-_{\psi^*}(s)
  = 
\left(
\begin{array}{c c c}
 1 & 0 & 0 \\
 v_{2, \, 1}(s) & 1 & 0 \\
 v_{3, \, 1}(s) & v_{3, \, 2}(s) & 1 
\end{array}
\right) 
\qquad 
(\  v_{2, \, 1}(S), v_{3, \, 1}(S) ,  v_{3, \, 2}(S)  \in k[S]  \,) . 
\end{align*}
\end{fleqn} 
Using Lemma 3.3 (2), we have, for all $\gamma \in k \backslash \{ - 1 \}$, 
\begin{fleqn}[36pt] 
\begin{align*}
& \hspace{-24pt} 
\left(
\begin{array}{c c c}
 1 + v_{2, \, 1}(\gamma) + \frac{1}{2} \, v_{3, \, 1}(\gamma)  & 1 + \frac{1}{2} \, v_{3, \, 2}(\gamma)  & \frac{1}{2}  \\  [6pt] 
 v_{2, \, 1}(\gamma)  +  v_{3, \, 1}(\gamma) & 1 + v_{3, \, 2}(\gamma) &  1 \\  [6pt] 
 v_{3, \, 1}(\gamma) & v_{3, \, 2}(\gamma) & 1
\end{array}
\right) \\
 & \qquad 
\hspace{-24pt} 
= 
\left(
\begin{array}{c c }
 (1 + \gamma)^{2 \, p^{e_1}}
 & (1 + \gamma)^{p^{e_1}}  \\  [6pt] 
 v_{2, \, 1} \left( \frac{\gamma}{1 + \gamma} \right) \, (1 + \gamma)^{2 \, p^{e_1}}
 &  v_{2, \, 1} \left( \frac{\gamma}{1 + \gamma} \right) \, (1 + \gamma)^{p^{e_1}} + 1  \\ [6pt] 
 v_{3, \, 1} \left( \frac{\gamma}{1 + \gamma} \right) \, (1 + \gamma)^{2 \, p^{e_1}}
 &  v_{3, \, 1} \left( \frac{\gamma}{1 + \gamma} \right) \, (1 + \gamma)^{p^{e_1}}  + v_{3, \, 2}  \left( \frac{\gamma}{1 + \gamma} \right) 
\end{array}
\right.  \\ 
 & \hspace{180pt}
\left. 
\begin{array}{c}
 \frac{1}{2}  \\  [6pt] 
  \frac{1}{2} \,  v_{2, \, 1} \left( \frac{\gamma}{1 + \gamma} \right) +  \frac{1}{(1 + \gamma)^{ p^{e_1}} }  \\  [6pt] 
 \frac{1}{2} \,  v_{3, \, 1} \left( \frac{\gamma}{1 + \gamma} \right) 
    + v_{3, \, 2}  \left( \frac{\gamma}{1 + \gamma} \right) \,  \frac{1}{(1 + \gamma)^{ p^{e_1}} }  +  \frac{1}{(1 + \gamma)^{ 2 \, p^{e_1}} }  
\end{array}
\right) .  
\end{align*}
\end{fleqn} 
Comparing the $(1, 2)$-th entries of both sides of the above equality $(\ast)$, we have 
$1 +  v_{3, \, 2}(\gamma) / 2 = 1 + \gamma^{p^{e_1}} $ for all $\gamma \in k \backslash \{ - 1\}$, 
which implies $v_{3, \, 2}(S) = 2 \, S^{p^{e_1}}$.  
Comparing the $(2, 3)$-th entries of both sides of the equality $(\ast)$, we have 
$1 =  \frac{1}{2} \, v_{2, \, 1}\left( \frac{\gamma}{1 + \gamma} \right) 
 + \frac{1}{(1 + \gamma)^{p^{e_1}}} = 1 $ 
for all $\gamma \in k \backslash \{ - 1\}$, 
which implies $v_{2, \, 1}(S) = 2 \, S^{ p^{e_1} }$. 
Comparing the $(1, 1)$-th entries of both sides of the equality $(\ast)$, we have 
$ 1 + v_{2, \, 1}(\gamma) + v_{3, \, 1}(\gamma) / 2 =  (1 + \gamma)^{2 \, p^{e_1}}$, 
which implies $v_{3, \, 1}(S) = 2 \, S^{2 \, p^{e_1}}$.

(3) 
We already know the existence of $\psi^*$ (see the above assertion (1)). 
Let $\psi_1^*$, $\psi_2^*$ be two representations $\SL(2, k) \to \SL(3, k)$ of $\SL(2, k)$ such that 
$\psi_1^* \circ \iota = \psi_2^* \circ \iota = \varphi^*$. 
We know from the above assertion (2) that $u^-_{\psi_1^*} = u^-_{\psi_2^*}$, and then know from Lemma 1.10 that $\psi_1^* = \psi_2^*$. 
\QED

\subsubsection{(II)$^*$: $\ell_1 > 0$ and $\ell_3 = - \ell_1$}

\begin{lem}
Let $\varphi^* : \G_a \rtimes \G_m \to \SL(3, k)$ be a representation with the form 
\begin{fleqn}[36pt] 
\begin{align*} 
\varphi^*
\left(
\begin{array}{c c}
 a & b \\
 0 & d 
\end{array}
\right) 
 & =
\left(
\begin{array}{c c c}
 a^{\ell_1} & 0 & 0 \\
 0 & 1 & 0 \\
 0 & 0 & d^{ - \ell_3} 
\end{array} 
\right)
\qquad (\, \ell_1 > 0 > \ell_3 \, ) . 
\end{align*}
\end{fleqn}
Then $\varphi^*$ is not fundamental. 
\end{lem}

\Proof Suppose, to the contrary, that $\varphi^*$ is fundamental. 
Since $u_{\varphi^*}$ is trivial, we know from Lemma 1.9 that $\ell_1 =  0$. This contradicts $\ell_1 > 0$. 
\QED

\subsubsection{(III)$^*$: $(\ell_1, \ell_2, \ell_3) = (2 \, p^{e_1}, 0, - 2 \, p^{e_1})$}

\begin{lem}
Let $\varphi^* : \G_a \rtimes \G_m \to \SL(3, k)$ be a representation with the form 
\begin{fleqn}[36pt] 
\begin{align*}
\varphi^*
\left(
\begin{array}{c c}
 a & b \\
 0 & d 
\end{array}
\right) 
& = 
\left(
\begin{array}{c c c}
 a^{2 \, p^{e_1}} & a^{p^{e_1}} \, b^{p^{e_1}}  & 0 \\
 0 & 1 & 0 \\
 0 & 0 & d^{ 2 \, p^{e_1}} 
\end{array} 
\right) 
 \qquad  ( \, e_1 \geq 0 \, ) . 
\end{align*}
\end{fleqn} 
Then $\varphi^*$ is not fundamental. 
\end{lem}

\Proof 
We have 
\begin{fleqn}[36pt] 
\begin{align*}
u_{\psi^*}(t) 
& =
\left(
\begin{array}{c c c}
 1 & t^{p^{e_1}} & 0  \\
 0 & 1 & 0 \\
 0 & 0 & 1 
\end{array} 
\right) . 
\end{align*}
\end{fleqn} 
Suppose, to the contrary, that $\varphi^*$ is fundamental. 
We can express $u^-_{\psi^*} : \G_a \to \SL(3, k)$ as  
\begin{fleqn}[36pt] 
\begin{align*}
 u^-_{\psi^*}(s)
  = 
\left(
\begin{array}{c c c}
 1 & 0 & 0 \\
 v_{2, \, 1}(s) & 1 & 0 \\
 v_{3, \, 1}(s) & v_{3, \, 2}(s) & 1 
\end{array}
\right) 
\qquad 
(\  v_{2, \, 1}(S), v_{3, \, 1}(S) ,  v_{3, \, 2}(S)  \in k[S]  \,) . 
\end{align*}
\end{fleqn} 
Using Lemma 3.3 (2), we have, for all $\gamma \in k \backslash \{ - 1 \}$, 
\begin{fleqn}[36pt] 
\begin{align*}
& 
\left(
\begin{array}{c c c}
 1 + v_{2, \, 1}(\gamma) & 1 & 0 \\
 v_{2, \, 1}(\gamma) & 1 &  0 \\
 v_{3, \, 1}(\gamma) & v_{3, \, 2}(\gamma) & 1
\end{array}
\right) \\
 & \qquad = 
\left(
\begin{array}{c c c}
 (1 + \gamma)^{2 \, p^{e_1}} & (1 + \gamma)^{p^{e_1}} & 0 \\  [6pt] 
 v_{2, \, 1} \left( \frac{\gamma}{1 + \gamma} \right) \, (1 + \gamma)^{2 \, p^{e_1}}
 &  v_{2, \, 1} \left( \frac{\gamma}{1 + \gamma} \right) \, (1 + \gamma)^{p^{e_1}} + 1 
 & 0 \\ [6pt] 
 v_{3, \, 1} \left( \frac{\gamma}{1 + \gamma} \right) \, (1 + \gamma)^{2 \, p^{e_1}}
 &  v_{3, \, 1} \left( \frac{\gamma}{1 + \gamma} \right) \, (1 + \gamma)^{p^{e_1}} 
    + v_{3, \, 2}  \left( \frac{\gamma}{1 + \gamma} \right) 
 & \frac{1}{(1 + \gamma)^{2 \, p^{e_1}} }
\end{array}
\right) . 
\end{align*}
\end{fleqn} 
Comparing the $(3, 3)$-th entries of both sides of the above equality, we have $1 = 1 / (1 + \gamma)^{2 \, p^{e_1}}$ 
for all $\gamma \in k \backslash \{ - 1\}$. 
Considering an element $\gamma$ of $k$ such that $\gamma \notin k \backslash \{-1 , \, 0 \}$ and $\gamma^{p^{e_1}} \ne -2$, 
we have a contradiction. 
\QED

\subsubsection{(IV)$^*$: $(\ell_1, \ell_2, \ell_3) = (p^{e_2}, 0, - p^{e_2})$}

\begin{lem}
Let $\varphi^* : \G_a \rtimes \G_m \to \SL(3, k)$ be a representation with the form 
\begin{fleqn}[36pt] 
\begin{align*}
\varphi^*
\left(
\begin{array}{c c}
 a & b \\
 0 & d 
\end{array}
\right) 
& =
\left(
\begin{array}{c c c}
 a^{p^{e_2}} & 0 &   b^{p^{e_2}}  \\
 0 & 1 & 0 \\
 0 & 0 &  d^{p^{e_2}} 
\end{array} 
\right) 
\qquad  ( \, e_2 \geq 0 \, ) . 
\end{align*}
\end{fleqn} 
Then the following assertions {\rm (1)}, {\rm (2)}, {\rm (3)} hold true: 
\begin{enumerate}[label = {\rm (\arabic*)}]
\item $\varphi^*$ is fundamental. 

\item Let $\psi^* : \SL(2, k) \to \SL(3, k)$ be a representation of $\SL(2, k)$ so that $\varphi^* = \psi^* \circ \iota$. 
Then we have 
\begin{fleqn}[36pt] 
\begin{align*}
 u^-_{\psi^*}(s) 
 = 
\left(
\begin{array}{c c c}
 1 & 0 & 0 \\
 0 & 1 & 0 \\
 s^{p^{e_2}} & 0 & 1 
\end{array}
\right) .
\end{align*}
\end{fleqn} 

\item There exists a unique representation $\psi^* : \SL(2, k) \to \SL(3, k)$ such that 
$\varphi^* = \psi^* \circ \iota$. 
\end{enumerate} 
\end{lem}

\Proof 
(1) 
Consider the representation $\psi^* : \SL(2, k) \to \SL(3, k)$ defined by 
\begin{fleqn}[36pt] 
\begin{align*} 
 \psi^*
\left(
\begin{array}{c c}
 a & b \\
 c & d
\end{array}
\right) 
 := 
\left(
\begin{array}{c c c}
 a^{p^{e_2}} & 0 &   b^{p^{e_2}}  \\
 0 & 1 & 0 \\
 c^{p^{e_2}} & 0 &  d^{p^{e_2}} 
\end{array} 
\right) . 
\end{align*}
\end{fleqn}

(2) We have 
\begin{fleqn}[36pt] 
\begin{align*}
 u_{\psi^*}(t) 
  = 
\left(
\begin{array}{c c c}
 1 & 0 & t^{p^{e_2}}  \\
 0 & 1 & 0 \\
 0 & 0 & 1
\end{array}
\right) . 
\end{align*}
\end{fleqn} 
We can express $u^-_{\psi^*} : \G_a \to \SL(3, k)$ as  
\begin{fleqn}[36pt] 
\begin{align*}
 u^-_{\psi^*}(s)
  = 
\left(
\begin{array}{c c c}
 1 & 0 & 0 \\
 v_{2, \, 1}(s) & 1 & 0 \\
 v_{3, \, 1}(s) & v_{3, \, 2}(s) & 1 
\end{array}
\right) 
\qquad 
(\  v_{2, \, 1}(S), v_{3, \, 1}(S) ,  v_{3, \, 2}(S)  \in k[S]  \,) . 
\end{align*}
\end{fleqn} 
Using Lemma 3.3 (2), we have, for all $\gamma \in k \backslash \{ - 1 \}$, 
\begin{fleqn}[36pt] 
\begin{align*}
& 
\left(
\begin{array}{c c c}
 1 + v_{3, \, 1}(\gamma) & v_{3, \, 2}(\gamma) & 1 \\
 v_{2, \, 1}(\gamma) & 1 &  0 \\
 v_{3, \, 1}(\gamma) & v_{3, \, 2}(\gamma) & 1
\end{array}
\right) \\
 & \qquad = 
\left(
\begin{array}{c c c}
 (1 + \gamma)^{p^{e_2}} & 0 & 1 \\  [6pt] 
 v_{2, \, 1} \left( \frac{\gamma}{1 + \gamma} \right) \, (1 + \gamma)^{ p^{e_2}}
 & 1 
 &  v_{2, \, 1} \left( \frac{\gamma}{1 + \gamma} \right)  \\ [6pt] 
 v_{3, \, 1} \left( \frac{\gamma}{1 + \gamma} \right) \, (1 + \gamma)^{p^{e_2}}
 &  v_{3, \, 2}  \left( \frac{\gamma}{1 + \gamma} \right) 
 & v_{3, \, 1}  \left( \frac{\gamma}{1 + \gamma} \right) +  \frac{1}{(1 + \gamma)^{p^{e_2}} }
\end{array}
\right) . 
\tag{$\ast$}
\end{align*}
\end{fleqn}  
Comparing the $(1, 2)$-th entries of both sides of the above equality $(\ast)$, we have $v_{3, \, 2}(\gamma) = 0$ for all 
$\gamma \in k \backslash \{ - 1 \}$, which implies $v_{3,\, 2}(S) = 0$. 
Comparing the $(2, 3)$-th entries of both sides of the above equality $(\ast)$, 
we have $0 = v_{2, \, 1}(\frac{\gamma}{1 + \gamma}) $ for all 
$\gamma \in k \backslash \{ - 1 \}$, which implies $ v_{2, \, 1}(S) = 0$. 
Comparing the $(1, 1)$-th entries of both sides of the equality $(\ast)$, we have 
$1 +  v_{3, \, 1} (\gamma) = (1 + \gamma)^{p^{e_2}}$ for all $\gamma \in k \backslash \{ - 1 \}$, 
which implies $ v_{3, \, 1} (S) = S^{ p^{e_2} }$.

(3) 
The proof is similar to the proof of assertion (3) of Lemma 3.4 (use the above assertions (1) and (2)).  
\QED

\subsubsection{(V)$^*$: $p = 2$ and $(\ell_1, \ell_2, \ell_3) = (2 \, p^{e_1}, 0, - 2 \, p^{e_1})$}

\begin{lem}
If $p = 2$, we let $\varphi^* : \G_a \rtimes \G_m \to \SL(3, k)$ be a representation with the form 
\begin{fleqn}[36pt] 
\begin{align*}
\varphi^* 
\left(
\begin{array}{c c}
 a & b \\
 0 & d 
\end{array}
\right) 
  = 
\left(
\begin{array}{c c c}
 a^{2 \, p^{e_1}}  &  a^{p^{e_1}} \, b^{p^{e_1}}   &  b^{2 \, p^{e_1}}\\
 0  &  1  &  0  \\
 0 & 0 & d^{ 2 \, p^{e_1}} 
\end{array}  
\right) 
 \qquad  ( \, e_1 \geq 0 \, ) . 
\end{align*}
\end{fleqn} 
Then the following assertions {\rm (1)}, {\rm (2)}, {\rm (3)} hold true: 
\begin{enumerate}[label = {\rm (\arabic*)}]
\item $\varphi^*$ is fundamental. 

\item Let $\psi^* : \SL(2, k) \to \SL(3, k)$ be a representation of $\SL(2, k)$ so that $\varphi^* = \psi^* \circ \iota$. 
Then we have 
\begin{fleqn}[36pt] 
\begin{align*}
 u_{\psi^*}^-(s) 
 = 
\left(
\begin{array}{c c c}
 1 & 0 & 0 \\
 0 & 1 & 0 \\
 s^{2 \, p^{e_1}} & s^{p^{e_1}}  & 1 
\end{array}
\right) .
\end{align*}
\end{fleqn}

\item There exists a unique representation $\psi^* : \SL(2, k) \to \SL(3, k)$ such that 
$\varphi^* = \psi^* \circ \iota$. 
\end{enumerate} 
\end{lem}

\Proof  
(1) 
Consider the representation $\psi^* : \SL(2, k) \to \SL(3, k)$ defined by 
\begin{fleqn}[36pt] 
\begin{align*} 
\psi^*
\left(
\begin{array}{c c}
 a & b \\
 c & d 
\end{array}
\right) 
  = 
\left(
\begin{array}{c c c}
 a^{2 \, p^{e_1}}  & a^{p^{e_1}} \, b^{p^{e_1}}   &   b^{2 \, p^{e_1}}\\
 0  &  1  &  0  \\
 c^{ 2 \, p^{e_1}}  & c^{p^{e_1}} \, d^{p^{e_1}}  & d^{ 2 \, p^{e_1}} 
\end{array}  
\right) . 
\end{align*}
\end{fleqn}

(2) We have 
\begin{fleqn}[36pt] 
\begin{align*}
 u_{\psi^*}(t) 
  = 
\left(
\begin{array}{c c c}
 1 & t^{p^{e_1}} & t^{2 \, p^{e_1}}  \\
 0 & 1 & 0 \\
 0 & 0 & 1
\end{array}
\right) . 
\end{align*}
\end{fleqn} 
We can express $u^-_{\psi^*} : \G_a \to \SL(3, k)$ as  
\begin{fleqn}[36pt] 
\begin{align*}
 u^-_{\psi^*}(s)
  = 
\left(
\begin{array}{c c c}
 1 & 0 & 0 \\
 v_{2, \, 1}(s) & 1 & 0 \\
 v_{3, \, 1}(s) & v_{3, \, 2}(s) & 1 
\end{array}
\right) 
\qquad 
(\  v_{2, \, 1}(S), v_{3, \, 1}(S) ,  v_{3, \, 2}(S)  \in k[S]  \,) . 
\end{align*}
\end{fleqn} 
Using Lemma 3.3 (2), we have, for all $\gamma \in k \backslash \{ - 1 \}$, 
\begin{align*}
& 
\left(
\begin{array}{c c c}
 1 + v_{2, \, 1}(\gamma) + v_{3,\, 1}(\gamma) & 1 + v_{3, \, 2}(\gamma) & 1 \\
 v_{2, \, 1}(\gamma) & 1 & 0 \\
 v_{3, \, 1}(\gamma) & v_{3, \, 2}(\gamma) & 1 
\end{array}
\right) \\ 
 & \qquad = 
\left(
\begin{array}{c c c}
 (1 + \gamma)^{2 \, p^{e_1}} & (1 + \gamma)^{p^{e_1}} & 1 \\  [6pt] 
 v_{2, \, 1} \left( \frac{\gamma}{1 + \gamma} \right) \, (1 + \gamma)^{2 \, p^{e_1}}
 &  v_{2, \, 1} \left( \frac{\gamma}{1 + \gamma} \right) \, (1 + \gamma)^{p^{e_1}} + 1 
 & v_{2, \, 1} \left( \frac{\gamma}{1 + \gamma} \right) \\ [6pt] 
 v_{3, \, 1} \left( \frac{\gamma}{1 + \gamma} \right) \, (1 + \gamma)^{2 \, p^{e_1}}
 &  v_{3, \, 1} \left( \frac{\gamma}{1 + \gamma} \right) \, (1 + \gamma)^{p^{e_1}} 
    + v_{3, \, 2}  \left( \frac{\gamma}{1 + \gamma} \right) 
 & v_{3, \, 1} \left( \frac{\gamma}{1 + \gamma} \right) + \frac{1}{(1 + \gamma)^{2 \, p^{e_1}} }
\end{array}
\right) . 
\tag{$\ast$}
\end{align*}
Comparing the $(1, 2)$-th entries of both sides of the above equality $(\ast)$, 
we have $1 + v_{3, \, 2}(\gamma) = (1 + \gamma)^{p^{e_1}}$ for all 
$\gamma \in k \backslash \{ - 1 \}$, which implies $v_{3, \, 2}(S) = S^{p^{e_1}}$. 
Comparing the $(2, 3)$-th entries of both sides of the above equality $(\ast)$, 
we have $0 = v_{2,\, 1}(\frac{\gamma}{1 + \gamma})$ for all 
$\gamma \in k \backslash \{ -1 \}$, which implies $v_{2,\, 1}(S) = 0$. 
Comparing the $(1, 1)$-th entries of both sides of the equality $(\ast)$, we have 
$1 + v_{2, \, 1}(\gamma) + v_{3, \, 1}(\gamma) = (1 + \gamma)^{2 \,p^{e_1}}$ for all $\gamma \in k \backslash \{ - 1 \}$, 
which implies $v_{3, \, 1}(S) = S^{2 \, p^{e_1} }$ (since $p = 2$).

(3) The proof is straightforward. 
\QED

\subsubsection{(VI)$^*$: $(\ell_1, \ell_2, \ell_3) = (2 \, p^{e_1}, 0, - 2 \, p^{e_1})$}

\begin{lem}
Let $\varphi^* : \G_a \rtimes \G_m \to \SL(3, k)$ be a representation with the form 
\begin{fleqn}[36pt] 
\begin{align*}
\varphi^*
\left(
\begin{array}{c c}
 a & b \\
 0 & d 
\end{array}
\right) 
& =
\left(
\begin{array}{c c c}
 a^{2 \, p^{e_1}} & 0  & 0 \\
 0 & 1 &  b^{p^{e_1}} \, d^{p^{e_1}}  \\
 0 & 0 & d^{ 2 \, p^{e_1} } 
\end{array} 
\right) 
\qquad 
(\, e_1 \geq 0  \, ) . 
\end{align*}
\end{fleqn} 
Then $\varphi^*$ is not fundamental. 
\end{lem}

\Proof 
We have 
\begin{fleqn}[36pt] 
\begin{align*}
u_{\psi^*}(t) 
& =
\left(
\begin{array}{c c c}
 1 & 0 & 0  \\
 0 & 1 & t^{p^{e_1}} \\
 0 & 0 & 1 
\end{array} 
\right) . 
\end{align*}
\end{fleqn} 
Suppose, to the contrary, that $\varphi^*$ is fundamental. 
We can express $u^-_{\psi^*} : \G_a \to \SL(3, k)$ as  
\begin{fleqn}[36pt] 
\begin{align*}
 u^-_{\psi^*}(s)
  = 
\left(
\begin{array}{c c c}
 1 & 0 & 0 \\
 v_{2, \, 1}(s) & 1 & 0 \\
 v_{3, \, 1}(s) & v_{3, \, 2}(s) & 1 
\end{array}
\right) 
\qquad 
(\  v_{2, \, 1}(S), v_{3, \, 1}(S) ,  v_{3, \, 2}(S)  \in k[S]  \,) . 
\end{align*}
\end{fleqn} 
Using Lemma 3.3 (2), we have, for all $\gamma \in k \backslash \{ - 1 \}$, 
\begin{fleqn}[36pt] 
\begin{align*}
& 
\left(
\begin{array}{c c c}
 1 & 0 & 0 \\
 v_{2, \, 1}(\gamma) + v_{3, \, 1}(\gamma)  & 1 + v_{3, \, 2}(\gamma) & 1 \\
 v_{3, \, 1}(\gamma) & v_{3, \, 2}(\gamma) & 1
\end{array}
\right) \\
 & \qquad = 
\left(
\begin{array}{c c c}
 (1 + \gamma)^{2 \, p^{e_1}} & 0 & 0 \\  [6pt] 
 v_{2, \, 1} \left( \frac{\gamma}{1 + \gamma} \right) \, (1 + \gamma)^{2 \, p^{e_1}}
 &  1 
 &   \frac{1}{(1 + \gamma)^{p^{e_1}} }\\ [6pt] 
 v_{3, \, 1} \left( \frac{\gamma}{1 + \gamma} \right) \, (1 + \gamma)^{2 \, p^{e_1}}
 &  v_{3, \, 2}  \left( \frac{\gamma}{1 + \gamma} \right) 
 & v_{3, \, 2}  \left( \frac{\gamma}{1 + \gamma} \right)  \, \frac{1}{(1 + \gamma)^{p^{e_1}} } + \frac{1}{(1 + \gamma)^{p^{2 \, e_1}} }
\end{array}
\right) . 
\end{align*}
\end{fleqn} 
Comparing the $(1, 1)$-th entries of both sides of the above equality, we have $1 = (1 + \gamma)^{2 \, p^{e_1}}$ 
for all $\gamma \in k \backslash \{ - 1\}$. 
Considering an element $\gamma$ of $k$ such that $\gamma \notin k \backslash \{-1 , \, 0 \}$ and $\gamma^{p^{e_1}} \ne -2$, 
we have a contradiction. 
\QED

\subsubsection{(VII)$^*$: $p = 2$ and $(\ell_1, \ell_2, \ell_3) = (2 \, p^{e_1}, 0, - 2 \, p^{e_1})$}

\begin{lem}
If $p = 2$, we let $\varphi^* : \G_a \rtimes \G_m \to \SL(3, k)$ be a representation with the form 
\begin{fleqn}[36pt] 
\begin{align*}
\varphi^*
\left(
\begin{array}{c c}
 a & b \\
 0 & d 
\end{array}
\right) 
  = 
\left(
\begin{array}{c c c}
 a^{2 \, p^{e_1}}  &  0  &  b^{2 \, p^{e_1}}   \\
 0  &  1  &  b^{p^{e_1}} \, d^{p^{e_1}}  \\
 0 & 0 & d^{2 \, p^{e_1}} 
\end{array}  
\right) 
\qquad 
(\, e_1 \geq 0  \, ) . 
\end{align*}
\end{fleqn} 
Then the following assertions {\rm (1)}, {\rm (2)}, {\rm (3)} hold true: 
\begin{enumerate}[label = {\rm (\arabic*)}]
\item $\varphi^*$ is fundamental. 

\item Let $\psi^* : \SL(2, k) \to \SL(3, k)$ be a representation of $\SL(2, k)$ so that $\varphi^* = \psi^* \circ \iota$. 
Then we have 
\begin{fleqn}[36pt] 
\begin{align*}
 u_{\psi^*}^-(s) 
 = 
\left(
\begin{array}{c c c}
 1 & 0 & 0 \\
  s^{p^{e_1}}  & 1 & 0 \\
 s^{2 \, p^{e_1}} & 0 & 1 
\end{array}
\right) .
\end{align*}
\end{fleqn}

\item There exists a unique representation $\psi^* : \SL(2, k) \to \SL(3, k)$ such that 
$\varphi^* = \psi^* \circ \iota$. 
\end{enumerate} 
\end{lem}

\Proof 
(1) 
Consider the representation $\psi^* : \SL(2, k) \to \SL(3, k)$ defined by 
\begin{fleqn}[36pt] 
\begin{align*} 
\psi^*
\left(
\begin{array}{c c}
 a & b \\
 c & d 
\end{array}
\right) 
  = 
\left(
\begin{array}{c c c}
 a^{2 \, p^{e_1}}  &  0  &   b^{2 \, p^{e_1}}   \\
 a^{p^{e_1}} \, c^{ p^{e_1}}   &  1  &  b^{p^{e_1}} \, d^{ p^{e_1}}  \\
 c^{ 2 \, p^{e_1}}  & 0 & d^{ 2 \, p^{e_1}} 
\end{array}  
\right) . 
\end{align*}
\end{fleqn}

(2) 
We have 
\begin{fleqn}[36pt] 
\begin{align*}
 u_{\psi^*}(t) 
  = 
\left(
\begin{array}{c c c}
 1 & 0 & t^{p^{e_2}}  \\
 0 & 1 & 0 \\
 0 & 0 & 1
\end{array}
\right) . 
\end{align*}
\end{fleqn} 
We can express $u^-_{\psi^*} : \G_a \to \SL(3, k)$ as  
\begin{fleqn}[36pt] 
\begin{align*}
 u^-_{\psi^*}(s)
  = 
\left(
\begin{array}{c c c}
 1 & 0 & 0 \\
 v_{2, \, 1}(s) & 1 & 0 \\
 v_{3, \, 1}(s) & v_{3, \, 2}(s) & 1 
\end{array}
\right) 
\qquad 
(\  v_{2, \, 1}(S), v_{3, \, 1}(S) ,  v_{3, \, 2}(S)  \in k[S]  \,) . 
\end{align*}
\end{fleqn} 
Using Lemma 3.3 (2), we have, for all $\gamma \in k \backslash \{ - 1 \}$, 
\begin{align*}
& 
\left(
\begin{array}{c c c}
 1 +  v_{3,\, 1}(\gamma) & v_{3, \, 2}(\gamma) & 1 \\
 v_{2, \, 1}(\gamma) +  v_{3,\, 1}(\gamma) & 1 + v_{3, \, 2}(\gamma) & 1 \\
 v_{3, \, 1}(\gamma) & v_{3, \, 2}(\gamma) & 1 
\end{array}
\right) \\ 
 & \qquad = 
\left(
\begin{array}{c c c}
 (1 + \gamma)^{2 \, p^{e_1}} & 0 & 1 \\  [6pt] 
 v_{2, \, 1} \left( \frac{\gamma}{1 + \gamma} \right) \, (1 + \gamma)^{2 \, p^{e_1}}
 &   1 
 & v_{2, \, 1} \left( \frac{\gamma}{1 + \gamma} \right)  + \frac{1}{(1 + \gamma)^{p^{e_1}} }  \\ [6pt] 
 v_{3, \, 1} \left( \frac{\gamma}{1 + \gamma} \right) \, (1 + \gamma)^{2 \, p^{e_1}}
 &   v_{3, \, 2}  \left( \frac{\gamma}{1 + \gamma} \right) 
 & v_{3, \, 1} \left( \frac{\gamma}{1 + \gamma} \right) 
+ v_{3, \, 2}  \left( \frac{\gamma}{1 + \gamma} \right)  \, \frac{1}{(1 + \gamma)^{p^{e_1}} }
 + \frac{1}{(1 + \gamma)^{2 \, p^{e_1}} }
\end{array}
\right) . 
\tag{$\ast$}
\end{align*}
Comparing the $(1, 1)$-th entries of both sides of the above equality $(\ast)$, 
we have $1 + v_{3, \, 1}(\gamma) = (1 + \gamma)^{2 \, p^{e_1}}$ for all 
$\gamma \in k \backslash \{ - 1 \}$, which implies $v_{3, \, 1}(S) = S^{2 \, p^{e_1}}$. 
Comparing the $(1, 2)$-th entries of both sides of the above equality $(\ast)$, we have $ v_{3,\, 2}( \gamma) = 0 $ for all 
$\gamma \in k \backslash \{ -1 \}$, which implies $ v_{3,\, 2}(S) = 0$. 
Comparing the $(2, 3)$-th entries of both sides of the equality $(\ast)$, 
we have $1 = v_{2, \, 1}(\frac{\gamma}{1 + \gamma}) +  \frac{1}{ (1 + \gamma)^{p^{e_1}}}$ for all $\gamma \in k \backslash \{ - 1 \}$, 
which implies $v_{2, \, 1}(S) = S^{p^{e_1} }$ (since $p = 2$).

(3) 
The proof is straightforward.  
\QED

\subsubsection{(XII)$^*$}

\begin{lem}
Let $\varphi^* : \G_a \rtimes \G_m \to \SL(3, k)$ be a representation with the form 
\begin{fleqn}[36pt] 
\begin{align*}
\varphi^*
\left(
\begin{array}{c c}
 a & b \\
 0 & d 
\end{array}
\right) 
& =
\left(
\begin{array}{c c c}
  1 & 0 & 0 \\
 0 & 1& 0 \\
 0 & 0 & 1 
\end{array}  
\right) . 
\end{align*}
\end{fleqn} 
Then the following assertions {\rm (1)}, {\rm (2)}, {\rm (3)} hold true: 
\begin{enumerate}[label = {\rm (\arabic*)}]
\item $\varphi^*$ is fundamental. 

\item Let $\psi^* : \SL(2, k) \to \SL(3, k)$ be a representation of $\SL(2, k)$ so that $\varphi^* = \psi^* \circ \iota$. 
Then we have 
\begin{fleqn}[36pt] 
\begin{align*}
 u_{\psi^*}^-(s) 
 = 
\left(
\begin{array}{c c c}
 1 & 0 & 0 \\
 0 & 1 & 0 \\
 0 & 0 & 1 
\end{array}
\right) .
\end{align*}
\end{fleqn}

\item There exists a unique representation $\psi^* : \SL(2, k) \to \SL(3, k)$ such that 
$\varphi^* = \psi^* \circ \iota$. 
\end{enumerate} 
\end{lem}

\Proof 
(1) Consider the representation $\psi^* : \SL(2, k) \to \SL(3, k)$ defined by 
\begin{fleqn}[36pt] 
\begin{align*} 
\psi^*
\left(
\begin{array}{c c}
 a & b \\
 c & d 
\end{array}
\right) 
  = 
\left(
\begin{array}{c c c}
  1 & 0 & 0 \\
 0 & 1& 0 \\
 0 & 0 & 1 
\end{array}  
\right) . 
\end{align*}
\end{fleqn}

(2) Clearly, $u_{\varphi^*}$ is trivial. See Lemma 1.9.

(3) 
The proof is straightforward.  
\QED

\subsection{A classification of fundamental representations of $\G_a \rtimes \G_m$ into $\SL(3, k)$}

Based on the above Subsection 3.1, 
we can define five fundamental representations $\varphi^\sharp : \G_a \rtimes \G_m \to \SL(3, k)$, as follows: 

\begin{enumerate}[label = {\rm (\Roman*)$^\sharp$}, leftmargin = *]  

\item In the case where $p \geq 3$, 
\begin{fleqn}[36pt] 
\begin{align*} 
\varphi^\sharp
\left(
\begin{array}{c c}
 a & b \\
 0 & d 
\end{array}
\right) 
 & =
\left(
\begin{array}{c c c}
 a^{2 \, p^{e_1}}  &  a^{p^{e_1}} \, b^{p^{e_1}}  & \frac{1}{2} \,  b^{2 \, p^{e_1}} \\
 0  &  1  &   b^{p^{e_1}} \,  d^{p^{e_1}} \\
 0 & 0 & d^{2 \, p^{e_1}} 
\end{array}  
\right) 
\qquad 
(\, e_1 \geq 0 \, ) . 
\end{align*}
\end{fleqn}

\item[\rm (IV)$^\sharp$]  In the case where $p \geq 2$, 
\begin{fleqn}[36pt] 
\begin{align*}
\varphi^\sharp
\left(
\begin{array}{c c}
 a & b \\
 0 & d 
\end{array}
\right) 
& =
\left(
\begin{array}{c c c}
 a^{p^{e_2}} & 0 &   b^{p^{e_2}}  \\
 0 & 1 & 0 \\
 0 & 0 &  d^{p^{e_2}} 
\end{array} 
\right) 
\qquad  ( \, e_2 \geq 0 \, ) . 
\end{align*}
\end{fleqn}

\item[\rm (V)$^\sharp$] 
In the case where $p = 2$, 
\begin{fleqn}[36pt] 
\begin{align*}
\varphi^\sharp 
\left(
\begin{array}{c c}
 a & b \\
 0 & d 
\end{array}
\right) 
  = 
\left(
\begin{array}{c c c}
 a^{2 \, p^{e_1}}  &  a^{p^{e_1}} \, b^{p^{e_1}}   &  b^{2 \, p^{e_1}}\\
 0  &  1  &  0  \\
 0 & 0 & d^{ 2 \, p^{e_1}} 
\end{array}  
\right) 
 \qquad  ( \, e_1 \geq 0 \, ) . 
\end{align*}
\end{fleqn}

\item[\rm (VII)$^\sharp$]   
In the case where $p = 2$, 
\begin{fleqn}[36pt] 
\begin{align*}
\varphi^\sharp
\left(
\begin{array}{c c}
 a & b \\
 0 & d 
\end{array}
\right) 
  = 
\left(
\begin{array}{c c c}
 a^{2 \, p^{e_1}}  &  0  &  b^{2 \, p^{e_1}}   \\
 0  &  1  &  b^{p^{e_1}} \, d^{p^{e_1}}  \\
 0 & 0 & d^{2 \, p^{e_1}} 
\end{array}  
\right) 
\qquad 
(\, e_1 \geq 0  \, ) . 
\end{align*}
\end{fleqn}

\item[\rm (XII)$^\sharp$]  
In the case where $p \geq 2$, 
\begin{fleqn}[36pt] 
\begin{align*}
\varphi^\sharp
\left(
\begin{array}{c c}
 a & b \\
 0 & d 
\end{array}
\right) 
  = 
\left(
\begin{array}{c c c}
  1 & 0 & 0 \\
 0 & 1& 0 \\
 0 & 0 & 1 
\end{array}  
\right) . 
\end{align*}
\end{fleqn} 

\end{enumerate}

Let
\begin{align*}
\Lambda^\sharp & := 
\bigl\{ \, 
 {\rm (I)^\sharp}, \;  {\rm (IV)^\sharp}, \;  {\rm (V)^\sharp}, \;  {\rm (VII)^\sharp}, \;  {\rm (XII)^\sharp} 
\, \bigr\}.
\end{align*}

For any $\lambda \in \Lambda^\sharp$, we can define a set $R^\sharp(\lambda)$ as 
\begin{align*}
 R^\sharp(\lambda) 
 := \{ \, \varphi : \G_a \rtimes \G_m \to \SL(3, k) \mid \text{$\varphi$ is a fundamental representation with the form $\lambda$} \, \} . 
\end{align*}

\begin{thm}
The following assertions {\rm (1)} and {\rm (2)} hold true: 
\begin{enumerate}[label = {\rm (\arabic*)}]
\item Any fundamental representation $\varphi : \G_a \rtimes \G_m \to \SL(3, k)$ is equivalent to 
a fundamental representation $\varphi^\sharp : \G_a \rtimes \G_m \to \SL(3, k)$ with one of the forms 
$\lambda$ of $\Lambda^\sharp$. 

\item Let $\phi \in R^\sharp(\lambda)$ and $\psi \in R^\sharp(\mu)$, where $\lambda, \mu \in \Lambda^\sharp$. 
Write $h_\phi(z) = (z^{\ell_1},  z^{\ell_2},  z^{\ell_3})$ and $h_{\psi}(z) = (z^{m_1},  z^{m_2},  z^{m_3})$.   
\begin{enumerate}[label = {\rm (2.\arabic*)}]
\item 
If $\lambda \ne \mu$, then $\phi$ and $\psi$ are not equivalent. 

\item If $\lambda = \mu$, then the following conditions {\rm (i)}, {\rm (ii)}, {\rm (iii)}, {\rm (iv)} are equivalent: 
\begin{enumerate}[label = {\rm (\roman*)}]
\item $\phi$ and $\psi$ are equivalent. 

\item $h_\phi$ and $h_{\psi}$ are equivalent. 

\item $(\ell_1, \, \ell_2, \, \ell_3) = (m_1, \, m_2, \, m_3)$. 

\item $\phi = \psi$. 
\end{enumerate}
 
\end{enumerate} 

\end{enumerate} 
\end{thm}

\Proof 
(1) See Lemmas 3.2, 3.4, 3.5, 3.6, 3.7, 3.8, 3.9, 3.10, 3.11.

(2) See assertion (2) of Theorem 2.3. 
\QED

\begin{cor}
There exists a one-to-one correspondence between the set of all fundamental representations of $\G_a \rtimes \G_m$ 
into $\SL(3, k)$ up to equivalence and the set $\coprod_{\lambda \in \Lambda^\sharp} R^\sharp(\lambda)$, i.e., 
\[
\{\, \varphi : \G_a \rtimes \G_m \to \SL(3, k) \mid \text{$\varphi$ is a fundamental representation}\, \} / \sim 
 \; \cong  \: 
 \coprod_{\lambda \in \Lambda^\sharp} R^\sharp(\lambda) . 
\]
\end{cor}

\Proof The proof is straightforward (see the proof of Corollary 2.4 and use the above Theorem 3.12). 
\QED

\section{Fundamental representations of $\G_a$ into $\SL(3, k)$}

A representation $u :  \G_a \to \SL(n, k)$ of $\G_a \rtimes \G_m$ 
is said to be {\it fundamental} if there exists a representation $\psi : \SL(2, k) \to \SL(n, k)$ such that 
the following diagram commutes: 
\[
\xymatrix@R=36pt@C=36pt@M=8pt{
  \G_a \ar[r]^(.43)u \ar@{^(->}[d]_{\iota^+} &  \SL(3, k) \\
 \SL(2, k) \ar[ru]_\psi 
}
\]
where $\hspace{-6pt} \xymatrix@C=12pt@M=6pt{\iota^+ : \G_a \ar@{^(->}[r] & B}\hspace{-6pt} $ is the inclusion map defined by 
\[
\iota^+(t)
 := 
\left(
\begin{array}{c c}
 1 & t \\
 0 & 1
\end{array}
\right) . 
\]

\begin{thm}
Let $u : \G_a \to \SL(3, k)$ be a morphism of affine $k$-varieties. 
Then the following assertions {\rm (1)} and {\rm (2)} hold true: 
\begin{enumerate}[label = {\rm (\arabic*)}]
\item If $p = 2$, then $u : \G_a \to \SL(3, k)$ is a fundamental representation if and only if 
$u$ is equivalent to a fundamental representation $u^\sharp : \G_a \to \SL(3, k)$ with one of the following forms {\rm (1.1)}, 
{\rm (1.2)}, {\rm (1.3)}, {\rm (1.4)}: 
\begin{enumerate}[label = {\rm (1.\arabic*)}, leftmargin = *]
\item 
$
u^\sharp(t)
 = 
\left( 
\begin{array}{c c c}
 1 & 0 & t^{p^e} \\
 0 & 1 & 0 \\
 0 & 0 & 1
\end{array}
\right)
$
\qquad 
$(\, e \geq 0 \,)$.

\item 
$
u^\sharp(t) 
 =
\left( 
\begin{array}{c c c}
 1 & t^{p^e} & t^{2 \, p^e} \\
 0 & 1 & 0 \\
 0 & 0 & 1
\end{array}
\right)
$
\qquad 
$(\, e \geq 0 \,) $.

\item 
$
u^\sharp(t) 
 =
\left( 
\begin{array}{c c c}
 1 & 0 & t^{2 \, p^e} \\
 0 & 1 & t^{p^e} \\
 0 & 0 & 1
\end{array}
\right)
$
\qquad 
$(\, e \geq 0 \,) $.

\item 
$
u^\sharp(t) 
 =
\left( 
\begin{array}{c c c}
 1 & 0 & 0\\
 0 & 1 & 0 \\
 0 & 0 & 1
\end{array}
\right) 
$. 
\end{enumerate}

\item If $p \geq 3$, then $u : \G_a \to \SL(3, k)$ is a fundamental representation if and only if 
$u$ is equivalent to a fundamental representation $u^\sharp : \G_a \to \SL(3, k)$ with one of the following forms 
{\rm (2.1)}, {\rm (2.2)}, {\rm (2.3)}: 
\begin{enumerate}[label = {\rm (2.\arabic*)}, leftmargin = *]
\item 
$
u^\sharp(t) 
 =
\left( 
\begin{array}{c c c}
 1 & t^{p^e} & \frac{1}{2} \, t^{2 \, p^e} \\
 0 & 1 & t^{p^e} \\
 0 & 0 & 1
\end{array}
\right) 
$
\qquad 
$(\, e \geq 0 \, ) $.

\item 
$
u^\sharp(t) 
 =
\left( 
\begin{array}{c c c}
 1 & 0 & t^{p^e} \\
 0 & 1 & 0 \\
 0 & 0 & 1
\end{array}
\right)
$
\qquad 
$(\, e \geq 0 \,)$.

\item 
$
u^\sharp(t) 
 =
\left( 
\begin{array}{c c c}
 1 & 0 & 0\\
 0 & 1 & 0 \\
 0 & 0 & 1
\end{array}
\right) 
$. 
\end{enumerate}

\end{enumerate}   
\end{thm}

\Proof Assume that $u : \G_a \to \SL(3, k)$ is a fundamental representation of $\G_a$. 
So, there exists a representation $\psi : \SL(2, k) \to \SL(3, k)$ such that $\psi \circ \iota^+ = u$. 
Let $\varphi : \G_a \rtimes \G_m \to \SL(3, k)$ be the representation defined by $\varphi := \psi \circ \iota$, 
where $\iota : \G_a \rtimes \G_m \to \SL(2, k)$ is the injective homomorphism. 
Clearly, $\varphi$ is a fundamental representation. 
We can apply Theorem 3.12 to the $\varphi$.

We first prove assertion (1).  
If $p = 2$, there exists a representation $\varphi^\sharp : \G_a \rtimes \G_m \to \SL(3, k)$ such that 
$\varphi^\sharp$ is equivalent to $\varphi$ and $\varphi^\sharp$  
has one of the forms (IV)$^\sharp$, (V)$^\sharp$, (VII)$^\sharp$, (XII)$^\sharp$. 
Thus $u_{\varphi^\sharp} : \G_a \to \SL(3, k)$ has one of the forms (1.1), (1.2), (1.3), (1.4).

We next prove assertion (2).  
If $p \geq 3$, there exists a representation $\varphi^\sharp : \G_a \rtimes \G_m \to \SL(3, k)$ such that 
$\varphi^\sharp$ is equivalent to $\varphi$ and $\varphi^\sharp$  
has one of the forms (I)$^\sharp$, (IV)$^\sharp$, (XII)$^\sharp$. 
Thus $u_{\varphi^\sharp} : \G_a \to \SL(3, k)$ has one of the forms (2.1), (2.2), (2.3).  
\QED

\begin{cor}
The following assertions {\rm (1)} and {\rm (2)} hold true: 
\begin{enumerate}[label = {\rm (\arabic*)}]
\item If $p = 2$, then 
there exists a one-to-one correspondence between the set of all fundamental representations 
$u : \G_a \to \SL(3, k)$ up to equivalence and 
the set $\Z_{\geq 0} \, \amalg \, \Z_{\geq 0} \, \amalg \, \Z_{\geq 0} \, \amalg \, \{ I_3 \}$, 
i.e., 
\[
 \{ \, u : \G_a \to \SL(3, k) \mid \text{$u$ is fundamental} \, \} / \sim 
 \; \cong \; 
 \Z_{\geq 0} \, \amalg \, \Z_{\geq 0} \, \amalg \, \Z_{\geq 0} \, \amalg \,  \{\, I_3 \,\} . 
\]

\item If $p \geq 3$, then 
there exists a one-to-one correspondence between the set of all fundamental representations 
$u : \G_a \to \SL(3, k)$ up to equivalence and the set $\Z_{\geq 0} \, \amalg \, \Z_{\geq 0} \, \amalg \, \{ I_3 \}$, 
i.e., 
\[
 \{ \, u : \G_a \to \SL(3, k) \mid \text{$u$ is fundamental} \, \} / \sim 
 \; \cong \; 
 \Z_{\geq 0} \, \amalg \, \Z_{\geq 0} \, \amalg \, \{ \, I_3 \, \} . 
\]
\end{enumerate}
\end{cor}

\Proof 
The proofs of assertions (1) and (2) are straightforward. 
\QED

We also have the following corollary, which can yield new non-fundamental representations of $\G_a$ into $\SL(3, k)$.

\begin{cor}
Let $u : \G_a \to \SL(3, k)$ be a morphism of affine $k$-varieties with the form  
\begin{align*}
u(t) 
 = 
\left(
\begin{array}{c c c}
 1 & a_{1, 2}(t) & a_{1, 3}(t) \\
 0 & 1 & a_{2, 3}(t) \\
 0 & 0 & 1 
\end{array}
\right) 
\qquad 
(\, a_{1, 2}(T), a_{1, 3}(T), a_{2, 3}(T) \in k[T]  \,) . 
\end{align*}
Then the following assertions {\rm (1)} and {\rm (2)} hold true: 
\begin{enumerate}[label = {\rm (\arabic*)}]
\item If $p = 2$, then $u$ is a fundamental representation of $\G_a$ if and only if  $(a_{1, 2}(T),  \, a_{1, 3}(T) , \, a_{2, 3}(T) )$ 
satisfies one of the following conditions {\rm (1.1.a)}, {\rm (1.1.b)}, {\rm (1.2)}, {\rm (1.3)}, {\rm (1.4)}: 
\begin{enumerate}[label = {\rm (1.\arabic*)}, start = 2]

\item[\rm (1.1.a)] 
$(a_{1, 2}(T),  \, a_{1, 3}(T) , \, a_{2, 3}(T) ) = (\lambda \, T^{p^e}, \, \nu \, T^{p^e}, \, 0)$ 
for some $\lambda, \nu \in k$ satisfying $(\lambda, \nu ) \ne (0, 0)$ and for some $e \geq 0$.

\item[\rm (1.1.b)] 
$(a_{1, 2}(T),  \, a_{1, 3}(T) , \, a_{2, 3}(T) ) = (0, \, \nu \, T^{p^e}, \, \lambda \, T^{p^e})$ 
for some $\lambda, \nu \in k$ satisfying $ (\lambda, \nu ) \ne (0, 0)$ and for some $e \geq 0$.

\item 
$(a_{1, 2}(T),  \, a_{1, 3}(T) , \, a_{2, 3}(T) ) = (\lambda \, T^{p^e} + \mu \, T^{2 \, p^e}, \, \nu \, T^{p^e} + \xi \, T^{2 \, p^e}, \, 0)$ 
for some $\lambda, \mu, \nu, \xi \in k$ satisfying $\lambda \xi - \mu \nu \ne 0$ and for some $e \geq 0 $.

\item 
$(a_{1, 2}(T),  \, a_{1, 3}(T) , \, a_{2, 3}(T) ) = (0, \, \nu \, T^{p^e} + \xi \, T^{2 \, p^e}, \, \lambda \, T^{p^e} + \mu \, T^{2 \, p^e} )$
for some $\lambda, \mu, \nu, \xi \in k$ satisfying $\lambda \xi - \mu \nu \ne 0$ and for some $e \geq 0$.

\item $(a_{1, 2}(T),  \, a_{1, 3}(T) , \, a_{2, 3}(T) ) =  (0, \, 0, \, 0)$. 
\end{enumerate}

\item If $p \geq 3$, then $u$ is a fundamental representation of $\G_a$ if and only if $(a_{1, 2}(T),  \, a_{1, 3}(T) , \, a_{2, 3}(T) )$ 
satisfies one of the following conditions {\rm (2.1)}, {\rm (2.2.a)}, {\rm (2.2.b)}, {\rm (2.3)}: 
\begin{enumerate}[label = {\rm (2.\arabic*)}]
\item 
$(a_{1, 2}(T),  \, a_{1, 3}(T) , \, a_{2, 3}(T) ) = (c \, T^{p^e}, \, \frac{1}{2} \, \lambda \, c^2 \, T^{2 \, p^e}, \, \lambda \, c \, T^{p^e})$ 
for some $ c, \lambda \in k \backslash \{ 0 \}$ and for some $e \geq 0$.

\item[\rm (2.2.a)] 
$(a_{1, 2}(T),  \, a_{1, 3}(T) , \, a_{2, 3}(T) ) = (\lambda \, T^{p^e}, \, \nu \, T^{p^e}, \, 0)$ 
for some $\lambda, \nu \in k$ satisfying $(\lambda, \nu ) \ne (0, 0)$ and for some $e \geq 0$.

\item[\rm (2.2.b)] 
$(a_{1, 2}(T),  \, a_{1, 3}(T) , \, a_{2, 3}(T) ) = (0, \, \nu \, T^{p^e}, \, \lambda \, T^{p^e})$ 
for some $\lambda, \nu \in k$ satisfying $ (\lambda, \nu ) \ne (0, 0)$ and for some $e \geq 0$.

\item[\rm (2.3)] $(a_{1, 2}(T),  \, a_{1, 3}(T) , \, a_{2, 3}(T) ) =  (0, \, 0, \, 0)$. 
\end{enumerate} 

\end{enumerate} 
\end{cor}

\Proof 
(1) Assume $p = 2$ and $u$ is fundamental. 
There exists a regular matrix $P$ of $\GL(3, k)$ such that 
$u^\sharp(t) := P^{-1}  u(t) P$ has one of the forms (1.1), (1.2), (1.3), (1.4) given in assertion (1) of Theorem 4.1. 
Clearly, 
\begin{align*}
 (u(t) - I_3 ) \, P = P \, (u^\sharp(t) - I_3) . 
\tag{$\ast$}
\end{align*}
Write $P = (p_{i, j})_{1 \leq i, j \leq 3}$. 

In the case where $u^\sharp$ has the form (1.1) of Theorem 4.1, we obtain from $(\ast)$ that  
\[
\left(
\begin{array}{c c c}
 0 & a_{1, 2}(t) & a_{1, 3}(t) \\
 0 & 0 & a_{2, 3}(t) \\
 0 & 0 & 0 
\end{array}
\right) 
\left(
\begin{array}{c c c}
 p_{1, 1} & p_{1, 2} & p_{1, 3} \\
 p_{2, 1} & p_{2, 2} & p_{2, 3} \\
 p_{3, 1} & p_{3, 2} & p_{3, 3} 
\end{array}
\right) 
 = 
\left(
\begin{array}{c c c}
 p_{1, 1} & p_{1, 2} & p_{1, 3} \\
 p_{2, 1} & p_{2, 2} & p_{2, 3} \\
 p_{3, 1} & p_{3, 2} & p_{3, 3} 
\end{array}
\right) 
\left( 
\begin{array}{c c c}
 0 & 0 & t^{p^e} \\
 0 & 0 & 0 \\
 0 & 0 & 0
\end{array}
\right) . 
\]
We argue by separating the following two cases:  
\begin{enumerate}[label = {(\alph*)}]
\item $a_{2, 3}(T) = 0$. 

\item $a_{2, 3}(T) \ne 0$
\end{enumerate} 

In the case (a), 
\begin{align*}
\left(
\begin{array}{c c c}
 p_{2, 1} \, a_{1, 2}(t) + p_{3, 1} \, a_{1, 3}(t) 
 & p_{2, 2} \,  a_{1, 2}(t) + p_{3, 2}  \, a_{1, 3}(t) 
 & p_{2, 3} \,  a_{1, 2}(t) + p_{3, 3} \, a_{1, 3}(t)  \\
 0 & 0 & 0 \\
 0 & 0 & 0 
\end{array}
\right) \\
 = 
\left( 
\begin{array}{c c c}
 0 & 0 & p_{1, 1} \,  t^{p^e} \\
 0 & 0 & p_{2, 1} \,  t^{p^e} \\
 0 & 0 & p_{3, 1} \,  t^{p^e}
\end{array}
\right) , 
\end{align*}
which implies $p_{2, 1} = p_{3, 1} = 0$. Since $P$ is regular, we have  
\[
\left(
\begin{array}{c c}
 a_{1, 2}(t) & a_{1, 3}(t)
\end{array}
\right)
=
\left(
\begin{array}{c c}
 0 & p_{1, 1} \, t^{p^e} 
\end{array}
\right)
\left(
\begin{array}{c c}
 p_{2, 2} & p_{2, 3} \\
 p_{3, 2} & p_{3, 3} 
\end{array}
\right)^{-1} . 
\]
Thus $(a_{1, 2}(T), a_{1, 3}(T), a_{2, 3}(T))$ satisfies the condition (1.1.a).

In the case (b), we have $a_{1, 2}(T) = 0$ since $p = 2$ (see Lemma 1.4). 
Thereby  
\[
\left(
\begin{array}{c c c}
 p_{3, 1} \, a_{1, 3}(t) & p_{3, 2} \, a_{1, 3}(t) & p_{3, 3} \, a_{1, 3}(t)  \\
 p_{3, 1} \, a_{2, 3}(t) & p_{3, 2} \, a_{2, 3}(t) & p_{3, 3} \, a_{2, 3}(t)  \\
 0 & 0 & 0
\end{array}
\right) 
 = 
\left( 
\begin{array}{c c c}
 0 & 0 & p_{1, 1} \,t^{p^e} \\
 0 & 0 & p_{2, 1} \, t^{p^e} \\
 0 & 0 & p_{3, 1} \, t^{p^e} 
\end{array}
\right) . 
\]
Since $a_{2, 3}(T) \ne 0$, we have $p_{3, 1} = p_{3, 2} = 0$. 
Since $P$ is regular, we have $(p_{1, 1}, p_{2, 1}) \ne (0, 0)$ and $p_{3, 3} \ne 0$. 
Thus $(a_{1, 2}(T), a_{1, 3}(T), a_{2, 3}(T))$ satisfies the condition (1.1.b).

In the case where $u^\sharp$ has the form (1.2) of Theorem 4.1, we obtain from $(\ast)$ that  
\[
\left(
\begin{array}{c c c}
 0 & a_{1, 2}(t) & a_{1, 3}(t) \\
 0 & 0 & a_{2, 3}(t) \\
 0 & 0 & 0 
\end{array}
\right) 
\left(
\begin{array}{c c c}
 p_{1, 1} & p_{1, 2} & p_{1, 3} \\
 p_{2, 1} & p_{2, 2} & p_{2, 3} \\
 p_{3, 1} & p_{3, 2} & p_{3, 3} 
\end{array}
\right) 
 = 
\left(
\begin{array}{c c c}
 p_{1, 1} & p_{1, 2} & p_{1, 3} \\
 p_{2, 1} & p_{2, 2} & p_{2, 3} \\
 p_{3, 1} & p_{3, 2} & p_{3, 3} 
\end{array}
\right) 
\left( 
\begin{array}{c c c}
 0 & t^{p^e} & t^{2 \, p^e} \\
 0 & 0 & 0 \\
 0 & 0 & 0
\end{array}
\right) . 
\]
We already know from Lemmas 2.5 and 2.6 that $a_{2, 3}(T) = 0$. 
So, 
\begin{align*}
\left(
\begin{array}{c c c}
 p_{2, 1} \, a_{1, 2}(t) + p_{3, 1} \, a_{1, 3}(t) 
 & p_{2, 2} \, a_{1, 2}(t) + p_{3, 2} \, a_{1, 3}(t) 
 & p_{2, 3} \, a_{1, 2}(t) + p_{3, 3} \, a_{1, 3}(t)  \\
 0
 & 0 
 & 0 \\
0 & 0 & 0 
\end{array}
\right) \\
  = 
\left(
\begin{array}{c c c}
 0 & p_{1, 1} \, t^{p^e}  & p_{1, 1} \, t^{2 \, p^e} \\
 0 & p_{2, 1} \, t^{p^e} & p_{2, 1} \, t^{2 \, p^e}  \\
 0 & p_{3, 1} \, t^{p^e} & p_{3, 1} \, t^{2 \, p^e}  
\end{array}
\right)  , 
\end{align*}
which implies $p_{2, 1} = p_{3, 1} = 0$. 
Since $P$ is regular, we have $p_{1, 1} \ne 0$ and 
\[
\left(
\begin{array}{c c}
 a_{1, 2}(t) & a_{1, 3}(t)
\end{array}
\right) 
 =
\left(
\begin{array}{c c}
p_{1, 1} \, t^{p^e}  & p_{1, 1} \, t^{2 \, p^e} 
\end{array}
\right) 
\left(
\begin{array}{c c}
 p_{2, 2}  &  p_{2, 3}  \\
 p_{3, 2}  & p_{3, 3}
\end{array} 
\right)^{- 1} .  
\]
Thus $(a_{1, 2}(T), a_{1, 3}(T), a_{2, 3}(T))$ satisfies the condition (1.2).

In the case where $u^\sharp$ has the form (1.3) of Theorem 4.1, we obtain from $(\ast)$ that  
\[
\left(
\begin{array}{c c c}
 0 & a_{1, 2}(t) & a_{1, 3}(t) \\
 0 & 0 & a_{2, 3}(t) \\
 0 & 0 & 0 
\end{array}
\right) 
\left(
\begin{array}{c c c}
 p_{1, 1} & p_{1, 2} & p_{1, 3} \\
 p_{2, 1} & p_{2, 2} & p_{2, 3} \\
 p_{3, 1} & p_{3, 2} & p_{3, 3} 
\end{array}
\right) 
 = 
\left(
\begin{array}{c c c}
 p_{1, 1} & p_{1, 2} & p_{1, 3} \\
 p_{2, 1} & p_{2, 2} & p_{2, 3} \\
 p_{3, 1} & p_{3, 2} & p_{3, 3} 
\end{array}
\right) 
\left( 
\begin{array}{c c c}
 0 & 0 & t^{2 \, p^e} \\
 0 & 0 & t^{p^e}  \\
 0 & 0 & 0
\end{array}
\right) . 
\]
We already know from Lemmas 2.5 and 2.6 that $a_{1, 2}(T) = 0$. 
So, 
\[
\left(
\begin{array}{c c c}
 p_{3, 1} \, a_{1, 3}(t) 
 & p_{3, 2} \, a_{1, 3}(t) 
 & p_{3, 3} \, a_{1, 3}(t) \\
 p_{3, 1} \, a_{2, 3}(t) 
 & p_{3, 2} \, a_{2, 3}(t) 
 & p_{3, 3} \, a_{2, 3}(t) \\
 0 & 0 & 0 
\end{array}
\right) 
 = 
\left( 
\begin{array}{c c c}
 0 & 0 & p_{1, 1} \, t^{2 \, p^e} + p_{1, 2} \, t^{p^e}  \\
 0 & 0 & p_{2, 1} \, t^{2 \, p^e} + p_{2, 2} \, t^{p^e}  \\
 0 & 0 & p_{3, 1} \, t^{2 \, p^e} + p_{3, 2} \, t^{p^e} 
\end{array}
\right) , 
\]
which implies $p_{3, 1} = p_{3, 2} = 0$. Since $P$ is regular, we have $p_{3, 3} \ne 0$. So, 
\[
\left(
\begin{array}{c}
  a_{1, 3}(t) \\
  a_{2, 3}(t) 
\end{array} 
\right)
 = 
\frac{1}{p_{3, 3}} 
\left(
\begin{array}{c}
p_{1, 1} \, t^{2 \, p^e} + p_{1, 2} \, t^{p^e}  \\
p_{2, 1} \, t^{2 \, p^e} + p_{2, 2} \, t^{p^e}  
\end{array} 
\right) . 
\]
Thus $(a_{1, 2}(T), a_{1, 3}(T), a_{2, 3}(T))$ satisfies the condition (1.3).

In the case where $u^\sharp$ has the form (1.4) of Theorem 4.1, 
$(a_{1, 2}(T), a_{1, 3}(T), a_{2, 3}(T))$ clearly satisfies the condition (1.4).  
\medskip

Conversely assume that $(a_{1, 2}(T), a_{1, 3}(T), a_{2, 3}(T))$ satisfies one of the conditions (1.1.a), (1.1.b), (1.2), (1.3), (1.4). 
Using Theorem 4.1, we can prove that $u$ is fundamental. 
\\


(2) Assume $p \geq 3$ and $u$ is fundamental. 
There exists a regular matrix $P$ of $\GL(3, k)$ such that 
$u^\sharp(t) := P^{-1}  u(t) P$ has one of the forms (2.1), (2.2), (2.3) given in assertion (2) of Theorem 4.1. 
Clearly, 
\begin{align*}
 (u(t) - I_3 ) \, P = P \, (u^\sharp(t) - I_3) . 
\tag{$\ast\ast$}
\end{align*}
Write $P = (p_{i, j})_{1 \leq i, j \leq 3}$. 

In the case where $u^\sharp$ has the form (2.1) of Theorem 4.1, we obtain from $(\ast\ast)$ that  
\[
\left(
\begin{array}{c c c}
 0 & a_{1, 2}(t) & a_{1, 3}(t) \\
 0 & 0 & a_{2, 3}(t) \\
 0 & 0 & 0 
\end{array}
\right) 
\left(
\begin{array}{c c c}
 p_{1, 1} & p_{1, 2} & p_{1, 3} \\
 p_{2, 1} & p_{2, 2} & p_{2, 3} \\
 p_{3, 1} & p_{3, 2} & p_{3, 3} 
\end{array}
\right) 
 = 
\left(
\begin{array}{c c c}
 p_{1, 1} & p_{1, 2} & p_{1, 3} \\
 p_{2, 1} & p_{2, 2} & p_{2, 3} \\
 p_{3, 1} & p_{3, 2} & p_{3, 3} 
\end{array}
\right) 
\left( 
\begin{array}{c c c}
 0 & t^{p^e} & \frac{1}{2} \, t^{2 \, p^e} \\
 0 & 0 & t^{p^e} \\
 0 & 0 & 0
\end{array}
\right) . 
\]
Thereby 
\begin{align*}
\left(
\begin{array}{c c c}
 p_{2, 1} \, a_{1, 2}(t) + p_{3, 1} \, a_{1, 3}(t) 
 & p_{2, 2} \, a_{1, 2}(t) + p_{3, 2} \, a_{1, 3}(t) 
 & p_{2, 3} \, a_{1, 2}(t) + p_{3, 3} \, a_{1, 3}(t)  \\
 p_{3, 1} \, a_{2, 3}(t)
 & p_{3, 2} \, a_{2, 3}(t)
 & p_{3, 3} \, a_{2, 3}(t) \\
0 & 0 & 0 
\end{array}
\right) \\
 = 
\left( 
\begin{array}{c c c}
 0 & p_{1, 1} \, t^{p^e}  & \frac{1}{2} \, p_{1, 1} \, t^{2 \, p^e} + p_{1, 2} \, t^{p^e}  \\
 0 & p_{2, 1} \, t^{p^e}  & \frac{1}{2} \, p_{2, 1} \, t^{2 \, p^e} + p_{2, 2} \, t^{p^e} \\
 0 & p_{3, 1} \, t^{p^e}  & \frac{1}{2} \, p_{3, 1} \, t^{2 \, p^e} + p_{3, 2} \, t^{p^e}
\end{array}
\right) . 
\end{align*} 
We already know from Lemmas 2.5 and 2.6 that $a_{1, 2}(T) \ne 0$ and $a_{2, 3}(T) \ne 0$. 
Comparing the $(i, j)$-th entries ($(i, j) = (1, 1), (2, 1), (3, 3)$) of both sides of the above equality, 
we have $p_{2, 1} = p_{3, 1} = p_{3, 2} = 0$. 
Since $P$ is regular, we have $p_{1, 1} \ne 0$, $p_{2, 2} \ne 0$ and $p_{3, 3} \ne 0$. 
Comparing the $(i, j)$-th entries ($(i, j) = (1, 2), (2, 3), (1, 3)$) 
of both sides of the above equality, we know that $(a_{1, 2}(T),  \, a_{1, 3}(T) , \, a_{2, 3}(T) )$ satisfies the condition (2.1).

In the case where $u^\sharp$ has the form (2.2) of Theorem 4.1, 
we know that $(a_{1, 2}(T),  \, a_{1, 3}(T) , \, a_{2, 3}(T) )$ satisfies one of the conditions (2.2.a) and (2.2.b) 
(see the proof of the above assertion (1)).

In the case where $u^\sharp$ has the form (2.3) of Theorem 4.1, 
$(a_{1, 2}(T), a_{1, 3}(T), a_{2, 3}(T))$ clearly satisfies the condition (2.3).  
\medskip

Conversely assume that $(a_{1, 2}(T), a_{1, 3}(T), a_{2, 3}(T))$ satisfies one of the conditions 
(2.1), (2.2.a), (2.2.b), (2.3). Using Theorem 4.1, we can prove that $u$ is fundamental. 
\QED

\vspace{1cm}

\begin{flushright}
\begin{tabular}{l}
 Faculty of Education,\\
 Shizuoka University,\\
 836 Ohya, Suruga-ku,\\
 Shizuoka 422-8529, Japan\\
e-mail: tanimoto.ryuji@shizuoka.ac.jp
\end{tabular}
\end{flushright}

\end{document}